\documentclass[12pt, twoside]{article}
\usepackage[english]{babel}
\usepackage[utf8]{inputenc}
\usepackage[a4paper,top=2.5cm,bottom=2.5cm,left=1.8cm,right=1.8cm]{geometry}
\usepackage{indentfirst}
\usepackage{amsmath,amsfonts,amsthm,mathtools,amssymb,physics}
\usepackage{dsfont}
\usepackage{url}
\usepackage{tikz}
\usepackage{graphicx}
\usepackage{float}
\usepackage{enumitem}
\usepackage{xcolor}
\usepackage{hyperref}
\hypersetup{colorlinks, allcolors={black}}
\usepackage{cleveref}

\newtheorem{theorem}{Theorem}[section]
\newtheorem{corollary}{Corollary}[theorem]
\newtheorem{proposition}[theorem]{Proposition}
\newtheorem{lemma}[theorem]{Lemma}\theoremstyle{plain}
\newtheorem{introtheorem}{Theorem}

\theoremstyle{definition}
\newtheorem{definition}[theorem]{Definition}

\theoremstyle{remark}
\newtheorem{remark}[theorem]{Remark}
\theoremstyle{remark}
\newtheorem{example}[theorem]{Example}

\newcommand{\R}{\mathbb R}
\newcommand{\N}{\mathbb N}

\newcommand{\prob}{\mathcal P}

\newcommand{\Lip}{\operatorname{Lip}}
\newcommand{\M}{\mathcal{M}\bigl(\R^d\bigr)}
\newcommand{\probtwo}{\mathcal{P}_2\bigl(\R^d\bigr)}

\newcommand{\cone}{\mathfrak{C}}
\newcommand{\dc}{\mathrm{d}_\cone}
\newcommand{\dcwfr}{\mathrm{d}_{\cone,\frac{\pi}{2}}}

\newcommand{\y}{\mathfrak{y}}
\newcommand{\h}{\mathfrak{h}}
\newcommand{\vertex}{\mathfrak{o}}

\newcommand{\Geo}{\operatorname{Geo}}
\newcommand{\bfmu}{\boldsymbol{\mu}}
\newcommand{\bfy}{\boldsymbol{\mathfrak{y}}}
\newcommand{\bfx}{\boldsymbol{x}}
\newcommand{\bfr}{\boldsymbol{r}}
\newcommand{\bfm}{\boldsymbol{m}}

\newcommand{\bfomega}{\boldsymbol{\omega}}

\newcommand{\bfgamma}{\boldsymbol{\gamma}}

\newcommand{\oSC}{\operatorname{OptSC}}

\newcommand{\proj}[1]{\pi^{#1}}
\newcommand{\pf}{_{\#}}

\newcommand{\WFR}{\mathrm{WFR}}

\newcommand{\barlog}{\overline{\operatorname{Log}}}

\newcommand{\Log}{\operatorname{Log}}

\title{A Regular First-Order Calculus for the Wasserstein--Fisher--Rao Geometry}
\author{Enrico Sartor}
\date{}

\begin{document}

\maketitle

\begin{abstract}
We introduce a notion of intrinsic differentiability for functionals on the space of finite nonnegative measures equipped with the Wasserstein--Fisher--Rao (WFR) distance, also known as the Hellinger--Kantorovich (HK) distance. The resulting static first-order calculus accounts simultaneously for spatial transport and mass reweighting. We characterize regular barycentric logarithms through optimal semicouplings and measurable selections at the transport threshold, and prove their tangency at arbitrary finite source measures. The defining expansion extends to non-optimal semicouplings, with remainder of smaller order than the square root of their cost. This yields differentiation formulas for simultaneous transport and reweighting and uniqueness of the gradients. We also give a sufficient criterion in terms of regular flat derivatives. 
We characterize differentiability of one half the squared WFR distance by a common optimal Monge transport--growth representation with no target
mass created from the cone vertex. Equivalently, the regular barycentric
logarithm is unique and its norm equals the distance; the gradient is
minus this logarithm. Finally, we formulate the corresponding calculus on probability measures. On $\probtwo$, WFR differentiability implies intrinsic Wasserstein differentiability, with the transport component of the WFR gradient coinciding with the Wasserstein gradient. Whenever a flat derivative exists, one quarter of the reaction component agrees almost everywhere with its centered representative. As a byproduct we obtain a Wasserstein differentiability result for the squared $\WFR$ distance.
\end{abstract}

\tableofcontents

\section{Introduction}

Optimal transport has become a fundamental framework for variational problems on spaces of probability measures, with applications ranging from mean-field control and games to data analysis, sampling, and generative modeling. A key ingredient behind many of these developments is the differential structure of the quadratic Wasserstein space $\probtwo$. This can be formulated intrinsically, through the Wasserstein tangent structure and the corresponding notions of differential and gradient \cite{ambrosio-gigli-savare2005}, or through the closely related Lions differential, obtained by lifting functionals of probability measures to a Hilbert space of random variables \cite{cardaliaguet2010}. These tools make it possible to describe the sensitivity of functionals and formulate first-order optimality conditions directly at the level of probability distributions.

The Wasserstein geometry is, however, intrinsically mass preserving. While this is appropriate for evolutions governed by spatial transport alone (e.g., continuity equations), it becomes restrictive when spatial transport is accompanied by local creation, destruction, or reweighting of mass. Such mechanisms arise in reaction--transport equations, population dynamics, and partially observed control, where observations induce Bayesian reweighting of conditional distributions. A further motivation comes from mean-field control problems involving weighted particle systems, in which both particle positions and weights evolve over time. These examples motivate a first-order calculus that accounts simultaneously for spatial displacements and weight variations.

Unbalanced optimal transport extends classical transport by incorporating mechanisms for mass variation. Different choices of these mechanisms lead to different geometries, including the generalized Wasserstein distances of Piccoli and Rossi, developed in connection with transport equations with source terms \cite{piccoli-rossi2014}. In this work, we focus on the Wasserstein--Fisher--Rao (WFR) distance, also known as the Hellinger--Kantorovich (HK) distance. In its dynamic formulation, the continuity equation includes a source term, and the action penalizes both spatial velocity and relative mass variation \cite{chizat-peyre-schmitzer-vialard2018,kondratyev-monsaigneon-vorotnikov2016}. Equivalent static formulations use semicouplings or relax the marginal constraints through entropy penalties \cite{chizat-peyre-schmitzer-vialard2018unbalanced-optimal-transport,liero-mielke-savare2018}. Infinitesimally, WFR combines the transport geometry of the quadratic Wasserstein metric with the mass variation geometry of the Fisher--Rao metric. This connection also holds at the metric level: De Ponti, Sodini, and Tamanini \cite{deponti-sodini-tamanini2025} characterize the Hellinger--Kantorovich distance as the metric infimal convolution of the Hellinger and Wasserstein distances. The metric and geometric foundations of WFR are well established. Its cone representation provides a detailed description of geodesics and their relation to optimal transport on the underlying cone space \cite{liero-mielke-savare2016,laschos-mielke2019,liero-mielke-savare2018}. The associated gradient-flow theory has also received considerable attention \cite{gallouet-monsaigneon2017,laschos-mielke2026,mielke-zhu2025}.

Passing from these geometric foundations to a static first-order calculus nevertheless raises several issues. At a general source measure, optimal transport need not be described by a map, and geodesics may be non-unique. In particular, at the critical transport distance, different combinations of transport and reaction can have the same optimal cost. Moreover, creation of mass from the cone vertex is not represented by a transport--reaction field carried by the source measure. Its contribution must nevertheless be controlled in a first-order expansion whose remainder is measured by the WFR distance. A source-based differential calculus therefore requires a precise identification of the regular initial directions, a proof that the barycentric fields representing them belong to the tangent space, and estimates that remain valid in the presence of creation, annihilation, and geodesic non-uniqueness.

We address these issues by developing a regular intrinsic WFR differential at an arbitrary finite nonnegative measure. Its defining expansion pairs a tangent field with barycentric logarithms based at the source, with a remainder of smaller order than the WFR distance. We prove that this expansion extends to non-optimal semicouplings, allowing differentiation along explicit transport and reweighting perturbations. The calculus also applies to functionals defined on probability measures: its transport component recovers the Wasserstein gradient whenever the reference measure has finite second moment, while one quarter of its reaction component gives the centered flat derivative whenever the latter exists.

\subsection{Main contributions}

The main results concern endpoint tangency, strong first-order expansions,
differentiability of the squared WFR distance, and comparison with
Wasserstein and flat differentiability.

\subsubsection{Regular barycentric logarithms and tangency}

Let $\M$ denote the space of finite nonnegative Radon
measures on $\R^d$. For $\mu\in\M$, consider
$H(\mu)\coloneqq L^2(\mu;\R^d\times\R)$ with scalar product
\begin{equation*}
    \langle
        (\xi_1,\zeta_1),(\xi_2,\zeta_2)
    \rangle_\mu
    \coloneqq
    \int_{\R^d}
    \left(
        \xi_1(x)\cdot\xi_2(x)
        +\frac14\zeta_1(x)\zeta_2(x)
    \right)\,d\mu(x),
\end{equation*}
and define the regular tangent space by
\begin{equation*}
    T_\mu\M
    \coloneqq
    \overline{
        \left\{
            (\nabla\varphi,4\varphi):
            \varphi\in C_c^\infty(\R^d)
        \right\}
    }^{H(\mu)}.
\end{equation*}
A regular barycentric logarithm associated with a
constant-speed WFR geodesic $(\bfmu_t)_{t\in[0,1]}$ starting
from $\mu$ is a field $(\xi,\zeta)\in T_\mu\M$ satisfying
\begin{equation*}
    \lim_{t\downarrow0}
    \frac{\int_{\R^d}\varphi\,d\bfmu_t
          -\int_{\R^d}\varphi\,d\mu}{t}
    =
    \int_{\R^d}
        (\nabla\varphi\cdot\xi+\varphi\zeta)\,d\mu
    \qquad
    \text{for every }\varphi\in C_c^\infty(\R^d).
\end{equation*}
Thus the logarithm represents the initial variation of
smooth observables along the geodesic.

We characterize these fields through the semicoupling
formulation of WFR. Writing $\operatorname{OptSC}(\mu,\nu)$
for the optimal semicouplings between $\mu$ and $\nu$,
we associate with each $\bfgamma\in\operatorname{OptSC}(\mu,\nu)$
an explicit barycentric logarithm
$\barlog_{\bfgamma,u}$.
Here the measurable selection
$u\colon\R^d\times\R^d\to[0,1]$ accounts for the
non-uniqueness of the transport component at the
threshold $\norm{x-y}=\pi/2$; see
Definition~\ref{def:semicoupling-logarithm}.

\begin{introtheorem}[Characterization and tangency]
Let $\mu,\nu\in\M$. Every constant-speed WFR geodesic
connecting $\mu$ to $\nu$ admits a unique regular
barycentric logarithm. The set of all such logarithms is
\begin{equation*}
    \left\{
        \barlog_{\bfgamma,u}:
        \bfgamma\in\operatorname{OptSC}(\mu,\nu),\
        u\colon\R^d\times\R^d\to[0,1]
        \text{ measurable}
    \right\}.
\end{equation*}
In particular, every field in this set belongs to
$T_\mu\M$ and has $H(\mu)$-norm bounded by $\WFR(\mu,\nu)$.
\end{introtheorem}

The tangency statement holds at every prescribed finite source
measure, without absolute continuity or a Monge representation,
and for every optimal semicoupling and admissible threshold
selection. Its proof combines the cone representation with a
superposition principle for the continuity equation with reaction.

The explicit semicoupling construction also applies
to non-optimal semicouplings. For every
$\bfgamma\in\operatorname{SC}(\mu,\nu)$ and every
measurable selection $u$, its barycentric logarithm
belongs to $H(\mu)$ and satisfies
\begin{equation*}
    \norm{
        \barlog_{\bfgamma,u}
    }_\mu^2
    \leq \mathcal C(\bfgamma),
\end{equation*}
where $\mathcal C$ denotes the semicoupling cost.
This estimate is fundamental to the first-order
expansions developed below.

\subsubsection{Strong first-order expansions}

In Section~\ref{sec:wfr-calculus}, we introduce the
fundamental notion of the paper.
Proposition~\ref{propo:irrelevance-of-the-choice} shows
that the choice of $u$ at the cut locus is irrelevant
to first-order WFR expansions. We therefore fix
$u\equiv1$, choosing the transport logarithm at the
threshold, and suppress $u$ from the notation.

We say that a functional $F\colon\M\to\R$ is WFR
differentiable at $\mu$ if there exists a regular
tangent field $\mathcal V\in T_\mu\M$ such that
\begin{equation*}
    F(\nu)-F(\mu)
    =
    \left\langle
        \mathcal V,
        \barlog_{\bfgamma}
    \right\rangle_\mu
    +o\bigl(\WFR(\mu,\nu)\bigr)
    \qquad\text{as }\WFR(\mu,\nu)\to0,
\end{equation*}
uniformly over optimal semicouplings $\bfgamma$
between $\mu$ and $\nu$.

A fundamental result for applications is that this
expansion also holds along non-optimal semicouplings
whose costs vanish.

\begin{introtheorem}[Strength and uniqueness]
Let $F\colon\M\to\R$ be WFR differentiable at $\mu$.
Its tangent representative is unique and is denoted
by
\begin{equation*}
    \nabla_\WFR F(\mu)=(\nabla^T F(\mu),\nabla^R F(\mu)).
\end{equation*}
For every sequence
$\bfgamma^n\in\operatorname{SC}(\mu,\mu_n)$ satisfying
$\mathcal C(\bfgamma^n)\to0$, one has
\begin{equation*}
    F(\mu_n)-F(\mu)
    =
    \left\langle
        \nabla_\WFR F(\mu),
        \barlog_{\bfgamma^n}
    \right\rangle_\mu
    +o\left(\sqrt{\mathcal C(\bfgamma^n)}\right).
\end{equation*}
\end{introtheorem}

This result allows one to differentiate along
explicit perturbations without identifying an
optimal semicoupling. In particular for every
$(T,\rho)\in H(\mu)$ the perturbation
\begin{equation*}
    \mu^\varepsilon
    \coloneqq
    (\operatorname{Id}+\varepsilon T)_\#
    \left(
        \left(1+\frac{\varepsilon}{2}\rho\right)^2\mu
    \right)
\end{equation*}
satisfies
\begin{equation*}
    \lim_{\varepsilon\to0}
    \frac{F(\mu^\varepsilon)-F(\mu)}{\varepsilon}
    =
    \int_{\R^d}
        \left(\nabla^TF(\mu)\cdot T+\frac14\nabla^RF(\mu)\rho\right)\,d\mu.
\end{equation*}

We also establish a chain rule and a sufficient
criterion for WFR differentiability based on flat
derivatives. The criterion assumes an integral
representation along linear interpolations, joint
continuity of the flat derivative and its spatial
gradient, and locally uniform $C_b^{1,1}$ bounds
in the spatial variable. It yields WFR differentiability
of smooth expectations, cylindrical functionals,
and interaction energies.

Finally, the calculus is compatible with weighted
empirical measures. If $F$ is WFR differentiable at
$\mu=\sum_{i=1}^N e^{\ell_i}\delta_{x_i}$, then $\widetilde F(\hat x, \hat \ell)=F(\mu)$ is Fr\'echet differentiable at
$(\hat x,\hat\ell)$, with
\begin{equation*}
    D\widetilde F(\hat x,\hat\ell)
    [\hat\xi,\hat\zeta]
    =
    \sum_{i=1}^N e^{\ell_i}
    \left(
        \nabla^TF(\mu)(x_i)\cdot \xi_i+\frac14\nabla^R F(\mu)(x_i)\zeta_i
    \right).
\end{equation*}
This formula remains valid when particle positions
coincide and gives a finite-dimensional interpretation
of both components of the WFR gradient.

\subsubsection{Differentiability of the squared WFR distance}

In Section 5, we characterize WFR differentiability of the squared-distance functional
\begin{equation*}
    F_{\hat\mu}(\mu)\coloneqq\frac12\WFR^2(\mu,\hat\mu)
\end{equation*}
for a fixed target measure $\hat\mu\in\M$.

\begin{introtheorem}
The following conditions are equivalent:
\begin{enumerate}
    \item[(i)] $F_{\hat\mu}$ is WFR differentiable at $\mu$.
    \item[(ii)] There is exactly one regular barycentric logarithm $\mathcal L$
    between $\mu$ and $\hat\mu$, and
    \begin{equation*}
        \norm{\mathcal L}_\mu=\WFR(\mu,\hat\mu).
    \end{equation*}
    \item[(iii)] Every optimal cone coupling between $\mu$ and $\hat\mu$
    has the same Monge transport--growth representation at positive
    source radii and carries no target mass over the source vertex,
    as in~\eqref{eq:squared-wfr-monge-assumption}.
\end{enumerate}
Whenever these conditions hold, writing
$T$ and $q$ for the corresponding transport map and radial
growth factor, the semicoupling
\begin{equation*}
    \gamma^\star
    =
    \bigl(
        (\operatorname{Id},T)\pf\mu,
        (\operatorname{Id},T)\pf(q^2\mu)
    \bigr)
\end{equation*}
is optimal between $\mu$ and $\hat\mu$.
Its barycentric logarithm is $\barlog_{\gamma^\star}=\mathcal L$, and
\begin{equation*}
    \nabla_{\WFR}F_{\hat\mu}(\mu)
    =
    -\barlog_{\gamma^\star}.
\end{equation*}
\end{introtheorem}

Necessity follows from the metric-slope characterization of WFR
gradients: differentiability forces equality in the logarithm norm
bound, excluding both loss through barycentric averaging and target
creation. Sufficiency combines an integrated first-variation estimate
on the cone with stability of asymptotically optimal cone couplings.

\subsubsection{Probability measures and comparison of differentials}

On $\prob(\R^d)$, the unit-mass constraint leads to
the tangent space
\begin{equation*}
    T_\mu\prob(\R^d)
    \coloneqq
    \left\{
        (\xi,\zeta)\in T_\mu\M:
        \int_{\R^d}\zeta\,d\mu=0
    \right\}.
\end{equation*}
The same static expansion defines WFR differentiability
for functionals on $\prob(\R^d)$, with a tangent
representative in this space. Strength and uniqueness
remain valid, and we denote the resulting gradient
by $\widetilde\nabla_\WFR F(\mu)$.

The restriction of an ambient WFR-differentiable
functional is differentiable in this sense, with
\begin{equation*}
    \widetilde\nabla_\WFR(F_{\vert\prob})(\mu)
    =
    \Pi_\mu\nabla_\WFR F(\mu),
    \qquad
    \Pi_\mu(\xi,\zeta)
    =
    \left(\xi,\zeta-\int_{\R^d}\zeta\,d\mu\right).
\end{equation*}
The transport and reaction components can then be
compared with Wasserstein and flat derivatives,
respectively.

\begin{introtheorem}[Comparison of differentials]
Let $\mu\in\prob_2(\R^d)$ and let
$F\colon\prob(\R^d)\to\R$ be WFR differentiable at
$\mu$, with
\begin{equation*}
    \widetilde\nabla_\WFR F(\mu)
    =
    \bigl(
        \widetilde\nabla^T F(\mu),
        \widetilde\nabla^R F(\mu)
    \bigr).
\end{equation*}
Then $F_{\vert\prob_2(\R^d)}$ is Wasserstein
differentiable at $\mu$, and
\begin{equation*}
    \nabla_{W_2}F(\mu)
    =
    \widetilde\nabla^T F(\mu).
\end{equation*}
If, in addition, $F$ admits a flat derivative
$\frac{\delta F}{\delta\mu}(\mu,\cdot)\in C_b(\R^d)$,
then
\begin{equation*}
    \widetilde\nabla^R F(\mu)
    =
    4\left(
        \frac{\delta F}{\delta\mu}(\mu,\cdot)
        -
        \int_{\R^d}
            \frac{\delta F}{\delta\mu}(\mu,y)\,d\mu(y)
    \right)
    \qquad\mu\text{-a.e.}
\end{equation*}
\end{introtheorem}

Thus the transport component recovers the
Wasserstein gradient, while one quarter of the
reaction component gives the zero-mean representative
of the flat derivative whenever the latter exists.
The probability-space calculus itself requires no
second-moment assumption.

Combining this comparison with our characterization of differentiability
of the squared WFR distance, we obtain Wasserstein differentiability
of $\frac12\WFR^2(\cdot,\hat\mu)$ at every $\mu\in\prob_2(\R^d)$ satisfying
the equivalent conditions of
Theorem~\ref{thm:squared-wfr-differentiability}, with Wasserstein
gradient equal to minus the transport component of the
regular barycentric logarithm from $\mu$ to $\hat\mu$.

WFR differentiability alone does not imply flat
differentiability. We identify a class of absolutely
continuous flat directions controlled by the WFR
gradient and exhibit a counterexample based on
singular mass creation, for which a linear mixture
perturbation has WFR size of order $\sqrt{\varepsilon}$.

\subsection{Related work}

Our approach builds on the intrinsic differential calculus of the quadratic Wasserstein space developed in \cite{ambrosio-gigli-savare2005}, where tangent fields, subdifferentials, and first-order expansions of functionals are described through optimal transport plans. This provides the model for our static formulation of WFR differentiability. The closely related Lions approach uses a lifting to random variables; see \cite{cardaliaguet2010}. Differential tools on spaces of measures also play an important role in the derivation of necessary optimality conditions for measure-valued control problems \cite{bonnet-frankowska2021}. For the squared Wasserstein distance, Alfonsi and Jourdain \cite{alfonsi-jourdain2020} characterize Lions differentiability in terms of uniqueness of a map-induced optimal coupling. Our characterization of differentiability of the squared WFR distance provides an unbalanced counterpart, in which the regular logarithm must additionally capture the full metric displacement.

The cone representation and geodesic lifting results of Liero, Mielke, and Savar\'e \cite{liero-mielke-savare2018} provide the starting point for our characterization of regular logarithms. Their earlier work already describes the non-uniqueness of critical-distance Dirac geodesics \cite[Section~5.2]{liero-mielke-savare2016}. Further geometric properties of cones and HK spaces are established in \cite{laschos-mielke2019}, while \cite{liero-mielke-savare2023} develops regularity results for optimal potentials and transport--growth representations from interior points of geodesics. Related regularity and Monge-type results are obtained by Gallouët, Ghezzi, and Vialard \cite{gallouet-ghezzi-vialard2025}.

The infinitesimal transport--reaction structure has been extensively used in the study of evolution equations. Kondratyev, Monsaigneon, and Vorotnikov \cite{kondratyev-monsaigneon-vorotnikov2016} develop a formal Riemannian calculus based on weighted Sobolev potentials. Variational constructions and qualitative properties of WFR gradient flows are studied, among others, in \cite{gallouet-monsaigneon2017,laschos-mielke2026,mielke-zhu2025}. These works establish the role of the transport--reaction gradient structure in variational evolution problems.

A particularly close reference is Flei{\ss}ner~\cite{fleissner2023}, who proves tangency of minimal velocity fields at almost every time along absolutely continuous curves, establishes a metric linearization theorem, and differentiates the squared HK distance along such curves. Her earlier work~\cite{fleissner2021} also identifies elements of the Fr\'echet subdifferential at $h=0$ of $h\mapsto-\tfrac12\mathrm{HK}^2(\mu_h,\nu)$ for transport--reaction perturbations
$\mu_h$, and provides a quadratic remainder estimate for smooth test integrals expressed through optimal cone plans. This estimate is closely related to the smooth-observable expansions used here. Our endpoint tangency theorem concerns the specific barycentric fields associated with every optimal semicoupling at an arbitrary prescribed source measure, including all admissible selections at the
transport threshold.

Logarithmic maps have also been developed in the linearized optimal transport literature. Cai, Cheng, Schmitzer, and Thorpe \cite{cai-cheng-schmitzer-thorpe2022} introduce HK logarithmic and exponential maps for absolutely continuous reference measures, retaining an additional measure-valued component for singular creation, and discuss barycentric approximations when transport maps are unavailable. Sarrazin and Schmitzer \cite{sarrazin-schmitzer2024} explicitly define barycentric logarithmic maps for nondeterministic soft plans, without requiring optimality in their definition. Our semicoupling construction retains the pairing information at the critical distance and incorporates measurable selections describing the possible regular initial transport components. We use this construction to characterize the full set of regular source logarithms and prove tangency for every optimal choice.

Other recent works develop complementary differential structures on HK spaces. Dello Schiavo and Sodini \cite{delloschiavo-sodini2025} establish universal infinitesimal Hilbertianity and density of cylinder functions in the corresponding Sobolev spaces. For a canonical reference measure on the space of measures, they also identify a natural Dirichlet form with the Cheeger energy. Saidi, Mena, and Gunsilius \cite{saidi-mena-gunsilius2026} relate the local differential geometry of HK to Wasserstein geometry on the cone, obtaining tangent-space isometries and applications to parallel transport under suitable regularity assumptions. In particular, their admissible class yields lifts uniformly separated from the cone vertex. Our results concern pointwise static expansions at arbitrary finite source measures, allowing singular sources and vertex contributions in the representing plans.

Finally, our comparison with flat derivatives is related to the calculus of functions on spaces of finite measures. Ren and Wang \cite{ren-wang2021} establish relationships between extrinsic, linear functional, intrinsic, and $L$-derivatives under appropriate regularity assumptions. Our sufficient criterion identifies conditions under which a flat derivative yields the static WFR expansion required by our definition, including uniformity over optimal semicouplings and logarithmic selections.

\section{The Wasserstein--Fisher--Rao geometry}\label{sec:wfr-geometry}

In this section we recall some formulations and basic properties of the Wasserstein--Fisher--Rao distance. The material is mainly borrowed from \cite{chizat-peyre-schmitzer-vialard2018unbalanced-optimal-transport,liero-mielke-savare2018}.

\subsection{The semicoupling formulation}

\begin{definition}\label{def:semi-coupling}
Let $\mu,\nu\in\M$ be given. A \emph{semicoupling} between $\mu$ and $\nu$ is a pair
\begin{equation*}
    \bfgamma=(\gamma_0,\gamma_1)\in \mathcal M(\R^d\times\R^d)^2
\end{equation*}
such that
\begin{equation*}
\proj{1}\pf\gamma_0=\mu
\qquad\text{and}\qquad
\proj{2}\pf\gamma_1=\nu.
\end{equation*}
We denote the set of semicouplings between $\mu$ and $\nu$
by $\operatorname{SC}(\mu,\nu)$.
\end{definition}

We introduce the cost
$c\colon\R^d\times\R^d\times[0,\infty)\times[0,\infty)\to[0,\infty)$ as
\begin{equation}\label{eq:semi-coupling-cost}
c(x,y,a,b)
\coloneqq
a+b-2\sqrt{ab}\cos\left(\left(\norm{y-x}\wedge\frac{\pi}{2}\right)\right)
\end{equation}
The cost $c$ is jointly $1$-homogeneous in the variables $a$ and $b$, namely
\begin{equation*}
    c(x,y,\lambda a,\lambda b)
    =
    \lambda c(x,y,a,b)
    \qquad\text{for every }\lambda\geq 0.
\end{equation*}

Therefore, if $\lambda\in\mathcal M(\R^d\times\R^d)$ dominates both $\gamma_0$ and $\gamma_1$, and we set
\begin{equation*}
a\coloneqq\frac{d\gamma_0}{d\lambda},
\qquad
b\coloneqq\frac{d\gamma_1}{d\lambda},
\end{equation*}
then the semicoupling cost
\begin{equation*}
\mathcal{C}(\bfgamma)
\coloneqq
\int_{\R^d\times\R^d}
c(x,y,a(x,y),b(x,y)) \; d\lambda(x,y)
\end{equation*}
is independent of the choice of $\lambda$ and only depends on $\bfgamma$. Equivalently,
\begin{equation}\label{eq:total-mass-formualtion}
    \mathcal C(\bfgamma) = \mu(\R^d)+\nu(\R^d) -2\int_{\R^d\times\R^d} \sqrt{a(x,y)b(x,y)}\cos\left(\left(\norm{y-x}\wedge\frac{\pi}{2}\right)\right)\,d\lambda(x,y).
\end{equation}

The Wasserstein--Fisher--Rao distance is defined as follows.

\begin{definition}\label{def:semi-coupling-formulation}
Given $\mu,\nu\in\M$, we define their Wasserstein--Fisher--Rao distance as
\begin{equation}\label{eq:semi-coupling-formulation}
\WFR^2(\mu,\nu)
\coloneqq
\inf_{\bfgamma}\mathcal C(\bfgamma),
\end{equation}
where the infimum ranges over all semicouplings between $\mu$ and $\nu$. We denote the set of minimizers
in \eqref{eq:semi-coupling-formulation} by $\operatorname{OptSC}(\mu,\nu)$.
\end{definition}

\begin{remark}[Interpretation of semicouplings]
Unlike classical transport plans, semicouplings allow the
amount of mass leaving a point $x$ to differ from the
amount arriving at $y$, thereby accounting for mass
variation during transport.

To interpret the cost, let $a=d\gamma_0/d\lambda$ and
$b=d\gamma_1/d\lambda$. The cost integrand can be decomposed as
\begin{equation*}
\begin{split}
    a+b-2\sqrt{ab}\cos(\norm{y-x}\wedge\frac\pi2)
    =(\sqrt a-\sqrt b)^2+2\sqrt{ab}
    \bigl(1-\cos(\norm{y-x}\wedge\frac\pi2)\bigr).
\end{split}
\end{equation*}
The first term penalizes mass variation, while the
second accounts for spatial displacement. In particular,
for small displacements,
\begin{equation*}
    2\sqrt{ab}\bigl(1-\cos(\norm{y-x})\bigr)
    =\sqrt{ab}\norm{y-x}^2
    +O\bigl(\sqrt{ab}\norm{y-x}^4\bigr),
\end{equation*}
recovering the classical quadratic transport cost at
leading order when $a=b$.
Finally, if either density vanishes, the cost describes
pure annihilation or creation of mass, while for
$\norm{y-x}\ge\pi/2$ the interaction vanishes, so transporting
mass is no more advantageous than annihilating it at
the source and recreating it at the target.
\end{remark}

\subsection{The cone formulation}

\subsubsection{The cone geometry}

Let $Y\coloneqq\R^d\times[0,\infty)$. Consider the function $\dc\colon Y\times Y\to[0,\infty)$ defined by
\begin{equation}\label{eq:distance}
\dc^2\bigl((x,r),(y,s)\bigr)
\coloneqq
r^2+s^2-2rs\cos\bigl(\norm{x-y}\wedge \pi\bigr).
\end{equation}
The function $\dc$ is nonnegative, symmetric and satisfies the triangle inequality. It is only a semidistance on $Y$, since all the points of the form $(x,0)$ have zero distance from each other. We therefore consider the equivalence relation
\begin{equation}\label{eq:cone-equivalence-relation}
(x,r)\sim(y,s)
\quad\Longleftrightarrow\quad
r=s=0
\quad\text{or}\quad
x=y,\;r=s,
\end{equation}
and define the cone
\begin{equation*}
    \cone\coloneqq Y/\sim
\end{equation*}

equipped with the induced distance, which we still denote by $\dc$. We denote a generic element of $\cone$ by $\y=[x,r]$,
while the cone vertex is denoted by $\vertex$. The  projection map is $\mathfrak p\colon Y \to \cone$:
\begin{equation*}
    \mathfrak p\left((x,r)\right)= [x,r],
\end{equation*}
which admits right inverses $\mathfrak i_{\bar x}\colon \cone\to Y$ defined as
\begin{equation}\label{eq:right-inverses}
    \mathfrak i_{\bar x}(\y):=
    \begin{cases}
    (x,r),&\y=[x,r],\ r>0,\\
    (\bar x,0),&\mathfrak y=\vertex.
    \end{cases}
\end{equation}
The radial component of $i_{\bar x}$ is continuous,
while $i_{\bar x}$ itself is in general only Borel
measurable due to its behavior at the vertex.

\subsubsection{Wasserstein--Fisher--Rao as Wasserstein on the cone}

It will also be useful to consider the truncated cone distance
\begin{equation}\label{eq:truncated-cone-distance}
\dcwfr^2\bigl([x,r],[y,s]\bigr)
\coloneqq
r^2+s^2-2rs\cos\bigl(\norm{x-y}\wedge \frac{\pi}{2}\bigr).
\end{equation}
Notice that
\begin{equation}\label{eq:semicoupling-cone-cost}
c(x,y,a,b)
=
\dcwfr^2\bigl([x,\sqrt a],[y,\sqrt b]\bigr).
\end{equation}
The distinction between $\dc$ and $\dcwfr$ is useful: the first one is the canonical cone distance and is the natural one for the description of geodesics, whereas the second one is directly related to the semicoupling formulation. Nevertheless, the two distances induce the same Wasserstein--Fisher--Rao distance on $\M$, as recalled below.

Let $A\in\mathcal M(\cone^N)$ and write
\begin{equation*}
        \y_i=[x_i,r_i],
    \qquad i=1,\ldots,N.
\end{equation*}

We say that $A\in\mathcal M _2(\cone^N)$ if

$$
    \int_{\cone^N}\sum_{i=1}^N r_i^2\,dA<+\infty.
$$

For $i=1,\ldots,N$, we define the $i$-th \emph{homogeneous marginal} of $A$ as
\begin{equation}\label{eq:homogeneous-marginal}
\h_i(A)
\coloneqq
(x_i)\pf(r_i^2A)\in\M.
\end{equation}
Here the value assigned to $x_i$ when $\y_i=\vertex$ is irrelevant, since the measure $r_i^2A$ does not charge the vertex. Equivalently, for every bounded and continuous function $\varphi\colon\R^d\to\R$,
\begin{equation}\label{eq:homogeneous-marginal-dual}
\int_{\R^d}\varphi(x)\,d\h_i(A)(x)
=
\int_{\cone^N}r_i^2\varphi(x_i)\,dA(\y_1,\ldots,\y_N).
\end{equation}
When $N=1$ we simply write $\h(A)$.

For $i=1,\ldots,N$, define the coordinate projections on
$Y^N$ by
\begin{equation*}
    \mathsf{x}_i\bigl((x_j,r_j)_{j=1}^N\bigr)\coloneqq x_i,
    \qquad
    \mathsf{r}_i\bigl((x_j,r_j)_{j=1}^N\bigr)\coloneqq r_i.
\end{equation*}
Fix $\bar x\in\mathbb R^d$ and let
$\mathfrak i_{\bar x}^{\otimes N}$ and
$\mathfrak p^{\otimes N}$ denote the product maps induced by
the right inverse $\mathfrak i_{\bar x}$ and the quotient map
$\mathfrak p$. For $A\in\mathcal M_2(\mathfrak C^N)$, its
coordinate lift
$\widetilde A\coloneqq(\mathfrak i_{\bar x}^{\otimes N})_\#A$
satisfies
\begin{equation*}
    (\mathfrak p^{\otimes N})_\#\widetilde A=A,
    \qquad
    h_i(A)=(\mathsf{x}_i)_\#
        \bigl(\mathsf{r}_i^2\widetilde A\bigr),
    \quad i=1,\ldots,N.
\end{equation*}
Conversely, for any finite nonnegative measure $\widetilde A$
on $(\mathbb R^d\times[0,\infty))^N$ satisfying
$\int\sum_{i=1}^N\mathsf{r}_i^2\,d\widetilde A<\infty$,
\begin{equation*}
    h_i\bigl((\mathfrak p^{\otimes N})_\#\widetilde A\bigr)
    =(\mathsf{x}_i)_\#
        \bigl(\mathsf{r}_i^2\widetilde A\bigr),
    \quad i=1,\ldots,N.
\end{equation*}
These identities are independent of the spatial representative
at the vertex, since the weight $\mathsf{r}_i^2$ vanishes there.
Thus passing to the quotient, or choosing a coordinate lift,
preserves all homogeneous marginals.

\begin{definition}\label{def:cone-coupling}
Let $\mu,\nu\in\M$. A measure $A\in \mathcal M_2(\cone^2)$ is called a \emph{cone coupling} between $\mu$ and $\nu$ if
\begin{equation*}
        \h_1(A)=\mu
    \qquad\text{and}\qquad
    \h_2(A)=\nu.
\end{equation*}
We denote the set of cone couplings between $\mu$ and $\nu$ by
\begin{equation*}
    \Gamma_{\cone}(\mu,\nu)
    \coloneqq
    \left\{
        A\in\mathcal M_2(\cone^2):
        \h_1(A)=\mu,\;
        \h_2(A)=\nu
    \right\}.
\end{equation*}
\end{definition}

We have the following characterization of the Wasserstein--Fisher--Rao distance; see \cite[Problem 7.4 and Lemma 7.19]{liero-mielke-savare2018}.

\begin{theorem}[Cone formulation]\label{thm:wfr-cone-formulation}
For every $\mu,\nu\in\M$ it holds
\begin{equation}\label{eq:cone-formulation}
\WFR^2(\mu,\nu)
=
\min_{A\in\Gamma_{\cone}(\mu,\nu)}
\int_{\cone^2}
\dc^2(\y_0,\y_1)\,dA(\y_0,\y_1).
\end{equation}
Moreover,
\begin{equation}\label{eq:truncated-cone-formulation}
\WFR^2(\mu,\nu)
=
\min_{A\in\Gamma_{\cone}(\mu,\nu)}
\int_{\cone^2}
\dcwfr^2(\y_0,\y_1)\,dA(\y_0,\y_1).
\end{equation}
\end{theorem}

Accordingly, we say that a cone coupling is \emph{optimal} if it attains the minimum in \eqref{eq:cone-formulation}. We denote by
\begin{equation*}
    \mathrm{Opt\Gamma}_{\cone}(\mu,\nu)
\end{equation*}
the collection of optimal cone couplings between $\mu$ and $\nu$.

From the previous theorem (see again \cite[Lemma 7.19]{liero-mielke-savare2018}), if we introduce the set of \emph{ineffective pairs}
\begin{equation}\label{eq:ineffective-pairs}
    \mathfrak B\coloneqq \left\{\bigl([x,r], [y,s]\bigr)\in \cone\times\cone\colon rs>0, \norm{x-y}>\frac{\pi}{2}\right\}
\end{equation}
we get that if $A\in \mathrm{Opt\Gamma}_{\cone}(\mu,\nu)$, then $A(\mathfrak B)=0$, that is, any optimal cone plan does not pair points whose spatial coordinates are more than $\frac{\pi}{2}$ apart. We say that a pair $(\y_0,\y_1)\in \cone^2$ is \emph{reduced} if it does not belong to $\mathfrak B$.

\begin{remark}
The homogeneous formulation has an important geometric interpretation.
The cone variable $r$ should be thought of as the square root of the mass:
a cone point $[x,r]$ represents the measure $r^2\delta_x$.
Hence transport on the cone simultaneously encodes displacement of the
position $x$ and variation of the mass through the radial coordinate $r$.
Together with Theorem~\ref{thm:wfr-cone-formulation}, this provides a
close geometric relation between $\WFR$ and the Wasserstein distance
on the cone. This geometric relation also suggests a close connection
between the corresponding differential structures.
Nonetheless, applying standard Wasserstein differential calculus on
the cone requires care because of the singularity at the vertex:
the cone is not a smooth Riemannian manifold there.
This also prevents a direct application of the usual Lions lifting
approach, since $L^2(\Omega;\cone)$ is a metric space rather than
a Hilbert space with its natural cone-induced distance.
\end{remark}

\subsection{Wasserstein--Fisher--Rao geodesics}\label{subsec:wfr-geodesics}

We now recall the characterization of Wasserstein--Fisher--Rao geodesics in terms of geodesics on the cone, following \cite[Theorem 8.6]{liero-mielke-savare2018}, and record
corresponding formulas and estimates used in the remainder of the paper.

A constant-speed WFR geodesic connecting $\mu$ and $\nu$ is a
curve $\bfmu\colon[0,1]\to\M$ satisfying, for every $t_1,t_2\in[0,1]$,
\begin{equation*}
    \bfmu_0=\mu,\qquad \bfmu_1=\nu,\qquad
    \WFR(\bfmu_{t_2},\bfmu_{t_1})
    =|t_2-t_1|\WFR(\mu,\nu).
\end{equation*}

We denote by
\begin{equation*}
\Geo(\cone)
\coloneqq
\left\{
\bfy\colon[0,1]\to\cone
\colon
\dc\left(\bfy_{t_1},\bfy_{t_2}\right)
=
\lvert t_2-t_1\rvert
\dc\left(\bfy_1,\bfy_0\right),
\quad t_1,t_2\in [0,1]
\right\}
\end{equation*}
the set of constant-speed cone geodesics. With the uniform
topology, this is a closed subset of the Polish space
$C([0,1];\cone)$, hence itself Polish.
Since $\R^d$ is uniquely geodesic, so is $(\cone,\dc)$;
see \cite[Chapter~I.5, Corollary~5.11]{bridson1999}.
The endpoint map
\begin{equation*}
    E\colon\Geo(\cone)\to\cone\times\cone,
    \qquad E(\bfy)\coloneqq(\bfy_0,\bfy_1),
\end{equation*}
is therefore a continuous bijection between Polish spaces.
Its inverse $G$, which assigns to each pair its unique cone
geodesic, is Borel measurable by the Lusin--Souslin theorem
\cite[Corollary~15.2]{kechris1995}.

For $t\in[0,1]$, let
\begin{equation*}
    e_t\colon\Geo(\cone)\to\cone,
    \qquad e_t(\bfy)\coloneqq\bfy_t,
\end{equation*}
and define the Borel coordinate evaluation maps
\begin{equation*}
    \mathsf{x}_t\coloneqq\operatorname{pr}_x\circ\mathfrak i_0\circ e_t,
    \qquad
    \mathsf{r}_t\coloneqq\operatorname{pr}_r\circ\mathfrak i_0\circ e_t,
\end{equation*}
where $\operatorname{pr}_x,\operatorname{pr}_r$ are the coordinate
projections on $\R^d\times[0,\infty)$.
For $\bfy\in\Geo(\cone)$, write
\begin{equation*}
    \bfx_t\coloneqq\mathsf{x}_t(\bfy),\qquad
    \bfr_t\coloneqq\mathsf{r}_t(\bfy),\qquad
    \bfy_t=[\bfx_t,\bfr_t].
\end{equation*}
Both coordinate evaluation maps are Borel measurable.
The radial coordinate is intrinsic, whereas the spatial
coordinate depends on the chosen representative only at
the vertex. In particular, homogeneous expressions such
as $\bfr_t^2\varphi(\bfx_t)$ are independent of this choice.

For a finite measure $\boldsymbol{\Pi}\in\mathcal M(\Geo(\cone))$ satisfying
\begin{equation}\label{eq:Pi-integrability}
    (e_t)\pf\boldsymbol\Pi\in\mathcal M_2 (\cone)\qquad t\in[0,1],
\end{equation}
we define its homogeneous marginal at time $t$ by
\begin{equation}\label{eq:dynamic-homogeneous-marginal}
\h_t(\boldsymbol{\Pi})
\coloneqq
\h\bigl((e_t)\pf\boldsymbol{\Pi}\bigr).
\end{equation}
Equivalently, for every $\varphi\in C_b(\R^d)$,
\begin{equation}\label{eq:dynamic-homogeneous-marginal-dual}
\int_{\R^d}\varphi(x)\,d\h_t(\boldsymbol{\Pi})(x)
=
\int_{\Geo(\cone)}
\bfr_t^2\varphi(\bfx_t)\,d\boldsymbol{\Pi}(\bfy).
\end{equation}

The following characterization, which states that every Wasserstein--Fisher--Rao geodesic is the homogeneous marginal of a superposition of geodesics on the cone, can be found in \cite[Theorem 8.6]{liero-mielke-savare2018}.

\begin{theorem}[Characterization of WFR geodesics]
\label{thm:wfr-geodesics-cone}
Let $\mu_0,\mu_1\in\M$. A curve
$\bfmu\colon[0,1]\to\M$ is a constant-speed WFR geodesic
connecting $\mu_0$ to $\mu_1$ if and only if there exists
\begin{equation*}
    \boldsymbol{\Pi}\in\prob(\Geo(\cone))
\end{equation*}
such that
\begin{equation}\label{eq:wfr-geodesic-homogeneous-marginal}
    \bfmu_t=\h_t(\boldsymbol{\Pi})
    \qquad\text{for every }t\in[0,1],
\end{equation}
and
\begin{equation*}
    (e_0,e_1)\pf\boldsymbol{\Pi}
    \in\mathrm{Opt\Gamma}_{\cone}(\mu_0,\mu_1).
\end{equation*}
For every such representation and every
$0\leq t_1<t_2\leq1$, one also has
\begin{equation}\label{eq:intermediate-optimal-cone-plan}
    (e_{t_1},e_{t_2})\pf\boldsymbol{\Pi}
    \in\mathrm{Opt\Gamma}_{\cone}
    (\bfmu_{t_1},\bfmu_{t_2}). 
\end{equation}
Moreover, every
$A\in\mathrm{Opt\Gamma}_{\cone}(\mu_0,\mu_1)$ admits a
representation
\begin{equation*}
    A=(e_0,e_1)\pf\boldsymbol{\Pi}_A,
    \qquad
    \boldsymbol{\Pi}_A\in
    \mathcal M(\Geo(\cone)).
\end{equation*}
The curve $\bfmu_t\coloneqq\h_t(\boldsymbol{\Pi}_A)$
is a constant-speed WFR geodesic connecting $\mu_0$ to
$\mu_1$, and
\eqref{eq:intermediate-optimal-cone-plan} holds.
If $A\in\prob(\cone^2)$, then
$\boldsymbol{\Pi}_A$ is a probability measure.
\end{theorem}

We next record the cone-geodesic formulas and estimates used to
construct and characterize regular barycentric logarithms in
Section~\ref{sec:regular-logarithms}.

The set of reduced cone geodesics is
\begin{equation*}
    \Geo^*(\cone)\coloneqq \left\{\bfy\in\Geo(\cone)\colon (\bfy_0,\bfy_1)\notin\mathfrak B\right\}.
\end{equation*}

Every probability representation $\boldsymbol{\Pi}$ in
Theorem~\ref{thm:wfr-geodesics-cone} is concentrated on this set:
its optimal endpoint plan $A\coloneqq(e_0,e_1)\pf\boldsymbol{\Pi}$
satisfies $A(\mathfrak B)=0$, hence
\begin{equation*}
    \boldsymbol{\Pi}\bigl(\Geo^*(\cone)\bigr)=1.
\end{equation*}

Define the vector-valued sine function
\begin{equation*}
\mathbf{sin}(z)
\coloneqq
\begin{cases}
\sin(\norm{z})\dfrac{z}{\norm{z}}, & z\neq 0,\\[2mm]
0, & z=0.
\end{cases}
\end{equation*}
and, for $y\in\R^d$, the kernel
\begin{equation*}
    K_y(x)\coloneqq\cos\left(\norm{x-y}\wedge\pi\right).
\end{equation*}
Then $K_y\in C_b^{1,1}(\R^d)$ and, uniformly in $y$,
\begin{equation}\label{eq:kernel-bounds}
    |K_y(x)|^2+\norm{\nabla K_y(x)}^2=1,
    \qquad
    \Lip(\nabla K_y)\leq1.
\end{equation}
Indeed, on $0<\norm{x-y}<\pi$, the radial and tangential
eigenvalues of the Hessian are respectively
$-\cos(\norm{x-y})$ and
$-\sin(\norm{x-y})/\norm{x-y}$, both bounded in absolute
value by one. The gradient extends continuously across
$\norm{x-y}=0,\pi$ and vanishes on $\norm{x-y}\geq\pi$.

For the cone geodesic $\bfy_t=[\bfx_t,\bfr_t]$ joining
$\y_0=[x,r]$ to $\y_1=[y,s]$, the explicit formulas are as follows;
see \cite[Section~2.1.2, Eq.~(2.7)]{liero-mielke-savare2023}.

When $r,s>0$ and $\norm{x-y}<\pi$, they follow by
identifying the cone over the segment joining $x$ and $y$
with a Euclidean sector, in which the geodesic is a
straight line segment. When $\norm{x-y}\geq\pi$, the
geodesic passes through the vertex.

The radial component satisfies
\begin{equation*}
    \bfr_t^2=\bfm_t(\y_0,\y_1)
    \coloneqq (1-t)^2r^2+t^2s^2+2rst(1-t)K_y(x).
\end{equation*}
For the spatial component we need instead to distinguish different behaviors depending on the endpoints.

If $r,s>0$, set
\begin{equation*}
u\coloneqq \frac{s}{r}K_y(x)-1,
\qquad
v\coloneqq \frac{s}{r}\nabla K_y(x).
\end{equation*}
If $\norm{y-x}=0$, then $\bfx_t\equiv x$. If $0<\norm{y-x}<\pi$, then
\begin{equation*}
    \bfx_t
    =
    x+\operatorname{Arg}\left((1-t)r+ts e^{\mathrm{i}\norm{y-x}}\right)
    \frac{y-x}{\norm{y-x}},
    \qquad t\in[0,1],
\end{equation*}
where $\operatorname{Arg}$ denotes the principal argument,
which belongs to $[0,\pi]$ for the displayed complex
numbers.
If $\norm{y-x}\geq\pi$, setting $\tau\coloneqq r/(r+s)$, we have
\begin{equation*}
    \bfy_t
    =
    \begin{cases}
        [x,r-t(r+s)], & 0\leq t<\tau,\\
        \vertex, & t=\tau,\\
        [y,t(r+s)-r], & \tau<t\leq1.
    \end{cases}
\end{equation*}

Finally, if one of the endpoints is the vertex, the geodesic is simply radial:
\begin{equation*}
\bfy_t([x,r],\vertex)=[x,\bfr_t]=[x,(1-t)r],
\qquad
\bfy_t(\vertex,[y,s])=[y,\bfr_t]=[y,ts].
\end{equation*}
If both endpoints coincide with $\vertex$, the geodesic is constant: $\bfy_t\equiv\vertex$.

Therefore, apart from the trivial geodesic $\bfy_t\equiv \vertex$, no reduced geodesic hits the vertex in $(0,1)$.

The next lemma records the velocities along cone geodesics with
nonzero initial radius. Its proof and that of
Lemma~\ref{lemma:cone-observable-estimates} are given in
Appendix~\ref{app:cone-geodesic-calculations}.

\begin{lemma}[Velocities along cone geodesics]
\label{lemma:cone-geodesic-velocities}
Let $\bfy_t=[\bfx_t,\bfr_t]\in\Geo(\cone)$ be a cone
geodesic with $\bfy_0=[x,r]\neq\vertex$ and
$t\in[0,1)$ such that $\bfr_t>0$.
Derivatives at $t=0$ are understood as right derivatives.

If $\bfy_1=[y,s]\neq\vertex$, then, with $u,v$ defined above,
\begin{equation*}
\dot{\bfx}_t
=
\frac{v}{(1+tu)^2+t^2\norm{v}^2},
\qquad
\dot{\bfr}_t
=
r
\frac{(1+tu)u+t\norm{v}^2}
{\sqrt{(1+tu)^2+t^2\norm{v}^2}},
\end{equation*}
and
\begin{equation*}
\dot{\bfm}_t
=
-2(1-t)r^2+2ts^2
+2(1-2t)rs
K_y(x).
\end{equation*}
If $\bfy_1=\vertex$, then the geodesic is radial and, for
$t\in[0,1)$, $\dot{\bfx}_t=0$ and $\dot{\bfr}_t=-r$.
In both cases,
\begin{equation}\label{eq:constant-speed-norm}
    \bfr_t^2\norm{\dot{\bfx}_t}^2
    +|\dot{\bfr}_t|^2
    =\dc^2(\bfy_0,\bfy_1),
    \qquad t\in[0,1)
    \text{ with }\bfr_t>0.
\end{equation}
\end{lemma}

For $\y_0=[x,r]$ and $\y_1=[y,s]$ with $r,s>0$,
Lemma~\ref{lemma:cone-geodesic-velocities} identifies the cone
logarithm in spatial and radial coordinates as
\begin{equation*}
    \Log_{\cone}(\y_0,\y_1)
    \coloneqq
    \left(\dot{\bfx}_0,\dot{\bfr}_0\right)
    =
    \left(
        \frac{s}{r}
        \nabla K_y(x),
        s
        K_y(x)
        -r
    \right).
\end{equation*}
If $\y_1=\vertex$, the radial geodesic instead gives
\begin{equation*}
    \Log_{\cone}(\y_0,\vertex)=(0,-r).
\end{equation*}
We denote the components of the cone logarithm as follows:
\begin{equation*}
    \Log_{\cone}(\y_0,\y_1)=(\Log^T_{\cone}(\y_0,\y_1), \Log^R_{\cone}(\y_0,\y_1)).
\end{equation*}

The next lemma gives the initial derivative and a quadratic expansion
of homogeneous observables along cone geodesics.

\begin{lemma}
\label{lemma:cone-observable-estimates}
Let $\varphi\in C_b^{1,1}(\R^d)$ and define
\begin{equation*}
    \Phi_\varphi([x,r])\coloneqq r^2\varphi(x),
    \qquad
    \Phi_\varphi(\vertex)\coloneqq0.
\end{equation*}
Let $\bfy_t\in\Geo(\cone)$ be a cone geodesic joining
$\bfy_0=[x,r]$ and $\bfy_1=[y,s]$. The map $t\mapsto\Phi_\varphi(\bfy_t)$ has a right derivative at
$t=0$. If $r>0$, then
\begin{equation*}
    \left.\frac{d}{dt}\Phi_\varphi(\bfy_t)\right|_{t=0^+}
    =
    r^2\nabla\varphi(x)\cdot
    \Log_\cone^T(\bfy_0,\bfy_1)
    +
    2r\varphi(x)
    \Log_\cone^R(\bfy_0,\bfy_1),
\end{equation*}
while, if $\bfy_0=\vertex$, the right derivative is zero.

Moreover, the map $t\mapsto\Phi_\varphi(\bfy_t)$ is differentiable in $(0,1)$ with
\begin{equation*}
    \left|
        \frac{d}{dt}\Phi_\varphi(\bfy_t)
    \right|
    \leq B_\varphi(r^2+s^2)
\end{equation*}
at every interior time, where $B_\varphi\coloneqq 3\norm{\varphi}_\infty
        +\frac12\norm{\nabla\varphi}_\infty$, and consequently
\begin{equation}\label{eq:incremental-ratio-bound}
    \left|
    \frac{\Phi_\varphi(\bfy_t)-\Phi_\varphi(\bfy_0)}{t}
    \right|
    \leq B_\varphi(r^2+s^2),
    \qquad t\in(0,1].
\end{equation}
Finally,
\begin{equation*}
    \left|
        \Phi_\varphi(\bfy_t)-\Phi_\varphi(\bfy_0)
        -t\left.\frac{d}{d\tau}\right|_{\tau=0^+}
        \Phi_\varphi(\bfy_\tau)
    \right|
    \leq
    C_\varphi t^2 \dc(\bfy_0,\bfy_1)^2,
    \qquad t\in[0,1],
\end{equation*}
where one may take
\begin{equation}\label{eq:quantitative-dependence}
    C_\varphi
    =
    \norm{\varphi}_\infty
    +\frac12\norm{\nabla\varphi}_\infty
    +\frac12\Lip(\nabla\varphi).
\end{equation}
\end{lemma}

\section{Regular logarithms and tangency}\label{sec:regular-logarithms}

In this section we study the notion of regular Wasserstein--Fisher--Rao logarithms and characterize them in terms of optimal semicouplings and measurable selections at the cut locus $\norm{x-y}=\frac{\pi}{2}$. The key ingredient is the tangency result of Theorem~\ref{propo:tangency}, which shows that the barycentric logarithm associated with any optimal semicoupling, and any admissible selection at the cut locus, belongs to the WFR tangent space at the source measure. Combining this with the representation of WFR geodesics through probabilities supported on $\Geo(\cone)$ yields the characterization of regular logarithms stated in Theorem~\ref{thm:regular-barycentric-log-characterization}. We end the section by showing that
the choice at the cut locus is irrelevant for first-order expansions.

\subsection{Regular and semicoupling logarithms}

We first recall the definition of the tangent space \cite{kondratyev-monsaigneon-vorotnikov2016}. 
On $\R^d\times\R$, we consider the scalar product
\begin{equation}\label{eq:wfr-pointwise-scalar-product}
    \langle (\xi_1,\zeta_1),(\xi_2,\zeta_2)\rangle_{\WFR}
    \coloneqq
    \xi_1\cdot\xi_2+\frac14\zeta_1\zeta_2,
\end{equation}
with associated norm $\norm{\cdot}_{\WFR}$.
Then, for $\mu\in\M$, we introduce the Hilbert space 
\begin{equation*}
    H(\mu)\coloneqq L^2(\mu; \R^d\times\R)
\end{equation*}
equipped with the scalar product
\begin{equation}
\langle \mathcal V_1,\mathcal V_2\rangle_\mu
\coloneqq
\int_{\R^d}\langle
(\xi_1(x),\zeta_1(x)),(\xi_2(x),\zeta_2(x))\rangle_{\WFR}\,d\mu(x),
\end{equation}
where $\mathcal V_i=(\xi_i,\zeta_i)$, $i=1,2$.
The associated norm is denoted by
$\norm{\cdot}_\mu$. The tangent space to $\M$ at $\mu$ is then defined by
\begin{equation}
T_\mu\M
\coloneqq
\overline{
\left\{
(\nabla\varphi,4\varphi)
\colon
\varphi\in C_c^\infty(\R^d)
\right\}
}^{H(\mu)}.
\end{equation}

We now introduce the fundamental notion of this section.

\begin{definition}[Regular barycentric logarithm]
We say that $\mathcal V=(\xi,\zeta)\in T_\mu\M$ is a
\emph{regular barycentric logarithm} between $\mu$ and $\nu$
if there exists a constant-speed WFR geodesic
$\bfmu\colon[0,1]\to\M$, with $\bfmu_0=\mu$ and $\bfmu_1=\nu$,
such that, for every $\varphi\in C_c^\infty(\R^d)$,
\begin{equation}\label{eq:regular-barycentric-logarithm}
    \lim_{t\downarrow0}
    \frac{\int_{\R^d}\varphi\,d\bfmu_t
    -\int_{\R^d}\varphi\,d\mu}{t}
    =
    \int_{\R^d}
    \left(\nabla\varphi\cdot\xi+\varphi\zeta\right)d\mu.
\end{equation}
\end{definition}

For a fixed geodesic, there is at most one such field. Indeed,
the difference of two representatives belongs to $T_\mu\M$
and is orthogonal to every $(\nabla\varphi,4\varphi)$,
$\varphi\in C_c^\infty(\R^d)$, and therefore vanishes.

To characterize all regular barycentric logarithms between two measures, we introduce the following definition.

\begin{definition}[Semicoupling logarithmic map]
\label{def:semicoupling-logarithm}
Let $\bfgamma=(\gamma_0,\gamma_1)$ be a semicoupling between
$\mu$ and $\nu$, and let
$u\colon\R^d\times\R^d\to[0,1]$ be measurable.
Choose a finite positive measure $\lambda$ dominating
$\gamma_0$ and $\gamma_1$, and set
\begin{equation*}
    a\coloneqq\frac{d\gamma_0}{d\lambda},
    \qquad
    b\coloneqq\frac{d\gamma_1}{d\lambda}.
\end{equation*}
We define the
\emph{transport} and \emph{reaction} components,
$\gamma_0$-a.e., by
\begin{subequations}\label{eq:semicoupling-logarithmic-map}
\begin{align}
    \Log^{T}_{\bfgamma,u}(x,y)
    &\coloneqq
    \sqrt{\frac{b(x,y)}{a(x,y)}}\,
    \begin{cases}
        \mathbf{sin}(y-x),
            & \norm{y-x}<\pi/2,\\[2pt]
        \displaystyle u(x,y)\frac{y-x}{\norm{y-x}},
            & \norm{y-x}=\pi/2,\\[4pt]
        0,
            & \norm{y-x}>\pi/2,
    \end{cases}
    \label{eq:semicoupling-log-transport}
    \\[4pt]
    \Log^{R}_{\bfgamma}(x,y)
    &\coloneqq
    2\left(
        \sqrt{\frac{b(x,y)}{a(x,y)}}\,
        \cos\left(\norm{y-x}\wedge\frac{\pi}{2}\right)-1
    \right).
    \label{eq:semicoupling-log-reaction}
\end{align}
\end{subequations}
The \emph{semicoupling logarithmic map} is
\begin{equation*}
    \Log_{\bfgamma,u}
    \coloneqq
    \left(
        \Log^{T}_{\bfgamma,u},
        \Log^{R}_{\bfgamma}
    \right).
\end{equation*}

Disintegrating the first leg as
$\gamma_0(dx,dy)=\mu(dx)\kappa_x(dy)$, we define
the corresponding \emph{barycentric components} by
\begin{equation*}
\begin{aligned}
    \barlog^{T}_{\bfgamma,u}(x)
    \coloneqq
    \int_{\R^d}
        \Log^{T}_{\bfgamma,u}(x,y)
    \,\kappa_x(dy),\qquad
    \barlog^{R}_{\bfgamma}(x)
    \coloneqq
    \int_{\R^d}
        \Log^{R}_{\bfgamma}(x,y)
    \,\kappa_x(dy).
\end{aligned}
\end{equation*}
The \emph{barycentric logarithmic map} is then
\begin{equation*}
    \barlog_{\bfgamma,u}(x)
    \coloneqq
    \left(
        \barlog^{T}_{\bfgamma,u}(x),
        \barlog^{R}_{\bfgamma}(x)
    \right)
    =
    \int_{\R^d}
        \Log_{\bfgamma,u}(x,y)
    \,\kappa_x(dy).
\end{equation*}
\end{definition}

{\begin{remark}\label{rem:cone-log-normalization}
The semicoupling logarithm is the normalized version of the cone logarithm. For $r>0$ define
\begin{equation*}
    N_r(\xi,\zeta)
    \coloneqq
    \left(\xi,\frac{2\zeta}{r}\right).
\end{equation*}For every reduced cone pair $(\y_0,\y_1)$, we define the
normalized cone logarithm
\begin{equation*}
    \mathcal L_{\cone}(\y_0,\y_1)=(\mathcal L_{\cone}^T(\y_0,\y_1),\mathcal L_{\cone}^R(\y_0,\y_1))
    \coloneqq
    \begin{cases}
        N_r\bigl(\Log_{\cone}(\y_0,\y_1)\bigr),
        & \y_0=[x,r],\ r>0,\\
        (0,0), & \y_0=\vertex.
    \end{cases}
\end{equation*}
By the formulas of Section~\ref{subsec:wfr-geodesics}, for
$r,s>0$ and $\norm{y-x}\le\pi/2$, the normalized cone
logarithm satisfies
\begin{equation}\label{eq:normalized-cone-log-explicit}
\mathcal L_\cone([x,r],[y,s])
=
\left(
\frac{s}{r}{\bf sin}(y-x),
2\left(\frac{s}{r}\cos(\norm{y-x})-1\right)
\right),
\end{equation}
whereas
$\mathcal L_\cone([x,r],\vertex)=(0,-2)$.
Then, on $\{a>0,\ \norm{x-y}<\pi/2\}$,
\begin{equation}\label{eq:normalized-cone-log-semicoupling}
    \Log_{\bfgamma,u}(x,y)
    =
    \mathcal L_{\cone}
    \left([x,\sqrt{a(x,y)}],[y,\sqrt{b(x,y)}]\right),
\end{equation}
and the same identity holds at the cutoff $\norm{y-x}=\pi/2$
when $u(x,y)=1$.

The definition is independent of the dominating measure $\lambda$,
since the ratio $b/a$ is unchanged $\gamma_0$-a.e.; this follows by
comparing any two choices through a common dominating measure.
Moreover, the normalized cone logarithm is jointly $0$-homogeneous
in the radial variables:
\begin{equation*}
    \mathcal L_{\cone}([x,\theta r],[y,\theta s])
    =\mathcal L_{\cone}([x,r],[y,s]),
    \qquad \theta>0.
\end{equation*}
The dilation identities and the resulting normalization of finite
cone-geodesic measures are collected in
Lemma~\ref{lemma:cone-dilation-normalization}.

Finally, the definition is only needed $\gamma_0$-a.e. Since
$d\gamma_0=a\,d\lambda$, one has $\gamma_0(\{a=0\})=0$; on this set the
source cone point $[x,\sqrt a]$ is the vertex $\vertex$, and any value
assigned to the logarithm is therefore irrelevant.
\end{remark}}

\begin{remark}[On the role of the selection $u$]
Let $a,b>0$ and $x,y\in\R^d$ satisfy $\norm{y-x}=\pi/2$. By the cone-geodesic
formulas in Section~\ref{subsec:wfr-geodesics}, the curves
\begin{equation*}
    \bfmu^T_t\coloneqq m_t\delta_{x_t},
    \qquad
    \bfmu^R_t\coloneqq (1-t)^2a\delta_x+t^2b\delta_y,
\end{equation*}
where
\begin{equation*}
    m_t\coloneqq (1-t)^2a+t^2b,
    \qquad
    x_t\coloneqq x+\operatorname{Arg}
    \bigl((1-t)\sqrt a+\mathrm i t\sqrt b\bigr)\frac{(y-x)}{\norm{y-x}},
\end{equation*}
are constant-speed WFR geodesics between $\mu\coloneqq a\delta_x$
and $\nu\coloneqq b\delta_y$.
Indeed, they arise from the optimal cone couplings
\begin{equation*}
    A^T\coloneqq\delta_{([x,\sqrt a],[y,\sqrt b])},
    \qquad
    A^R\coloneqq\delta_{([x,\sqrt a],\vertex)}
                   +\delta_{(\vertex,[y,\sqrt b])},
\end{equation*}
both of cost $a+b=\WFR^2(\mu,\nu)$, by
Theorem~\ref{thm:wfr-geodesics-cone}.
The first curve is a traveling Dirac mass, whereas the second
consists of two static Dirac masses with changing weights.

For each $u\in[0,1]$, the coupling $uA^T+(1-u)A^R$ is optimal
as well, so the same theorem shows that
$\bfmu^u_t\coloneqq u\bfmu^T_t+(1-u)\bfmu^R_t$ is a WFR geodesic.
Differentiating smooth expectations at $t=0$ gives its regular
barycentric logarithm
\begin{equation*}
    \mathcal V^u(x)=\left(u\sqrt{\frac ba}\frac{(y-x)}{\norm{y-x}},-2\right).
\end{equation*}
Thus $u$ specifies the proportion assigned to the traveling
geodesic: it scales the transport component, while the reaction
component is unchanged. The measurable selection in
Definition~\ref{def:semicoupling-logarithm} makes this choice
separately at each pair on the cut locus.
\end{remark}

Barycentric semicoupling logarithms belong to $H(\mu)$,
and their squared norm is bounded by the semicoupling cost.

\begin{lemma}\label{lemma:energy-bound}
Let $\bfgamma\in\operatorname{SC}(\mu,\nu)$ and let
$u\colon\R^d\times\R^d\to[0,1]$ be measurable.
Then $\barlog_{\bfgamma,u}\in H(\mu)$
and
\begin{equation}\label{eq:energy-bound}
    \norm{\barlog_{\bfgamma,u}}_\mu^2
    \leq
    \int_{\R^d\times\R^d}
        \norm{\Log_{\bfgamma,u}(x,y)}_{\WFR}^2
    \,d\gamma_0(x,y)
    \leq \mathcal C(\bfgamma).
\end{equation}
\end{lemma}

\begin{proof}
Let $\lambda$ dominate $\gamma_0$ and $\gamma_1$, and write
$a=d\gamma_0/d\lambda$ and $b=d\gamma_1/d\lambda$.
On $\{a>0\}$, a direct computation gives
\begin{equation*}
    a\norm{\Log_{\bfgamma,u}(x,y)}_{\WFR}^2
=
\begin{cases}
a+b-2\sqrt{ab}\cos(\norm{x-y}),
    & \norm{x-y}<\pi/2,\\
a+u(x,y)^2b,
    & \norm{x-y}=\pi/2,\\
a,
    & \norm{x-y}>\pi/2.
\end{cases}
\end{equation*}
Since $0\leq u\leq1$, these expressions are bounded above by
$c(x,y,a,b)$.
Extending the logarithm by zero on $\{a=0\}$ and integrating
against $\lambda$ proves the second inequality in
\eqref{eq:energy-bound}.
The first follows from Jensen's inequality applied to the
disintegration $\gamma_0(dx,dy)=\mu(dx)\kappa_x(dy)$.
\end{proof}

The first inequality in \eqref{eq:energy-bound} reflects
barycentric averaging. The second reflects mass creation and,
at the cutoff, the attenuation of the transport component by $u$.
The following example isolates the first two effects and shows
that both inequalities can be strict even for an optimal semicoupling.

\begin{example}\label{ex:logarithm-information-loss}
Work in dimension one. Fix $0<\rho<\pi/2$, $R>\pi/2$,
and $m>0$, and set
\begin{equation*}
    \mu=\delta_0,
    \qquad
    \nu=\frac12\delta_{-\rho}+\frac12\delta_\rho+m\delta_R.
\end{equation*}
Consider the semicoupling
$\bfgamma^\star\coloneqq(\gamma_0^\star,\gamma_1^\star)$ with
\begin{equation*}
    \gamma_0^\star
    =\frac12\delta_{(0,-\rho)}+\frac12\delta_{(0,\rho)},
    \qquad
    \gamma_1^\star=\gamma_0^\star+m\delta_{(0,R)}.
\end{equation*}
This semicoupling is optimal. Indeed, let
$\bfgamma=(\gamma_0,\gamma_1)\in\operatorname{SC}(\mu,\nu)$
and set
\begin{equation*}
    a_\pm\coloneqq\gamma_0(\{(0,\pm \rho)\}),
    \qquad
    b_\pm\coloneqq\gamma_1(\{(0,\pm \rho)\}).
\end{equation*}
The atom at $R$ has zero interaction with the source,
while the marginal constraints yield
$a_-+a_+\leq1$ and $b_\pm\leq1/2$.
By \eqref{eq:total-mass-formualtion}, minimizing the cost
amounts to maximizing the interaction term
\begin{equation*}
    \int_{\R\times\R}
        \sqrt{a(x,y)b(x,y)}
        \cos\left(|y-x|\wedge\frac{\pi}{2}\right)
    \,d\lambda(x,y)
    =
    \cos(\rho)\left(\sqrt{a_-b_-}+\sqrt{a_+b_+}\right),
\end{equation*}
where $a=d\gamma_0/d\lambda$ and $b=d\gamma_1/d\lambda$
for a common dominating measure $\lambda$.
By Cauchy--Schwarz,
\begin{equation*}
    \cos(\rho)\left(\sqrt{a_-b_-}+\sqrt{a_+b_+}\right)
    \leq
    \cos(\rho)\sqrt{(a_-+a_+)(b_-+b_+)}
    \leq\cos(\rho),
\end{equation*}
with equality attained by $\bfgamma^\star$. Consequently,
\begin{equation*}
    \WFR^2(\mu,\nu)
    =\mathcal C(\bfgamma^\star)
    =2(1-\cos(\rho))+m.
\end{equation*}

The selection $u$ plays no role in this example, and
the logarithm at the two atoms of $\gamma_0^\star$ is
\begin{equation*}
    \Log_{\bfgamma^\star,u}(0,\pm \rho)
    =\bigl(\pm\sin(\rho),\,2(\cos(\rho)-1)\bigr).
\end{equation*}
Its barycenter is therefore
\begin{equation*}
    \barlog_{\bfgamma^\star,u}(0)
    =\bigl(0,\,2(\cos(\rho)-1)\bigr).
\end{equation*}
In particular,
\begin{equation*}
    \norm{\barlog_{\bfgamma^\star,u}}_\mu^2
    =(1-\cos(\rho))^2,
    \qquad
    \int_{\R\times\R}
        \norm{\Log_{\bfgamma^\star,u}}_{\WFR}^2
    \,d\gamma_0^\star
    =2(1-\cos(\rho)).
\end{equation*}
Thus the first inequality in \eqref{eq:energy-bound}
loses exactly $\sin^2(\rho)$ through the cancellation of
the two transport components. The second loses exactly
$m$, the cost of creating the atom at $R$, which is
invisible to the source-supported logarithm.
\end{example}

We can now state our main characterization theorem. For this purpose,
we introduce the barycentric field associated with an optimal
cone-geodesic representation.

Let $\bfmu$ be a constant-speed WFR geodesic connecting
$\mu$ to $\nu$, and let ${\boldsymbol{\Pi}}\in\prob(\Geo(\cone))$
be an optimal representation as in Theorem~\ref{thm:wfr-geodesics-cone}. Writing
$e_t(\bfy)=[\bfx_t,\bfr_t]$, define the finite measure
\begin{equation*}
    dQ^{\boldsymbol{\Pi}}(\bfy)\coloneqq \bfr_0^2\,d{\boldsymbol{\Pi}}(\bfy).
\end{equation*}
Since $(\mathsf x_0)\pf Q^{\boldsymbol{\Pi}}=\mu$, we may disintegrate
\begin{equation*}
    Q^{\boldsymbol{\Pi}}=\int_{\R^d}Q_x^{\boldsymbol{\Pi}}\,d\mu(x),
\end{equation*}
where $Q_x^{\boldsymbol{\Pi}}$ is a probability measure concentrated on
$\mathsf x_0^{-1}(\{x\})$ for $\mu$-a.e.\ $x$. We define the corresponding
barycentric field by
\begin{equation}\label{eq:cone-barycentric-field}
    \mathcal V^{\boldsymbol{\Pi}}(x)\coloneqq
    \int_{\Geo(\cone)}
        \mathcal L_{\cone}(\bfy_0,\bfy_1)
        \,dQ_x^{\boldsymbol{\Pi}}(\bfy).
\end{equation}

\begin{theorem}[Characterization of regular barycentric logarithms]
\label{thm:regular-barycentric-log-characterization}
Let $\mu,\nu\in\M$. The set of regular barycentric logarithms
between $\mu$ and $\nu$ is
\begin{equation}\label{eq:regular-logarithms-characterization}
    \left\{
        \barlog_{\bfgamma,u}:
        \bfgamma\in\operatorname{OptSC}(\mu,\nu),\
        u:\R^d\times\R^d\to[0,1]
        \text{ measurable}
    \right\}.
\end{equation}
Moreover, every constant-speed WFR geodesic connecting $\mu$ to
$\nu$ admits a unique regular barycentric logarithm. For every
optimal cone-geodesic representation $\boldsymbol{\Pi}$ of that geodesic, this
logarithm equals $\mathcal V^{\boldsymbol{\Pi}}$. In particular, $\mathcal V^{\boldsymbol{\Pi}}$ is independent,
up to $\mu$-a.e.\ equality, of the optimal representation.

Every regular barycentric logarithm $\mathcal V$ between $\mu$ and $\nu$
satisfies
\begin{equation*}
    \norm{\mathcal V}_\mu\leq\WFR(\mu,\nu).
\end{equation*}
\end{theorem}

\begin{proof}
Let 
$\bfgamma\in\operatorname{OptSC}(\mu,\nu)$ and let
$u:\R^d\times\R^d\to[0,1]$ be measurable.
By Proposition~\ref{propo:semicoupling-log-realization},
there exists an optimal cone-geodesic representation $\boldsymbol{\Pi}$ such that
\begin{equation*}
    \mathcal V^{\boldsymbol{\Pi}}=\barlog_{\bfgamma,u}.
\end{equation*}
Its homogeneous marginals form a constant-speed WFR geodesic
connecting $\mu$ to $\nu$, and
Proposition~\ref{propo:cone-barycentric-initial-derivative} gives the initial
derivative identity \eqref{eq:regular-barycentric-logarithm}.
By Theorem~\ref{propo:tangency}, the field
$\barlog_{\bfgamma,u}$ belongs to $T_\mu\M$.
Hence it is a regular barycentric logarithm.

Conversely, let $\mathcal V\in T_\mu\M$ be a regular barycentric logarithm
associated with a constant-speed WFR geodesic
$\bfmu$. Choose an optimal cone-geodesic
representation ${\boldsymbol{\Pi}}$ of $\bfmu$ using Theorem~\ref{thm:wfr-geodesics-cone}.
By Proposition~\ref{propo:cone-barycenter-semicoupling}, there exist
$\bfgamma\in\operatorname{OptSC}(\mu,\nu)$ and a measurable selection
$u\colon\R^d\times\R^d\to[0,1]$ such that
\begin{equation*}
    \mathcal V^{\boldsymbol{\Pi}}=\barlog_{\bfgamma,u}.
\end{equation*}
Theorem~\ref{propo:tangency} implies that
$\mathcal V^{\boldsymbol{\Pi}}\in T_\mu\M$. Moreover, by
Proposition~\ref{propo:cone-barycentric-initial-derivative}, both $\mathcal V$ and
$\mathcal V^{\boldsymbol{\Pi}}$ represent the same initial derivative. Therefore
\begin{equation*}
    \left\langle
        \mathcal V-\mathcal V ^{\boldsymbol{\Pi}},(\nabla\varphi,4\varphi)
    \right\rangle_\mu=0
    \qquad
    \text{for every }\varphi\in C_c^\infty(\R^d).
\end{equation*}
Since $\mathcal V-\mathcal V^{\boldsymbol{\Pi}}\in T_\mu\M$, the definition of the tangent space
gives $\mathcal V=\mathcal V^{\boldsymbol{\Pi}}$ in $H(\mu)$. This proves
\eqref{eq:regular-logarithms-characterization}.

The same argument applies to every constant-speed WFR geodesic
connecting $\mu$ to $\nu$ and to each of its optimal
cone-geodesic representations. It proves both uniqueness of
the associated regular barycentric logarithm and independence
of $\mathcal V^{\boldsymbol{\Pi}}$ from the representation.

Finally, the norm bound follows from
Proposition~\ref{propo:cone-barycentric-initial-derivative}.
\end{proof}

\subsection{Representation by cone geodesics}

We prove the representation results used in the preceding Theorem \ref{thm:regular-barycentric-log-characterization}, starting with the initial derivative along a WFR geodesic.

\begin{proposition}\label{propo:cone-barycentric-initial-derivative}
Let $\bfmu$ be a constant-speed WFR geodesic connecting
$\mu$ and $\nu$, and let
$\boldsymbol{\Pi}\in\prob(\Geo(\cone))$
be an optimal representation as in
Theorem \ref{thm:wfr-geodesics-cone}. Then
$\mathcal V^{\boldsymbol{\Pi}}
\coloneqq(\xi^{\boldsymbol{\Pi}},\zeta^{\boldsymbol{\Pi}})
\in H(\mu)$ satisfies
\begin{equation*}
    \norm{\mathcal V^{\boldsymbol{\Pi}}}_{\mu}^2
    \leq\WFR^2(\mu,\nu)
\end{equation*}
and, for every $\varphi\in C_c^\infty(\R^d)$,
\begin{equation}\label{eq:cone-barycentric-limit}
    \lim_{t\downarrow0}
    \frac{\int_{\R^d}\varphi\,d\bfmu_t
    -\int_{\R^d}\varphi\,d\mu}{t}
    =
    \int_{\R^d}
    \left(
        \nabla\varphi\cdot\xi^{\boldsymbol{\Pi}}
        +\varphi\zeta^{\boldsymbol{\Pi}}
    \right)d\mu.
\end{equation}
In particular, if
$\mathcal V^{\boldsymbol{\Pi}}\in T_\mu\M$,
then it is the regular barycentric logarithm associated
with $\bfmu$.
\end{proposition}

\begin{proof}

We first verify square integrability. By the cone
logarithmic energy identity \eqref{eq:constant-speed-norm}, on $\{\bfr_0>0\}$,
\begin{equation*}
    \bfr_0^2\|\mathcal L_\cone(\bfy_0,\bfy_1)\|_{\WFR}^2
=d_\cone^2(\bfy_0,\bfy_1).
\end{equation*}
Therefore,
\begin{equation*}
    \int_{\Geo(\cone)}
\|\mathcal L_\cone(\bfy_0,\bfy_1)\|_{\WFR}^2\,dQ^{\boldsymbol{\Pi}
}=
\int_{\{\bfr_0>0\}}
d_\cone^2(\bfy_0,\bfy_1)\,d\boldsymbol\Pi
\le\WFR^2(\mu,\nu).
\end{equation*}
Conditional Jensen's inequality proves that
$\mathcal V^{\boldsymbol\Pi}\in H(\mu)$ and gives
the stated norm bound.

Fix $\varphi\in C_c^\infty(\R^d)$ and write $\bfy_t=[\bfx_t,\bfr_t]$. By
\eqref{eq:dynamic-homogeneous-marginal-dual},
\begin{equation}\label{eq:lifted-limit}
    \frac{\int_{\R^d}\varphi\,d\bfmu_t
    -\int_{\R^d}\varphi\,d\mu}{t}
    =
    \int_{\Geo^*(\cone)}
    \frac{\bfr_t^2\varphi(\bfx_t)-\bfr_0^2\varphi(\bfx_0)}{t}
    \,d\boldsymbol{\Pi}(\bfy).
\end{equation}
By \eqref{eq:incremental-ratio-bound}, the absolute
value of the integrand is bounded by $B_\varphi (\bfr_0^2+\bfr_1^2)$. This bound is integrable: since $(e_0,e_1)\pf\boldsymbol{\Pi}\in\Gamma_\cone(\mu,\nu)$, 
\begin{equation*}
    \int(\bfr_0^2+\bfr_1^2)\,d\boldsymbol{\Pi}(\bfy)
    =
    \mu(\R^d)+\nu(\R^d)<\infty.
\end{equation*}

From Lemma \ref{lemma:cone-observable-estimates}, the pointwise limit of the difference
quotient is
\begin{equation*}
\bfr_0^2\left(
\nabla\varphi(\bfx_0)\cdot\mathcal L_\cone^T(\bfy_0,\bfy_1)
+\varphi(\bfx_0)\mathcal L_\cone^R(\bfy_0,\bfy_1)
\right)
\end{equation*}
on $\{\bfr_0>0\}$, and zero on $\{\bfr_0=0\}$.
Dominated convergence therefore yields
\begin{equation*}
    \begin{aligned}
\lim_{t\downarrow0}
\frac{\int_{\R^d}\varphi\,d\bfmu_t
-\int_{\R^d}\varphi\,d\mu}{t}
&=
\int_{\Geo(\cone)}
\bfr_0^2\left(
\nabla\varphi(\bfx_0)\cdot\mathcal L_\cone^T(\bfy_0,\bfy_1)
+\varphi(\bfx_0)\mathcal L_\cone^R(\bfy_0,\bfy_1)
\right)d\boldsymbol\Pi(\bfy)\\
&=
\int_{\R^d}
\left(
\nabla\varphi\cdot\xi^{\boldsymbol\Pi}
+\varphi\zeta^{\boldsymbol\Pi}
\right)d\mu,
\end{aligned}
\end{equation*}
where the last equality follows from disintegration.
\end{proof}

If $\boldsymbol\Pi$ is an optimal cone-geodesic
representation and $\mu=\h_0(\boldsymbol\Pi)$, its
barycentric field $\mathcal V^{\boldsymbol\Pi}$ satisfies
\begin{equation}\label{eq:cone-barycenter-testing}
\int_{\R^d}\psi(x)\mathcal V^{\boldsymbol\Pi}(x)\,d\mu(x)
=
\int_{\Geo^*(\cone)}
\psi(\bfx_0)\bfr_0^2\mathcal L_\cone(\bfy_0,\bfy_1)
\,d\boldsymbol\Pi(\bfy)
\end{equation}
for every bounded Borel function $\psi:\R^d\to\R$.
The integrals are vector-valued. This identity follows
directly from the source-weighted disintegration
defining $\mathcal V^{\boldsymbol\Pi}$.

We now compare the semicoupling construction with the
barycentric fields associated with optimal cone-geodesic
representations.

Both constructions below use the identity
\begin{equation}\label{eq:reduced-source-logarithm}
    r^2\mathcal L_\cone([x,r],[y,s])
    =\left(rs\mathbf{sin}(y-x),\,
        2rs\cos\left(\norm{y-x}\wedge\frac\pi2\right)-2r^2\right)
\end{equation}
for reduced cone pairs. It follows from
\eqref{eq:normalized-cone-log-explicit} when both radii are positive,
and also holds for annihilation and creation, with the source-weighted
expression set to zero when $r=0$. In particular, above the cutoff
reducedness gives $rs=0$.

\begin{proposition}[Realization of semicoupling logarithms]
\label{propo:semicoupling-log-realization}
Let $\mu,\nu\in\M$, let
$\bfgamma\in\operatorname{OptSC}(\mu,\nu)$, and let
$u\colon\R^d\times\R^d\to[0,1]$ be measurable.
There exists
\begin{equation*}
    \boldsymbol{\Pi}\in
    \prob(\Geo^*(\cone))
\end{equation*}
such that
\begin{equation*}
    (e_0,e_1)\pf\boldsymbol{\Pi}
    \in\operatorname{Opt\Gamma}_{\cone}(\mu,\nu)
\end{equation*}
and
\begin{equation}\label{eq:semicoupling-log-realized}
    \mathcal V^{\boldsymbol{\Pi}}
    =
    \barlog_{\bfgamma,u}
    \qquad\mu\text{-a.e.}
\end{equation}
In particular, $\barlog_{\bfgamma,u}$ satisfies
\eqref{eq:regular-barycentric-logarithm} along  $\bfmu_t\coloneqq\h_t(\boldsymbol{\Pi})$.
\end{proposition}

\begin{proof}
Choose a finite positive measure $\lambda$ dominating both
legs of $\bfgamma$, and write
\begin{equation*}
    \gamma_0=a\lambda,
    \qquad
    \gamma_1=b\lambda.
\end{equation*}
We fix finite nonnegative measurable versions of $a,b$.

For each $(x,y)$, define a finite measure on $[0,\infty)^2$ by
\begin{equation}\label{eq:realization-radius-kernel}
    \eta_{x,y}
    \coloneqq
    \begin{cases}
    \delta_{(\sqrt a,\sqrt b)},
    & \norm{x-y}<\pi/2,\\[1mm]
    u\,\delta_{(\sqrt a,\sqrt b)}
    +(1-u)\bigl[\delta_{(\sqrt a,0)}
    +\delta_{(0,\sqrt b)}\bigr],
    & \norm{x-y}=\pi/2,\\[1mm]
    \delta_{(\sqrt a,0)}+\delta_{(0,\sqrt b)},
    & \norm{x-y}>\pi/2,
    \end{cases}
\end{equation}
where $a,b,u$ are evaluated at $(x,y)$.
This is a measurable kernel whose total mass is finite.
Its radial moments satisfy
\begin{equation}\label{eq:realization-radius-moments}
    \int r^2\,d\eta_{x,y}(r,s)=a(x,y),
    \qquad
    \int s^2\,d\eta_{x,y}(r,s)=b(x,y).
\end{equation}

Define a finite measure $\widetilde A$ on
$\bigl(\R^d\times[0,\infty)\bigr)^2$ by
\begin{equation*}
    \widetilde A(dx,dr,dy,ds)
    \coloneqq
    \lambda(dx,dy)\eta_{x,y}(dr,ds),
\end{equation*}
and let $A\coloneqq(\mathfrak p^{\otimes 2})\pf\widetilde A$.
By \eqref{eq:realization-radius-moments}, its homogeneous marginals
are $\mu,\nu$, and it is concentrated on reduced pairs.
The mixed radial moment is
\begin{equation}\label{eq:realization-mixed-moment}
    m(x,y)\coloneqq\int rs\,d\eta_{x,y}(r,s)
    =\begin{cases}
        \sqrt{ab},&\norm{x-y}<\pi/2,\\
        u\sqrt{ab},&\norm{x-y}=\pi/2,\\
        0,&\norm{x-y}>\pi/2.
    \end{cases}
\end{equation}
Set $k(x,y)\coloneqq\cos(\norm{y-x}\wedge\pi/2)$.
Expanding the cone cost and integrating in the radial variables gives
\begin{equation}\label{eq:realization-cone-cost}
\begin{aligned}
    \int_{\cone^2}\dc^2(\y_0,\y_1)\,dA
    =\int_{\R^{2d}}(a+b-2mk)\,d\lambda=\mathcal C(\bfgamma)=\WFR^2(\mu,\nu),
\end{aligned}
\end{equation}
where we used \eqref{eq:realization-mixed-moment} and the fact that
$k=0$ at and above the cutoff. Thus
$A\in\operatorname{Opt\Gamma}_{\cone}(\mu,\nu)$.

Similarly, integrating \eqref{eq:reduced-source-logarithm} gives
\begin{equation}\label{eq:realization-source-moment}
\begin{aligned}
    \int r^2\mathcal L_\cone([x,r],[y,s])\,d\eta_{x,y}(r,s)
    =\bigl(m\mathbf{sin}(y-x),2(mk-a)\bigr)=a(x,y)\Log_{\bfgamma,u}(x,y),
\end{aligned}
\end{equation}
where the right-hand side is set to zero when $a=0$.
Indeed, below the cutoff $m=\sqrt{ab}$, at the cutoff
$m=u\sqrt{ab}$ and $k=0$, and above it $m=k=0$.
This also accounts for the radial annihilation and creation contributions.
Set $\widetilde{\boldsymbol{\Pi}}\coloneqq G\pf A$.
By \eqref{eq:realization-source-moment} and the disintegration of
$\gamma_0$, for every bounded Borel function $\psi$ we obtain
\begin{equation}\label{eq:unnormalized-equality}
\begin{aligned}
    \int_{\Geo(\cone)}
    \psi(\bfx_0)\bfr_0^2\mathcal L_\cone(\bfy_0,\bfy_1)
    \,d\widetilde{\boldsymbol{\Pi}}
    =\int_{\R^{2d}}\psi(x)a(x,y)
      \Log_{\bfgamma,u}(x,y)\,d\lambda=\int_{\R^d}\psi\,
      \barlog_{\bfgamma,u}\,d\mu.
\end{aligned}
\end{equation}
All integrals are finite, since $m\leq\sqrt{ab}$ and
$\int(a+b)\,d\lambda<\infty$.
By Lemma~\ref{lemma:cone-dilation-normalization}, we may normalize
$\widetilde{\boldsymbol{\Pi}}$ to a probability measure
$\boldsymbol{\Pi}$ without changing its homogeneous time marginals,
endpoint cost, or the source-weighted integrals in
\eqref{eq:unnormalized-equality}.
It remains concentrated on $\Geo^*(\cone)$, and its endpoint coupling
is optimal. Identity~\eqref{eq:cone-barycenter-testing} then proves
\eqref{eq:semicoupling-log-realized}.
Theorem~\ref{thm:wfr-geodesics-cone} and
Proposition~\ref{propo:cone-barycentric-initial-derivative} give the
asserted geodesic and initial derivative identity.
\end{proof}

Conversely, we recover an optimal semicoupling and a selection $u$
from an optimal cone-geodesic representation. The same radial
moment identities apply, but now optimality forces saturation
of the mixed-moment bound below the cutoff.

\begin{proposition}[Semicoupling representation of cone barycenters]
\label{propo:cone-barycenter-semicoupling}
Let $\mu,\nu\in\M$ and let
$\boldsymbol{\Pi}\in\prob(\Geo(\cone))$
satisfy
\begin{equation*}
    A\coloneqq(e_0,e_1)\pf\boldsymbol{\Pi}
    \in\operatorname{Opt\Gamma}_{\cone}(\mu,\nu).
\end{equation*}
Then there exist
$\bfgamma\in\operatorname{OptSC}(\mu,\nu)$ and a measurable
function $u\colon\R^d\times\R^d\to[0,1]$ such that
\begin{equation}\label{eq:cone-barycenter-semicoupling}
    \mathcal V^{\boldsymbol{\Pi}}
    =
    \barlog_{\bfgamma,u}
    \qquad\mu\text{-a.e.}
\end{equation}
\end{proposition}

\begin{proof}
The endpoint coupling $A$ is concentrated on reduced pairs.
We lift $A$ to a probability measure $\widetilde A\coloneqq (\mathfrak i _0^{\otimes 2})\pf A$ on $\bigl(\R^d\times[0,\infty)\bigr)^2$.

Let $\lambda$ be the $(x,y)$-marginal of $\widetilde A$.
Disintegrate
\begin{equation}\label{eq:A-disintegration}
    \widetilde A(dx,dr,dy,ds)
    =
    \lambda(dx,dy)\eta_{x,y}(dr,ds),
\end{equation}
where $\eta_{x,y}$ is a probability measure for
$\lambda$-a.e. $(x,y)$.
Define
\begin{equation}\label{eq:conditional-radius-moments}
    a\coloneqq\int r^2\,d\eta_{x,y},
    \qquad
    b\coloneqq\int s^2\,d\eta_{x,y},
    \qquad
    m\coloneqq\int rs\,d\eta_{x,y}.
\end{equation}
These functions are measurable and finite $\lambda$-a.e.,
since
\begin{equation*}
    \int(a+b)\,d\lambda
    =
    \mu(\R^d)+\nu(\R^d)<\infty.
\end{equation*}
Cauchy--Schwarz in $L^2(\eta_{x,y})$ gives
\begin{equation}\label{eq:conditional-mixed-moment-bound}
    0\leq m\leq\sqrt{ab}\qquad\lambda\text{-a.e.}
\end{equation}
Set $\gamma_0\coloneqq a\lambda$ and $\gamma_1\coloneqq b\lambda$.
The homogeneous marginal constraints on $A$ show that
$\bfgamma=(\gamma_0,\gamma_1)\in\operatorname{SC}(\mu,\nu)$.

For brevity, write
\begin{equation*}
    k(x,y)\coloneqq
    \cos\left(\norm{y-x}\wedge \pi/2\right).
\end{equation*}
Since $A$ is reduced, its coordinate lift satisfies
$rs=0$ almost everywhere on $\{\norm{y-x}>\pi/2\}$.
Consequently, $m=0$ $\lambda$-a.e. on $\{\norm{y-x}>\pi/2\}$
and, using \eqref{eq:A-disintegration}, its optimal cone cost can be written as
\begin{equation*}
    \WFR^2(\mu,\nu)
    =
    \int(a+b-2mk)\,d\lambda.
\end{equation*}
Using \eqref{eq:conditional-mixed-moment-bound} and the
semicoupling formulation of WFR, we obtain
\begin{align*}
    \WFR^2(\mu,\nu)
    \leq\mathcal C(\bfgamma)
    =\int(a+b-2\sqrt{ab}\,k)\,d\lambda
    \leq\int(a+b-2mk)\,d\lambda
    =\WFR^2(\mu,\nu).
\end{align*}
Thus $\bfgamma\in\operatorname{OptSC}(\mu,\nu)$ and
\begin{equation*}
    \int k(\sqrt{ab}-m)\,d\lambda=0.
\end{equation*}
The integrand is nonnegative. Since $k>0$ below the cutoff,
we conclude that
\begin{equation}\label{eq:conditional-moment-saturation}
    m=\sqrt{ab}
    \qquad\lambda\text{-a.e. on }\{\norm{y-x}<\pi/2\}.
\end{equation}

Define the measurable selection
\begin{equation}\label{eq:selection-from-radius-moments}
    u(x,y)\coloneqq
    \begin{cases}
    \displaystyle\frac{m(x,y)}{\sqrt{a(x,y)b(x,y)}},
    & \norm{y-x}=\pi/2,\ a(x,y)b(x,y)>0,\\[2mm]
    0,&\text{otherwise}.
    \end{cases}
\end{equation}
By \eqref{eq:conditional-mixed-moment-bound}, it takes values
in $[0,1]$.

For $a(x,y)>0$, integrating
\eqref{eq:reduced-source-logarithm} gives
\begin{equation}\label{eq:conditional-source-log}
\begin{aligned}
    \ell(x,y)&\coloneqq\frac1{a(x,y)}
       \int r^2\mathcal L_\cone([x,r],[y,s])\,d\eta_{x,y}(r,s)\\
    &=\left(\frac ma\mathbf{sin}(y-x),\,
             2\left(\frac ma k-1\right)\right).
\end{aligned}
\end{equation}
Below the cutoff, \eqref{eq:conditional-moment-saturation}
identifies this with $\Log_{\bfgamma,u}$.
At the cutoff, the same follows from
\eqref{eq:selection-from-radius-moments}; if $b=0$, then $m=0$
and both expressions equal $(0,-2)$. Above the cutoff, $m=k=0$,
so they again equal $(0,-2)$. Thus
$\ell=\Log_{\bfgamma,u}$ $\gamma_0$-a.e.
Set $\ell=0$ on $\{a=0\}$.

By \eqref{eq:cone-barycenter-testing} and the disintegrations of
$\widetilde A$ and $\gamma_0$, for every bounded Borel function $\psi$,
\begin{equation*}
\begin{aligned}
    \int_{\R^d}\psi\,\mathcal V^{\boldsymbol{\Pi}}\,d\mu
    &=\int_{\R^{2d}}\psi(x)a(x,y)\ell(x,y)\,d\lambda\\
    &=\int_{\R^d}\psi\,
      \barlog_{\bfgamma,u}\,d\mu
\end{aligned}
\end{equation*}
and all integrals are finite by the logarithmic energy bounds.
This proves \eqref{eq:cone-barycenter-semicoupling}.
\end{proof}

\subsection{Endpoint tangency}

\begin{theorem}[Tangency of optimal semicoupling logarithms]
\label{propo:tangency}
Let $\mu,\nu\in\M$,
$\bfgamma\in\operatorname{OptSC}(\mu,\nu)$, and let
$u\colon\R^d\times\R^d\to[0,1]$ be measurable.
Then
\begin{equation*}
    \barlog_{\bfgamma,u}\in T_\mu\M.
\end{equation*}
\end{theorem}
\begin{proof}
Set $\mathcal L\coloneqq\barlog_{\bfgamma,u}$.
By Proposition~\ref{propo:semicoupling-log-realization}, there exists
$\boldsymbol{\Pi}\in\prob(\Geo^*(\cone))$ such that
\begin{equation*}
    A\coloneqq(e_0,e_1)\pf\boldsymbol{\Pi}
    \in\operatorname{Opt\Gamma}_{\cone}(\mu,\nu),
    \qquad
    \mathcal V^{\boldsymbol{\Pi}}=\mathcal L.
\end{equation*}
Proposition~\ref{propo:cone-barycentric-initial-derivative} gives
$\mathcal L\in H(\mu)$. Moreover, \eqref{eq:cone-barycenter-testing} identifies
$\mathcal L$ as the source-weighted barycenter of the normalized
cone logarithm $\mathcal L_\cone$.

Let $\mathcal W=(w,\eta)\in T_\mu\M^\perp$. Thus
\begin{equation}\label{eq:normal-space}
    \int_{\R^d}
    \left(w\cdot\nabla\varphi+\eta\varphi\right)\,d\mu=0
    \qquad
    \text{for every }\varphi\in C_c^\infty(\R^d).
\end{equation}
It suffices to prove $\langle\mathcal L,\mathcal W\rangle_\mu=0$.

The orthogonality relation \eqref{eq:normal-space} means that
$\mathcal W$ generates a stationary solution of the inhomogeneous
continuity equation. After restricting to a bounded domain, we use
the superposition principle of~\cite{bredies-carion-fanzon2022} to
represent this solution by characteristics and choose an interior
time $\tau$ at which almost every characteristic with positive
radius has the prescribed derivatives.
We then normalize the positive first radii of $A$ to one and glue
its source coordinate to the characteristics at time $\tau$.
Moving along these characteristics perturbs the first cone endpoint;
a correction accounts for paths whose radius vanishes at $\tau$
and restores the first homogeneous marginal exactly. This yields
competitors $A_\varepsilon\in\Gamma_\cone(\mu,\nu)$ for both signs
of $\varepsilon$, with $A_0=A$ after normalization.

Optimality gives a local minimum of the competitor cost at zero.
To differentiate it, we use the stationary energy identity to obtain
a uniform $L^2$ bound on the interaction difference quotients,
then truncate the possibly unbounded target radius. Finally, we
identify the derivative with
$-2\langle\mathcal L,\mathcal W\rangle_\mu$, which proves the
required orthogonality.

\emph{Restriction to bounded domains}: By
Lemma~\ref{lemma:compact-support-approximation} and continuity
of the pairing in $H(\mu)$, we may assume that both components
of $\mathcal W$ vanish $\mu$-a.e.\ outside $B_R$, for some $R>0$.
We therefore choose finite-valued Borel representatives that vanish
everywhere outside $B_R$.

Let $\Omega=\overline B_{2R}$ and
$\widetilde\mu\coloneqq\mu\llcorner\Omega$.
Since
\begin{equation*}
    w\widetilde\mu=w\mu,
    \qquad
    \eta\widetilde\mu=\eta\mu,
\end{equation*}
equation~\eqref{eq:normal-space} gives
\begin{equation}\label{eq:localized-equation}
    \operatorname{div}(w\widetilde\mu)
    =\eta\widetilde\mu
    \qquad\text{in }\mathcal D'(\R^d).
\end{equation}
Consequently, the constant curve
$\widetilde{\bfmu}_t\equiv\widetilde\mu$ solves the
inhomogeneous continuity equation on $\Omega$.
The equation holds with no boundary flux, since both
coefficients vanish on a neighborhood of $\partial\Omega$.
Testing \eqref{eq:normal-space} with a smooth function equal
to one on $\overline B_R$ also gives
\begin{equation}\label{eq:normal-reaction-mean}
    \int_{\R^d}\eta\,d\mu
    =
    \int_{\Omega}\eta\,d\widetilde\mu
    =0.
\end{equation}
If $\widetilde\mu=0$, then $\mathcal W=0$ in $H(\mu)$ and
there is nothing to prove. Henceforth assume $\widetilde\mu\neq0$.

Let
\begin{equation*}
    \mathcal A_{\Omega}
    \coloneqq
    \{r^2\delta_x:r\geq0,\ x\in\Omega\}
    \subset \mathcal M^+(\Omega)\quad\text{and}\quad \mathbf C_{\Omega}
\coloneqq C([0,1];\mathcal A_{\Omega}),
\end{equation*}
where $\mathcal A_{\Omega}$ is endowed with the
bounded-Lipschitz distance and $\mathbf C_{\Omega}$ with
the corresponding uniform distance. For
$\bfomega\in\mathbf C_{\Omega}$, write
\begin{equation*}
    \bfomega_t
    =
    \bfr_t(\bfomega)^2\delta_{\bfx_t(\bfomega)},
\end{equation*}
where $\bfr_t(\bfomega)\geq0$,
$\bfx_t(\bfomega)\in\Omega$ whenever
$\bfr_t(\bfomega)>0$, and we set $\bfx_t(\bfomega)=0$ when
$\bfr_t(\bfomega)=0$, consistently with the choice of
$\mathfrak i_0$ as a right inverse of $\mathfrak p$.

We further set
\begin{equation*}
    \mathbf S_{\Omega}
    \coloneqq
    \left\{
        \bfomega\in\mathbf C_{\Omega}:
        \bfr_\cdot(\bfomega)\in AC^2([0,1]),\;
        \bfr_\cdot(\bfomega)\bfx_\cdot(\bfomega)
        \in AC^2([0,1];\R^d)
    \right\}.
\end{equation*}
For $\bfomega\in\mathbf S_{\Omega}$, consider the
characteristic equations
\begin{equation}\label{eq:stationary-characteristics}
    \dot\bfx_t(\bfomega)=w(\bfx_t(\bfomega)),
    \qquad
    \dot\bfr_t(\bfomega)
    =
    \frac12\eta(\bfx_t(\bfomega))\bfr_t(\bfomega),
\end{equation}
and define
\begin{equation*}
    \mathbf S_{\Omega}^{w,\eta}
    \coloneqq
    \left\{
        \bfomega\in\mathbf S_{\Omega}:
        \eqref{eq:stationary-characteristics}
        \text{ holds for a.e. }t\in\{\bfr_\cdot(\bfomega)>0\}
    \right\},
\end{equation*}
together with
\begin{equation*}
    \mathcal S
    \coloneqq
    \left\{
        \bfomega\in\mathbf S_{\Omega}^{w,\eta}:
        \int_0^1\bfr_t(\bfomega)^2\,dt=1
    \right\}.
\end{equation*}
By \cite[Theorem~4.3 and definition~(32)]{bredies-carion-fanzon2022},
there exists a finite positive Borel measure $\sigma$ on
$\mathcal S$ such that
\begin{equation}\label{eq:stationary-superposition}
    \widetilde\mu
    =
    \int_{\mathcal S}
        \bfr_t(\bfomega)^2
        \delta_{\bfx_t(\bfomega)}
    \,d\sigma(\bfomega)
    \qquad\text{for every }t\in[0,1].
\end{equation}
Moreover, for $\sigma$-a.e.\ $\bfomega$,
$\bfx_\cdot(\bfomega)$ is locally absolutely continuous on
$\{\bfr_\cdot(\bfomega)>0\}$.

Let $E$ be the set of pairs
$(t,\bfomega)\in(0,1)\times\mathcal S$ with
$\bfr_t(\bfomega)>0$ for which at least one of the two derivatives
in \eqref{eq:stationary-characteristics} does not exist or does
not have the prescribed value. To verify measurability, set
\begin{equation*}
\begin{aligned}
    Z(t,\bfomega)
    &\coloneqq(\bfx_t(\bfomega),\bfr_t(\bfomega)),\\
    G(t,\bfomega)
    &\coloneqq\left(w(\bfx_t(\bfomega)),
        \tfrac12\eta(\bfx_t(\bfomega))\bfr_t(\bfomega)\right).
\end{aligned}
\end{equation*}
The radius map is continuous, and the spatial coordinate map is
continuous where the radius is positive and Borel under our
convention at zero. Thus $Z$ and $G$ are jointly Borel measurable.
For each nonzero rational $h$, the difference quotient
\begin{equation*}
    \frac{Z(t+h,\bfomega)-Z(t,\bfomega)}{h}
\end{equation*}
is Borel on $\{(t,\bfomega):0<t,t+h<1\}$. Convergence of these
quotients to $G(t,\bfomega)$ as rational $h\to0$ is a Borel
condition: for every $k\geq1$, there must exist $n\geq1$ such that
the error is at most $1/k$ for every admissible rational $h$ with
$0<|h|<1/n$. At a point with $\bfr_t>0$, both coordinate paths
are continuous in a neighborhood of $t$, so this rational limit
exists with the prescribed value if and only if the ordinary
derivative does. Hence $E$ is Borel measurable.
For $\sigma$-a.e.\ $\bfomega$, its time section
is Lebesgue-null. Therefore, Tonelli's theorem gives
\begin{equation*}
    0
    =
    \int_{\mathcal S}\int_0^1
    \mathbf 1_E(t,\bfomega)\,dt\,d\sigma(\bfomega)
    =
    \int_0^1\sigma(E_t)\,dt,
\end{equation*}
where
\begin{equation*}
    E_t\coloneqq\{\bfomega\in\mathcal S\colon (t,\bfomega)\in E\}.
\end{equation*}
Thus $\sigma(E_t)=0$ for a.e.\ $t\in(0,1)$.

For $\bfomega\in\mathcal S$, define the Borel energy density
\begin{equation*}
    \ell_t(\bfomega)
    \coloneqq
    \bfr_t(\bfomega)^2
    \left(
        \norm{w(\bfx_t(\bfomega))}^2
        +\frac14|\eta(\bfx_t(\bfomega))|^2
    \right).
\end{equation*}
By \eqref{eq:stationary-characteristics}, this agrees a.e.\ in time
on $\{\bfr_t>0\}$ with the cone energy
$|\dot\bfr_t|^2+\bfr_t^2\norm{w(\bfx_t)}^2$.
The same holds on $\{\bfr_t=0\}$, since an absolutely continuous
function has derivative zero a.e.\ on each of its level sets.
By \eqref{eq:stationary-superposition},
\begin{equation}\label{eq:stationary-energy}
    \int_{\mathcal S}\ell_t\,d\sigma
    =\int_{\Omega}
      \left(\norm{w}^2+\frac14|\eta|^2\right)\,d\widetilde\mu
    =\norm{\mathcal W}_\mu^2
    \eqqcolon\mathcal E
    \qquad(t\in[0,1]).
\end{equation}

Choose $\tau\in(0,1)$ such that $\sigma(E_\tau)=0$.
Then, for $\sigma$-a.e.\ $\bfomega$ with $\bfr_\tau(\bfomega)>0$,
both coordinate paths are differentiable at $\tau$ and satisfy
\eqref{eq:stationary-characteristics}. This gives the pointwise
derivatives needed below; the energy identity will provide a
uniform $L^2$ bound on the corresponding difference quotients.

\emph{Building the competitors $A_\varepsilon$}: Write
$\y_0=[x,r]$ and $\y_1=[y,s]$. We may assume without loss of generality that $r=1$ for $A$-a.e.\ pair with $r>0$. Indeed, consider the dilation
\begin{equation*}
    \vartheta(\y_0,\y_1)
    \coloneqq
    \begin{cases}
        r,&r>0,\\
        1,&r=0,
    \end{cases}
\end{equation*}
which is admissible, since
\begin{equation*}
    \int\vartheta^2\,dA
    =\mu(\R^d)+A(\{r=0\})<\infty.
\end{equation*}
It sends $([x,r],[y,s])$ to $([x,1],[y,s/r])$ on
$\{r>0\}$ and leaves the source-vertex part unchanged, and by
Lemma~\ref{lemma:cone-dilation-normalization}, it preserves
both homogeneous marginals, the optimal cost, reducedness, and
the source-weighted barycenter $\mathcal L$. 

For $-\tau<\varepsilon<1-\tau$, we perturb the first endpoint
along the characteristics \eqref{eq:stationary-characteristics}.
Set
\begin{equation*}
    \cone_{\Omega}\coloneqq
    \{(\y_0,\y_1)\in\cone^2:
      \y_0=[x,r],\ r>0,\ x\in\Omega\},
    \qquad
    A_\Omega\coloneqq A\llcorner\cone_{\Omega}.
\end{equation*}
Since its first radius equals one,
$(\mathsf x)\pf A_\Omega=\widetilde\mu$. Disintegrate
\begin{equation*}
    A_\Omega(d\y_0,d\y_1)
    =\widetilde\mu(dx)Q_x(d\y_0,d\y_1).
\end{equation*}
For $\widetilde\mu$-a.e.\ $x$, the kernel $Q_x$ is a probability
measure concentrated on pairs with $\y_0=[x,1]$.
On $\mathcal S\times\cone^2$, define
\begin{equation}\label{eq:tangency-gluing}
    \Lambda(d\bfomega,d\y_0,d\y_1)
    \coloneqq
    \bfr_\tau(\bfomega)^2\,\sigma(d\bfomega)
    Q_{\bfx_\tau(\bfomega)}(d\y_0,d\y_1).
\end{equation}
This measure is concentrated on triples satisfying
\begin{equation*}
    \y_0=[\bfx_\tau(\bfomega),1],
    \qquad \bfr_\tau(\bfomega)>0,
\end{equation*}
and its marginals on $\mathcal S$ and $\cone^2$ are respectively
$\bfr_\tau^2\sigma$ and $A_\Omega$. In particular,
\begin{equation}\label{eq:tangency-gluing-moments}
    \Lambda(\mathcal S\times\cone^2)=\widetilde\mu(\Omega)<\infty,
    \qquad
    \int s^2\,d\Lambda
    =\int s^2\,dA_\Omega\leq\nu(\R^d).
\end{equation}
For $-\tau<\varepsilon<1-\tau$, set
\begin{equation*}
    a_\varepsilon(\bfomega)
    \coloneqq
    \frac{\bfr_{\tau+\varepsilon}(\bfomega)}
         {\bfr_\tau(\bfomega)}
    \qquad\text{on }\{\bfr_\tau>0\},
\end{equation*}
and define, $\Lambda$-a.e.,
\begin{equation*}
    F_\varepsilon(\bfomega,\y_0,\y_1)
    \coloneqq
    \left(
        [\bfx_{\tau+\varepsilon}(\bfomega), a_\varepsilon(\bfomega)],
        \y_1
    \right).
\end{equation*}
Then $F_\varepsilon$ leaves the second cone endpoint fixed and moves
the first one according to the characteristic indexed by $\bfomega$.

By \eqref{eq:tangency-gluing-moments}, $\Lambda$ is finite,
so we may push it forward to a finite measure on $\cone^2$ through $F_\varepsilon$. Its first homogeneous marginal is characterized,
for every bounded measurable $\varphi$, by
\begin{equation*}
\begin{aligned}
    \int \varphi(z)\,
    d\h_1\bigl((F_\varepsilon)\pf\Lambda\bigr)(z)
    &=
    \int_{\mathcal S\times\cone^2}
    a_\varepsilon(\bfomega)^2
    \varphi(\bfx_{\tau+\varepsilon}(\bfomega))
    \,d\Lambda\\
    &=
    \int_{\{\bfr_\tau>0\}}
    \bfr_{\tau+\varepsilon}(\bfomega)^2
    \varphi(\bfx_{\tau+\varepsilon}(\bfomega))
    \,d\sigma(\bfomega).
\end{aligned}
\end{equation*}
Therefore, by \eqref{eq:stationary-superposition},
\begin{equation}\label{eq:correction-term}
    \h_1\bigl((F_\varepsilon)\pf\Lambda\bigr)
    =
    \widetilde\mu-\beta_\varepsilon,
\end{equation}
where
\begin{equation*}
    \beta_\varepsilon
    \coloneqq
    \int_{\{\bfr_\tau=0\}}
    \bfr_{\tau+\varepsilon}(\bfomega)^2
    \delta_{\bfx_{\tau+\varepsilon}(\bfomega)}
    \,d\sigma(\bfomega)\in\mathcal M^+(\Omega)
\end{equation*}
is the contribution at time $\tau+\varepsilon$ of characteristics
whose radius vanishes at time $\tau$.

The second homogeneous marginal is unchanged. Indeed, for every
bounded measurable $\varphi$,
\begin{equation*}
\begin{aligned}
    \int \varphi(z)\,
    d\h_2\bigl((F_\varepsilon)\pf\Lambda\bigr)(z)
    &=
    \int_{\mathcal S\times\cone^2}
    s^2\varphi(y)\,d\Lambda\\
    &=
    \int_{\cone^2}
    s^2\varphi(y)\,dA_\Omega.
\end{aligned}
\end{equation*}
Hence
\begin{equation*}
    \h_2\bigl((F_\varepsilon)\pf\Lambda\bigr)
    =
    \h_2(A_\Omega).
\end{equation*}

Define
\begin{equation*}
    \mathfrak c(z)\coloneqq([z,1],\vertex).
\end{equation*}
Since
\begin{equation*}
    \h_1(\mathfrak c\pf\beta_\varepsilon)=\beta_\varepsilon,
    \qquad
    \h_2(\mathfrak c\pf\beta_\varepsilon)=0,
\end{equation*}
set
\begin{equation}\label{eq:tangency-competitor}
    A_\varepsilon
    \coloneqq(F_\varepsilon)\pf\Lambda
    +(A-A_\Omega)
    +\mathfrak c\pf\beta_\varepsilon.
\end{equation}
Using \eqref{eq:correction-term} and the preceding identity for the
second homogeneous marginal, we obtain
\begin{equation*}
\begin{aligned}
    \h_1(A_\varepsilon)
    &=
    (\widetilde\mu-\beta_\varepsilon)
    +\h_1(A-A_\Omega)
    +\beta_\varepsilon
    =\mu,\\
    \h_2(A_\varepsilon)
    &=
    \h_2(A_\Omega)
    +\h_2(A-A_\Omega)
    =\nu.
\end{aligned}
\end{equation*}
Moreover,
\begin{equation*}
    \int r^2\,dA_\varepsilon=\mu(\R^d),
    \qquad
    \int s^2\,dA_\varepsilon=\nu(\R^d),
\end{equation*}
so that $A_\varepsilon\in\Gamma_\cone(\mu,\nu)$.

Notice also that $\beta_0=0$ and $F_0$ leaves both cone endpoints
unchanged on the support of $\Lambda$. Since the $\cone^2$-marginal
of $\Lambda$ is $A_\Omega$,
\begin{equation*}
    (F_0)\pf\Lambda=A_\Omega.
\end{equation*}
Consequently $A_0=A$.

\emph{Differentiation}: Set
\begin{equation*}
    C(\varepsilon)
    \coloneqq
    \int_{\cone^2}\dc^2\,dA_\varepsilon.
\end{equation*}
Since $A_\varepsilon\in\Gamma_\cone(\mu,\nu)$ and $A=A_0$ is
optimal,
\begin{equation}\label{eq:tangency-cost-minimum}
    C(\varepsilon)
    \geq C(0)
    =
    \WFR^2(\mu,\nu)
\end{equation}
for every $-\tau<\varepsilon<1-\tau$.

Expanding the cone cost \eqref{eq:distance} and using the fixed homogeneous marginals of $A_\varepsilon$, we obtain
\begin{equation}\label{eq:tangency-cost-difference}
    C(\varepsilon)-C(0)
    =
    -2\int_{\mathcal S\times\cone^2} s
    \left[
        a_\varepsilon K_y(\bfx_{\tau+\varepsilon})
        -K_y(\bfx_\tau)
    \right]\,d\Lambda.
\end{equation}
Indeed, $A-A_\Omega$ is unchanged, while the correction term has
target radius zero and hence gives no contribution to the interaction
term. On the moving part, the interaction term is
\begin{equation*}
    -2\int_{\mathcal S\times\cone^2} sa_\varepsilon
    K_y(\bfx_{\tau+\varepsilon})\,d\Lambda,
\end{equation*}
which yields \eqref{eq:tangency-cost-difference} after subtracting
the value at $\varepsilon=0$.

We now justify differentiation in
\eqref{eq:tangency-cost-difference}. For fixed $y$, set
\begin{equation*}
    \psi_y(t,\bfomega)
    \coloneqq
    \bfr_t(\bfomega)K_y(\bfx_t(\bfomega)).
\end{equation*}
The function $\psi_y(\cdot,\bfomega)$ is absolutely continuous, also
across times at which $\bfr_t=0$. Indeed, using
$|K_y|\leq1$, $\operatorname{Lip}(K_y)\leq1$, and
$\norm{\bfx_t}\leq2R$, we have
\begin{equation*}
\begin{aligned}
    |\psi_y(t_2,\bfomega)-\psi_y(t_1,\bfomega)|
    &\leq
    |\bfr_{t_2}-\bfr_{t_1}|
    +
    \bfr_{t_1}\norm{\bfx_{t_2}-\bfx_{t_1}}\\
    &\leq
    (1+2R)|\bfr_{t_2}-\bfr_{t_1}|
    +
    \norm{
        \bfr_{t_2}\bfx_{t_2}-\bfr_{t_1}\bfx_{t_1}
    }.
\end{aligned}
\end{equation*}
Since both $\bfr_\cdot$ and
$\bfr_\cdot\bfx_\cdot$ are absolutely continuous, the same holds
for $\psi_y$.

On $\{\bfr_t>0\}$, the chain rule and
\eqref{eq:stationary-characteristics} give
\begin{equation*}
    \partial_t \psi_y
    =
    \dot\bfr_tK_y(\bfx_t)
    +
    \bfr_t\nabla K_y(\bfx_t)\cdot w(\bfx_t).
\end{equation*}
Hence, by Cauchy--Schwarz and
\eqref{eq:kernel-bounds},
\begin{equation*}
\begin{aligned}
    |\partial_t \psi_y|
    &\leq
    \left(
        |\dot\bfr_t|^2
        +
        \bfr_t^2\norm{w(\bfx_t)}^2
    \right)^{1/2}
    \left(
        |K_y(\bfx_t)|^2
        +
        \norm{\nabla K_y(\bfx_t)}^2
    \right)^{1/2}=\sqrt{\ell_t}.
\end{aligned}
\end{equation*}
On $\{\bfr_t=0\}$, $\psi_y=0$, and therefore
$\partial_t \psi_y=0$ for a.e.\ such $t$. Consequently,
\begin{equation}\label{eq:tangency-observable-increment}
    |\psi_y(t_2,\bfomega)-\psi_y(t_1,\bfomega)|
    \leq
    \int_{\min\{t_1,t_2\}}^{\max\{t_1,t_2\}}
    \sqrt{\ell_\theta(\bfomega)}\,d\theta.
\end{equation}

For $\varepsilon\neq0$, define
\begin{equation*}
    U_\varepsilon
    \coloneqq
    \frac{
        a_\varepsilon K_y(\bfx_{\tau+\varepsilon})
        -K_y(\bfx_\tau)
    }{\varepsilon}.
\end{equation*}
Since $\Lambda$ is concentrated on $\{\bfr_\tau>0\}$,
\begin{equation*}
    U_\varepsilon
    =
    \frac{
        \psi_y(\tau+\varepsilon,\bfomega)
        -
        \psi_y(\tau,\bfomega)
    }{
        \varepsilon\bfr_\tau(\bfomega)
    }
    \qquad\Lambda\text{-a.e.}
\end{equation*}
By the choice of $\tau$, for $\Lambda$-a.e.\ triple both
$\bfx_\cdot$ and $\bfr_\cdot$ are differentiable at $\tau$ and
satisfy \eqref{eq:stationary-characteristics}. Therefore
\begin{equation}\label{eq:tangency-quotient-limit}
    U_\varepsilon
    \longrightarrow
    U_0
    \coloneqq
    \nabla K_y(\bfx_\tau)\cdot w(\bfx_\tau)
    +
    \frac12K_y(\bfx_\tau)\eta(\bfx_\tau)
    \qquad\Lambda\text{-a.e.}
\end{equation}
The exceptional set is independent of $y$, since once
$\bfx_\cdot$ and $\bfr_\cdot$ are differentiable at $\tau$, the
ordinary chain rule applies to every $C^1$ function $K_y$.

Let $I_\varepsilon$ be the interval with endpoints
$\tau$ and $\tau+\varepsilon$. By
\eqref{eq:tangency-observable-increment},
\begin{equation*}
    |U_\varepsilon|
    \leq
    \frac1{|\varepsilon|\bfr_\tau}
    \int_{I_\varepsilon}\sqrt{\ell_t}\,dt.
\end{equation*}
Cauchy--Schwarz in time gives
\begin{equation*}
    |U_\varepsilon|^2
    \leq
    \frac1{|\varepsilon|^2\bfr_\tau^2}
    \left(\int_{I_\varepsilon}\sqrt{\ell_t}\,dt\right)^2
    \leq
    \frac1{|\varepsilon|\bfr_\tau^2}
    \int_{I_\varepsilon}\ell_t\,dt.
\end{equation*}
Integrating with respect to $\Lambda$, the factor
$\bfr_\tau^2$ in \eqref{eq:tangency-gluing} cancels the denominator.
Since $Q_x$ is a probability measure for
$\widetilde\mu$-a.e.\ $x$, we obtain
\begin{equation*}
\begin{aligned}
    \int |U_\varepsilon|^2\,d\Lambda
    &\leq
    \frac1{|\varepsilon|}
    \int_{I_\varepsilon}
    \int_{\{\bfr_\tau>0\}}
    \ell_t(\bfomega)\,d\sigma(\bfomega)\,dt\\
    &\leq
    \frac1{|\varepsilon|}
    \int_{I_\varepsilon}
    \int_{\mathcal S}
    \ell_t(\bfomega)\,d\sigma(\bfomega)\,dt\\
    &\leq
    \mathcal E,
\end{aligned}
\end{equation*}
where we used \eqref{eq:stationary-energy}. Thus
\begin{equation}\label{eq:tangency-uniform-quotient-bound}
    \sup_{\varepsilon}
    \norm{U_\varepsilon}_{L^2(\Lambda)}
    \leq\sqrt{\mathcal E}.
\end{equation}
Together with \eqref{eq:tangency-quotient-limit}, Fatou's lemma also
gives
\begin{equation*}
    \norm{U_0}_{L^2(\Lambda)}
    \leq\sqrt{\mathcal E}.
\end{equation*}

Since $\Lambda$ is finite,
\eqref{eq:tangency-uniform-quotient-bound} implies that
$\{U_\varepsilon\}$ is uniformly integrable in $L^1(\Lambda)$.
Indeed, for every measurable set $B$,
\begin{equation*}
    \int_B|U_\varepsilon|\,d\Lambda
    \leq
    \sqrt{\mathcal E}\,\Lambda(B)^{1/2}.
\end{equation*}
Hence Vitali's theorem and
\eqref{eq:tangency-quotient-limit} yield
\begin{equation}\label{eq:tangency-L1-convergence}
    U_\varepsilon\longrightarrow U_0
    \qquad\text{in }L^1(\Lambda).
\end{equation}

It remains to account for the possibly unbounded target radius
$s\in L^2(\Lambda)$; see \eqref{eq:tangency-gluing-moments}.
For every $N>0$,
\begin{equation*}
\begin{aligned}
    \left|\int_{\mathcal S\times\cone^2} s(U_\varepsilon-U_0)\,d\Lambda\right|
    &\leq N\norm{U_\varepsilon-U_0}_{L^1(\Lambda)}
    +\norm{s\mathbf1_{\{s>N\}}}_{L^2(\Lambda)}
     \norm{U_\varepsilon-U_0}_{L^2(\Lambda)}\\
    &\leq N\norm{U_\varepsilon-U_0}_{L^1(\Lambda)}
    +2\sqrt{\mathcal E}\,
     \norm{s\mathbf1_{\{s>N\}}}_{L^2(\Lambda)}.
\end{aligned}
\end{equation*}
Letting first $\varepsilon\to0$ and using
\eqref{eq:tangency-L1-convergence}, and then $N\to\infty$, gives
\begin{equation*}
    \int_{\mathcal S\times\cone^2}sU_\varepsilon\,d\Lambda
    \longrightarrow\int_{\mathcal S\times\cone^2} sU_0\,d\Lambda.
\end{equation*}
Dividing \eqref{eq:tangency-cost-difference} by $\varepsilon$ and
passing to the limit therefore yields the two-sided derivative
\begin{equation}\label{eq:tangency-cost-derivative}
\begin{aligned}
    C'(0)
    &=
    -2\int_{\mathcal S\times\cone^2} s
    \left[
        \nabla K_y(\bfx_\tau)\cdot w(\bfx_\tau)
        +
        \frac12K_y(\bfx_\tau)\eta(\bfx_\tau)
    \right]\,d\Lambda.
\end{aligned}
\end{equation}
\emph{Conclusion}: It remains to identify this expression with the
normal pairing. Since $A$ is concentrated on reduced pairs and
has first radius one on $\cone_\Omega$, the normalized cone
logarithm there is
\begin{equation}\label{eq:normalized-log-kernel}
    \mathcal L_\cone([x,1],[y,s])
    =\left(s\nabla K_y(x),\,2(sK_y(x)-1)\right).
\end{equation}
The formula also covers $s=0$.
Using \eqref{eq:cone-barycenter-testing}, which is preserved by
the dilation, and the fact that $\mathcal W$ vanishes outside
$\Omega$, we obtain
\begin{equation*}
\begin{aligned}
    \langle\mathcal L,\mathcal W\rangle_\mu
    &=\int_{\cone^2}
      \left[
        s\nabla K_y(x)\cdot w(x)
        +\frac12sK_y(x)\eta(x)-\frac12\eta(x)
      \right]\,dA_\Omega.
\end{aligned}
\end{equation*}
The last term vanishes, since the first spatial marginal of
$A_\Omega$ is $\widetilde\mu$ and
\begin{equation*}
    \int_{\cone^2}\eta(x)\,dA_\Omega
    =\int_{\Omega}\eta\,d\widetilde\mu=0
\end{equation*}
by \eqref{eq:normal-reaction-mean}. Since the $\cone^2$-marginal
of $\Lambda$ is $A_\Omega$, it follows that
\begin{equation*}
    \langle\mathcal L,\mathcal W\rangle_\mu
    =\int s\left[
        \nabla K_y(\bfx_\tau)\cdot w(\bfx_\tau)
        +\frac12K_y(\bfx_\tau)\eta(\bfx_\tau)
      \right]\,d\Lambda
    =-\frac12C'(0).
\end{equation*}
Since \eqref{eq:tangency-cost-minimum} shows that $C$ has a local
minimum at zero and $C'(0)$ exists as a two-sided derivative,
$C'(0)=0$. Hence
\begin{equation*}
    \langle\mathcal L,\mathcal W\rangle_\mu=0.
\end{equation*}

We have proved the desired orthogonality for every compactly
supported normal field. For a general
$\mathcal W\in T_\mu\M^\perp$,
Lemma~\ref{lemma:compact-support-approximation} provides compactly
supported normal fields $\mathcal W_n\to\mathcal W$ in $H(\mu)$.
Since $\mathcal L\in H(\mu)$, continuity of the pairing gives
$\langle\mathcal L,\mathcal W\rangle_\mu=0$.
Finally, since $T_\mu\M$ is a closed subspace of $H(\mu)$,
\begin{equation*}
    \mathcal L\in(T_\mu\M^\perp)^\perp=T_\mu\M.
    \qedhere
\end{equation*}
\end{proof}

\subsection{First-order independence of cut-locus selections}

For a fixed semicoupling $\bfgamma$, the choice of $u$ affects only the
transport component of the logarithm on the cut locus
\begin{equation}\label{eq:cut-locus}
\operatorname{Cut}
\coloneqq
\Bigl\{
(x,y)\in\R^d\times\R^d
\colon
\norm{x-y}=\frac{\pi}{2}
\Bigr\}.
\end{equation}
The following estimate shows that the difference between two
selections, paired with any fixed $L^2(\mu;\mathbb R^d)$ field,
is negligible at first order in the WFR distance.

\begin{proposition}\label{propo:irrelevance-of-the-choice}
Assume $\WFR(\mu,\mu_n)\to0$. Let $\bfgamma^n=(\gamma_0^n,\gamma_1^n)$
be an optimal semicoupling between $\mu$ and $\mu_n$, and let 
\begin{equation*}
    u^1_n,u^2_n\colon \R^d\times\R^d \to [0,1]
\end{equation*}
be measurable functions. Write
\begin{equation*}
    \gamma_0^n=a_n\lambda_n,
\qquad
\gamma_1^n=b_n\lambda_n,
\end{equation*}
for some common dominating measure $\lambda_n$. Set also $d_n\coloneqq\WFR(\mu,\mu_n)$ and
\begin{equation*}
    \Log_{\bfgamma_n, u_n^i} \coloneqq (\Log^{T}_{n,i}, \Log^{R}_{n,i}) \qquad i=1,2.
\end{equation*}
Then, for every $\xi\in L^2(\mu;\R^d)$,
\begin{equation}\label{eq:cut-choice-irrelevance}
\int_{\operatorname{Cut}}
\xi(x)\cdot
\bigl(\Log_{n,2}^T(x,y)-\Log_{n, 1}^T(x,y)\bigr)
\,d\gamma_0^n(x,y)
=
o(d_n).
\end{equation}
\end{proposition}
\begin{proof}
Define
\begin{equation*}
    \Delta_n(x,y)\coloneqq\norm{\Log_{n,2}^T(x,y)-\Log_{n, 1}^T(x,y)}
\end{equation*}
By definition, since $u_n^1,u_n^2\in[0,1]$ and the logarithms coincide outside of $\operatorname{Cut}$ we have
\begin{equation*}
    \Delta_n(x,y)\leq \mathds 1 _{\operatorname{Cut}}(x,y) \lvert u_n^2(x,y)- u_n^1(x,y)\rvert \sqrt{\frac{b_n(x,y)}{a_n(x,y)}}\leq \mathds 1 _{\operatorname{Cut}}(x,y)\sqrt{\frac{b_n(x,y)}{a_n(x,y)}}.
\end{equation*}
Moreover, on $\operatorname{Cut}$ it holds
\begin{equation*}
    c(x,y,a_n(x,y),b_n(x,y))=a_n(x,y)+b_n(x,y)
\end{equation*}
and thus
\begin{equation*}
    \int_{\R^d\times\R^d}\Delta_n^2\,d\gamma_0^n(x,y)\leq \int_{\operatorname{Cut}}b_n(x,y)\,d\lambda_n(x,y)\leq \int_{\operatorname{Cut}} c(x,y,a_n(x,y),b_n(x,y))\,d\lambda_n(x,y)\leq d_n^2
\end{equation*}

Using the Cauchy--Schwarz inequality we get
\begin{align}
\left\lvert
\int_{\operatorname{Cut}}
\xi(x)\cdot
\bigl(\Log_{n, 2}^T-\Log_{n, 1}^T\bigr)(x,y)
\,d\gamma_0^n(x,y)
\right\rvert
\leq
 2 d_n\,\left(
\int_{\operatorname{Cut}}\norm{\xi(x)}^2\,d\gamma_0^n(x,y)
\right)^{1/2}
\label{eq:cut-choice-CS}
\end{align}

It remains to prove that
\begin{equation}\label{eq:cut-source-vanishing}
\int_{\operatorname{Cut}}|\xi(x)|^2\,d\gamma_0^n(x,y)
\longrightarrow0.
\end{equation}
To this end, define $\eta_n\coloneqq{\proj{1}\pf}\bigl(\mathds1_{\operatorname{Cut}}\gamma_0^n\bigr)$.
Since $\pi^1\pf\gamma_0^n=\mu$, one has $\eta_n\leq\mu$. Moreover,
\begin{equation}
\eta_n(\R^d)
=
\gamma_0^n(\operatorname{Cut})
=
\int_{\operatorname{Cut}}a_n\,d\lambda_n
\nonumber\
\leq
\int_{\operatorname{Cut}}(a_n+b_n)\,d\lambda_n
\leq d_n^2,
\label{eq:cut-mass-bound}
\end{equation}
and therefore $\eta_n(\R^d)\to0$, hence Corollary \ref{cor:vanishing-dominated-measures-L2} proves \eqref{eq:cut-source-vanishing}. 

Combining \eqref{eq:cut-choice-CS} and \eqref{eq:cut-source-vanishing} yields \eqref{eq:cut-choice-irrelevance}.
\end{proof}

We henceforth set $u\equiv1$ and write $\Log_{\bfgamma}$
and $\barlog_{\bfgamma}$ for the resulting logarithms.
Proposition~\ref{propo:irrelevance-of-the-choice} ensures that first-order pairings are independent
of this convention.

\section{Regular Wasserstein--Fisher--Rao calculus}\label{sec:wfr-calculus}

We can now introduce our notion of differentiability.

\begin{definition}[WFR differentiability]\label{def:wfr-differentiability}
Let $F\colon\M\to\R$ and $\mu\in\M$. We say that $F$ is WFR differentiable at $\mu$ if there exists
\begin{equation*}
    \mathcal V=(\xi,\zeta)\in T_\mu\M
\end{equation*}
such that, for every sequence $(\mu_n)_{n\in\N}\subset\M$ satisfying $\WFR(\mu_n,\mu)\to0$, 
and every sequence of optimal semicouplings
\begin{equation*}
    \bfgamma^n=(\gamma_0^n,\gamma_1^n)
    \in\operatorname{OptSC}(\mu,\mu_n),
\end{equation*}
one has
\begin{equation}\label{eq:wfr-differentiability}
\begin{split}
    F(\mu_n)& = F(\mu) +\langle \mathcal{V}, \barlog_{\bfgamma^n}\rangle_\mu
    +o\bigl(\WFR(\mu,\mu_n)\bigr)
\end{split}
\end{equation}
as $n\to+\infty$.
\end{definition}

We recall that we have the following equality due to the definition of the barycentric logarithm
\begin{equation*}
    \langle\mathcal V,\barlog_{\bfgamma^n}\rangle_\mu=\int_{\R^d\times\R^d} \xi(x) \cdot \Log^T_{\bfgamma^n}{(x,y)}\,d \gamma_0^n(x,y)+\frac{1}{4}\int_{\R^d\times\R^d} \zeta(x) \Log^R_{\bfgamma^n}{(x,y)}\,d \gamma_0^n(x,y).
\end{equation*}

\subsection{Strength and uniqueness of the gradients}

Before presenting differentiability results for specific functionals,
we show that Wasserstein--Fisher--Rao gradients are strong, meaning
that, for WFR differentiable functionals, first-order expansions also
hold along non-optimal semicouplings. As a byproduct, we obtain
uniqueness of the gradients.

The main ingredient is the following expansion estimate for smooth
expectations.

\begin{proposition}
\label{propo:strong-linear-observables}
Let $\varphi\in C_b^{1,1}(\R^d)$ and define
\begin{equation*}
    F(\mu)\coloneqq\int_{\R^d}\varphi\,d\mu.
\end{equation*}
Then there exists a constant $C_\varphi>0$ such that, for every
$\mu,\nu\in\M$ and every semicoupling
$\bfgamma\in\operatorname{SC}(\mu,\nu)$, one has
\begin{equation*}
\left|
    F(\nu)-F(\mu)
    -
    \left\langle
        (\nabla\varphi,4\varphi),
        \barlog_{\bfgamma}
    \right\rangle_{\mu}
\right|
\leq C_\varphi\,\mathcal C(\bfgamma),
\end{equation*}
with $C_\varphi$ as in \eqref{eq:quantitative-dependence}.
\end{proposition}

\begin{proof}
Let $\lambda$ be a finite positive measure dominating
$\gamma_0$ and $\gamma_1$, and set as usual
\begin{equation*}
    a=\frac{d\gamma_0}{d\lambda},
    \qquad
    b=\frac{d\gamma_1}{d\lambda}.
\end{equation*}
Throughout, terms involving the logarithm are set to zero on $\{a=0\}$.
Since $(\pi^1)_\#\gamma_0=\mu$ and $(\pi^2)_\#\gamma_1=\nu$,
we have
\begin{equation*}
    F(\nu)-F(\mu)
    =
    \int_{\R^d\times\R^d}
    \bigl(b\varphi(y)-a\varphi(x)\bigr)\,d\lambda(x,y).
\end{equation*}
By the definition of the barycentric logarithm and of the WFR inner product,
\begin{equation*}
    \left\langle
        (\nabla\varphi,4\varphi),
        \barlog_{\bfgamma}
    \right\rangle_{\mu}
    =
    \int_{\R^d\times\R^d}
    a\left(
        \nabla\varphi(x)\cdot\Log_{\bfgamma}^T(x,y)
        +\varphi(x)\Log_{\bfgamma}^R(x,y)
    \right)\,d\lambda(x,y).
\end{equation*}
Define
\begin{equation*}
\begin{aligned}
    R_\varphi(x,y)\coloneqq{}b\varphi(y)-a\varphi(x)-a\left(
        \nabla\varphi(x)\cdot\Log_{\bfgamma}^T(x,y)
        +\varphi(x)\Log_{\bfgamma}^R(x,y)
    \right).
\end{aligned}
\end{equation*}
We first restrict to $\{a,b>0,\ \norm{y-x}\leq\pi/2\}$.
Consider the cone points
\begin{equation*}
    \y_0=[x,\sqrt{a(x,y)}],
    \qquad
    \y_1=[y,\sqrt{b(x,y)}].
\end{equation*}
Neither is the vertex, and $(\y_0,\y_1)\notin\mathfrak B$.
Recalling the choice $u\equiv1$, we have
\begin{equation*}
    \Log_{\bfgamma}^T(x,y)
    =\Log^T_\cone(\y_0,\y_1),
    \qquad
    \Log_{\bfgamma}^R(x,y)
    =\frac{2}{\sqrt{a(x,y)}}\Log^R_\cone(\y_0,\y_1).
\end{equation*}
Let $\bfy_t$ be the reduced cone geodesic joining these points and set
$h(t)\coloneqq\Phi_\varphi(\bfy_t)$.
Lemma~\ref{lemma:cone-observable-estimates} therefore gives
\begin{equation*}
    |R_\varphi(x,y)|
    =|h(1)-h(0)-h'(0^+)|
    \leq C_\varphi\dc^2(\y_0,\y_1)
    =C_\varphi\,c(x,y,a,b),
\end{equation*}
where we used $\dc^2(\y_0,\y_1)=\dcwfr^2(\y_0,\y_1)$
since $\norm{y-x}\leq\pi/2$.

On $\{a,b>0,\ \norm{y-x}>\pi/2\}$, one has
\begin{equation*}
    \Log_{\bfgamma}^T(x,y)=0,
    \qquad
    \Log_{\bfgamma}^R(x,y)=-2.
\end{equation*}
Hence $R_\varphi(x,y)=b\varphi(y)+a\varphi(x)$, and therefore
\begin{equation*}
    |R_\varphi(x,y)|
    \leq\norm{\varphi}_\infty(a+b)
    \leq C_\varphi\,c(x,y,a,b),
\end{equation*}
since $c(x,y,a,b)=a+b$ on this set.
Similarly,
\begin{equation*}
    R_\varphi(x,y)=b\varphi(y)\quad\text{on }\{a=0,\ b>0\},
    \qquad
    R_\varphi(x,y)=a\varphi(x)\quad\text{on }\{a>0,\ b=0\},
\end{equation*}
while $R_\varphi(x,y)=0$ on $\{a=b=0\}$.
Since $c(x,y,a,b)=a+b$ whenever $ab=0$, we conclude in every case that
\begin{equation*}
    |R_\varphi(x,y)|\leq C_\varphi\,c(x,y,a,b).
\end{equation*}
Integrating with respect to $\lambda$ and using \eqref{eq:semi-coupling-cost}
proves the claimed estimate.
\end{proof}

\begin{theorem}\label{thm:strongness}
Let $F\colon\M\to\R$ be WFR differentiable at $\mu\in\M$ and $\mathcal V\in T_\mu\M$ a tangent vector satisfying the expansion in Definition \ref{def:wfr-differentiability}. Then for every sequence $(\mu_n)_{n\in\N}\subset\M$
and every sequence of semicouplings
\begin{equation*}
    \bfgamma^n=(\gamma_0^n,\gamma_1^n)
    \in\operatorname{SC}(\mu,\mu_n),
\end{equation*}
with $\varepsilon^2_n\coloneqq \mathcal{C}(\bfgamma^n)\to 0$ one has
\begin{equation}\label{eq:strong-wfr-expansion}
\begin{split}
    F(\mu_n)= F(\mu) +\langle \mathcal{V}, \barlog_{\bfgamma^n}\rangle_\mu
    +o\bigl(\varepsilon_n\bigr)
\end{split}
\end{equation}    
\end{theorem}
\begin{proof}
For every $n$, choose $\tilde\bfgamma^n\in\oSC(\mu,\mu_n)$, 
and set $d_n\coloneqq\WFR(\mu,\mu_n)$.
By optimality, $d_n\leq\varepsilon_n$. The assumption $\varepsilon_n\to0$ therefore implies
$\mu_n\to\mu$ in WFR. 

If $\varepsilon_n=0$, then $\mu_n=\mu$ and the barycentric
logarithm associated with $\bfgamma^n$ vanishes.
Thus, \eqref{eq:strong-wfr-expansion} holds trivially for
these indices, and we may restrict to those for which
$\varepsilon_n>0$.

By WFR differentiability,
\begin{equation}\label{eq:optimal-expansion-for-strongness}
    F(\mu_n)
    =
    F(\mu)
    +
    \left\langle
        \mathcal V,
        \barlog_{\tilde\bfgamma^n}
    \right\rangle_{\mu}
    +
    r_n,
    \qquad
    r_n=o(d_n).
\end{equation}

It remains to compare the two barycentric logarithms. Define $T_n\colon T_\mu\M\to\R$
by
\begin{equation*}
        T_n\mathcal W
    \coloneqq
    \frac{
    \left\langle
        \mathcal W,
        \barlog_{\tilde\bfgamma^n}
        -
        \barlog_{\bfgamma^n}
    \right\rangle_{\mu}
    }{\varepsilon_n}.
\end{equation*}
By the energy estimate for barycentric logarithms \eqref{eq:energy-bound},
\begin{equation*}
        \left\|
        \barlog_{\tilde\bfgamma^n}
    \right\|_{\mu}
    \leq d_n,
    \qquad
    \left\|
        \barlog_{\bfgamma^n}
    \right\|_{\mu}
    \leq\varepsilon_n.
\end{equation*}
Thus
\begin{equation*}
    |T_n\mathcal W|\leq\frac{d_n+\varepsilon_n}{\varepsilon_n}
    \|\mathcal W\|_{\mu}\leq
    2\|\mathcal W\|_{\mu}.
\end{equation*}
In particular,
\begin{equation}\label{eq:Tn-uniform-bound}
    \sup_n\|T_n\|\leq2.
\end{equation}

Consider the dense subspace
\begin{equation*}
        D
    \coloneqq
    \left\{
        (\nabla\varphi,4\varphi):
        \varphi\in C_c^\infty(\R^d)
    \right\}
    \subset T_\mu\M.
\end{equation*}
Fix $ \mathcal W_\varphi=(\nabla\varphi,4\varphi)\in D$.
Applying Proposition~\ref{propo:strong-linear-observables} first to
$\bfgamma^n$ and then to $\tilde\bfgamma^n$, we obtain
\begin{equation*}
    \int\varphi\,d\mu_n-\int\varphi\,d\mu
    =
    \left\langle
        \mathcal W_\varphi,
        \barlog_{\bfgamma^n}
    \right\rangle_{\mu}
    +R_n,
\end{equation*}
 and
\begin{equation*}
        \int\varphi\,d\mu_n-\int\varphi\,d\mu
    =
    \left\langle
        \mathcal W_\varphi,
        \barlog_{\tilde\bfgamma^n}
    \right\rangle_{\mu}
    +\tilde R_n
\end{equation*}
with $|R_n|\leq C_\varphi\varepsilon_n^2$ and $    |\tilde R_n| \leq C_\varphi d_n^2$ respectively.
Subtracting yields
\begin{equation*}
    |T_n\mathcal W_\varphi|
    \leq
    C_\varphi
    \frac{\varepsilon_n^2+d_n^2}{\varepsilon_n}
    \longrightarrow0.
\end{equation*}
Hence $T_n\to0$ pointwise on $D$. Using the density of $D$ and the uniform bound \eqref{eq:Tn-uniform-bound} one readily extends the convergence to the full tangent space $T_\mu\M$.

Taking $\mathcal W=\mathcal V$ gives
\begin{equation*}
        \left\langle
        \mathcal V,
        \barlog_{\tilde\bfgamma^n}
        -
        \barlog_{\bfgamma^n}
    \right\rangle_{\mu}
    =
    \varepsilon_n T_n\mathcal V
    =
    o(\varepsilon_n).
\end{equation*}
Finally, $r_n=o(d_n)$ and $d_n\leq\varepsilon_n$ imply $r_n=o(\varepsilon_n)$. Substituting the preceding estimate into
\eqref{eq:optimal-expansion-for-strongness} proves
\eqref{eq:strong-wfr-expansion}.
\end{proof}

As a first consequence, we obtain first-order expansions along the
natural family of perturbations
\begin{equation}\label{eq:eps-perturbation}
    \mu^\varepsilon
    \coloneqq
    \left(\operatorname{Id}+\varepsilon T\right)_\#
    \left(1+\frac{\varepsilon}{2}\rho\right)^2\mu.
\end{equation}

\begin{corollary}\label{cor:natural-perturbations}
Let $F\colon\M\to\R$ be WFR differentiable at $\mu$, and let
$\mathcal V=(\xi,\zeta)\in T_\mu\M$ be a tangent vector realizing its
WFR differential. Then, for every
\begin{equation*}
    T\in L^2(\mu;\R^d),
    \qquad
    \rho\in L^2(\mu),
\end{equation*}
one has
\begin{equation}\label{eq:directional-derivative-natural-perturbation}
    \lim_{\varepsilon\to0}
    \frac{F(\mu^\varepsilon)-F(\mu)}{\varepsilon}
    =
    \int_{\R^d}\xi\cdot T\,d\mu
    +
    \frac14\int_{\R^d}\zeta\rho\,d\mu,
\end{equation}
where $\mu^\varepsilon$ is defined by
\eqref{eq:eps-perturbation}.
\end{corollary}

\begin{proof}
Apply Lemma~\ref{lemma:graph-linearization} with
$T_\varepsilon=T$ and $\rho_\varepsilon=\rho$.
The resulting semicouplings satisfy
\begin{equation*}
    \barlog_{\bfgamma^\varepsilon}
    =\varepsilon(T,\rho)+o_{H(\mu)}(|\varepsilon|),
    \qquad
    \sqrt{\mathcal C(\bfgamma^\varepsilon)}=O(|\varepsilon|).
\end{equation*}
The strong expansion consequently gives
\begin{equation*}
    F(\mu^\varepsilon)-F(\mu)
    =\varepsilon\langle\mathcal V,(T,\rho)\rangle_\mu
    +o(|\varepsilon|),
\end{equation*}
which proves the claim.
\end{proof}

Uniqueness of the tangent vector
representing the WFR differential now follows readily.

\begin{theorem}[Uniqueness of the WFR gradient]
\label{thm:uniqueness}
Assume that $F\colon\M\to\R$ is WFR differentiable at
$\mu\in\M$. Then there exists a unique $\mathcal V\in T_\mu\M$
for which the first-order expansion in the definition of WFR
differentiability holds. We denote this vector by
\begin{equation*}
        \nabla_{\WFR}F(\mu)=(\nabla^TF(\mu),\nabla^RF(\mu))
\end{equation*}
and call it the \emph{WFR gradient of $F$ at $\mu$}.
\end{theorem}

\begin{proof}
Suppose that $\mathcal V_1,\mathcal V_2\in T_\mu\M$
both satisfy the definition of WFR differentiability. Let
$\varphi\in C_c^\infty(\R^d)$ and consider the perturbation
\begin{equation*}
        \mu^\varepsilon
    =
    \left(
        \operatorname{Id}+\varepsilon\nabla\varphi
    \right)_\#
    (1+2\varepsilon\varphi)^2\mu,
\end{equation*}
corresponding to $T=\nabla\varphi$ and $\rho=4\varphi$.
Applying Corollary~\ref{cor:natural-perturbations} to each
$\mathcal V_i$ gives
\begin{equation*}
        \lim_{\varepsilon\to0}
    \frac{F(\mu^\varepsilon)-F(\mu)}{\varepsilon}
    =
    \left\langle
        \mathcal V_i,
        (\nabla\varphi,4\varphi)
    \right\rangle_{\mu}.
\end{equation*}
Since the left-hand side is independent of $i$,
\begin{equation*}
    \left\langle
    \mathcal V_1-\mathcal V_2,
    (\nabla\varphi,4\varphi)
    \right\rangle_{\mu}
    =0
\end{equation*}
for every $\varphi\in C_c^\infty(\R^d)$. By definition,
\begin{equation*}
        \left\{
        (\nabla\varphi,4\varphi):
        \varphi\in C_c^\infty(\R^d)
    \right\}
\end{equation*}
is dense in $T_\mu\M$. Since
$\mathcal V_1-\mathcal V_2\in T_\mu\M$, it follows that $\mathcal V_1-\mathcal V_2=0$ in $T_\mu\M$. Therefore $\mathcal V_1=\mathcal V_2$.
\end{proof}

We conclude this subsection by identifying the norm of the
WFR gradient with the descending metric slope. Unlike the
WFR differential, the latter is defined entirely in terms
of the metric and quantifies the maximal local rate of
decrease of a functional. Their equality shows that the
differential structure developed here recovers the intrinsic
metric notion of steepest descent, without regularity
assumptions on the reference measure.

\begin{proposition}\label{propo:metric-slope}
Let $F\colon \M\to\R$ be WFR differentiable at
$\mu\in\M$. Then
\begin{equation*}
    |\partial^- F|(\mu)
    =
    \|\nabla_{\WFR}F(\mu)\|_\mu,
\end{equation*}
where the descending metric slope is defined by
\begin{equation*}
    |\partial^-F|(\mu)
    \coloneqq
    \limsup_{\substack{\nu\to\mu\\\nu\neq\mu}}
    \frac{(F(\mu)-F(\nu))_+}{\WFR(\mu,\nu)}.
\end{equation*}
\end{proposition}
\begin{proof}
For any $\nu\to\mu$ and
$\bfgamma\in\oSC(\mu,\nu)$,
WFR differentiability and Lemma~\ref{lemma:energy-bound} give
\begin{equation*}
    F(\mu)-F(\nu)
\le
\|\nabla_\WFR F(\mu)\|_\mu\WFR(\mu,\nu)
+o(\WFR(\mu,\nu)).
\end{equation*}
Hence $|\partial^-F|(\mu)\le\|\nabla_\WFR F(\mu)\|_\mu$.

If $\nabla_\WFR F(\mu)=0$, this proves the claim. Otherwise, set
\begin{equation*}
    (T,\rho)\coloneqq -\frac{\nabla_\WFR F(\mu)}{\|\nabla_\WFR F(\mu)\|_\mu}
\end{equation*}
and consider
\begin{equation*}
    \mu^\varepsilon\coloneqq 
(\mathrm{Id}+\varepsilon T)_\#
\left(
\left(1+\frac{\varepsilon}{2}\rho\right)^2\mu
\right),
\qquad \varepsilon>0.
\end{equation*}
By Corollary~\ref{cor:natural-perturbations},
\begin{equation*}
    F(\mu)-F(\mu^\varepsilon)
=\varepsilon\|\nabla_\WFR F(\mu)\|_\mu+o(\varepsilon).
\end{equation*}
Moreover, Lemma~\ref{lemma:graph-linearization} provides semicouplings
$\bfgamma^\varepsilon$ between $\mu$ and
$\mu^\varepsilon$ such that
\begin{equation*}
    \mathcal C(\bfgamma^\varepsilon)
=\varepsilon^2+o(\varepsilon^2).
\end{equation*}
Therefore,
\begin{equation*}
    \begin{aligned}
|\partial^-F|(\mu)\ge
\limsup_{\varepsilon\downarrow0}
\frac{F(\mu)-F(\mu^\varepsilon)}
{\WFR(\mu,\mu^\varepsilon)}\ge
\lim_{\varepsilon\downarrow0}
\frac{\varepsilon\|\nabla_\WFR F(\mu)\|_\mu+o(\varepsilon)}
{\sqrt{\mathcal C(\bfgamma^\varepsilon)}}=\|\nabla_\WFR F(\mu)\|_\mu.
\end{aligned}\qedhere
\end{equation*}
\end{proof}

\subsection{Calculus rules and differentiability criteria}

\subsubsection{A differentiability criterion and chain rules}

We present a Wasserstein--Fisher--Rao differentiability criterion based on flat derivatives.

\begin{theorem}[From flat to WFR differentiability]
\label{thm:differentiability-criterion}
Let $F\colon\M\to\R$ and fix $\mu\in\M$. Assume that there exist
a convex WFR-open neighborhood $U$ of $\mu$ and a function
$f\colon U\times\R^d\to\R$ satisfying the following properties:
\begin{enumerate}
    \item For every $\mu_0,\mu_1\in U$,
    \begin{equation}\label{eq:flat-ftc}
        F(\mu_1)-F(\mu_0)
        =
        \int_0^1\int_{\R^d}
        f\bigl((1-t)\mu_0+t\mu_1,x\bigr)
        \,d(\mu_1-\mu_0)(x)\,dt.
    \end{equation}

    \item For every $\nu\in U$, the function
    $f(\nu,\cdot)$ belongs to $C_b^{1,1}(\R^d)$, and
    there exists $K<+\infty$ such that
    \begin{equation}\label{eq:flat-uniform-spatial-bound}
        \sup_{\nu\in U}
        \left(
            \|f(\nu,\cdot)\|_\infty
            +
            \frac12\|\nabla_x f(\nu,\cdot)\|_\infty
            +
            \frac12\operatorname{Lip}\bigl(\nabla_x f(\nu,\cdot)\bigr)
        \right)
        \le K.
    \end{equation}

    \item The maps
    \begin{equation*}
        (\nu,x)\longmapsto f(\nu,x),
        \qquad
        (\nu,x)\longmapsto\nabla_x f(\nu,x)
    \end{equation*}
    are continuous on $U\times\R^d$, where $U$ is endowed
    with the WFR topology.
\end{enumerate}
Here $C_b^{1,1}(\R^d)$ denotes the space of bounded continuously
differentiable functions with bounded, globally Lipschitz gradient.

Then $F$ is WFR differentiable at $\mu$ and 
\begin{equation}\label{eq:flat-wfr-representative}
    \nabla_\WFR F(\mu)
    =
    \bigl(\nabla_x f(\mu,\cdot),4f(\mu,\cdot)\bigr)
    \in T_\mu\M.
\end{equation}
More precisely, there exist $r_0>0$ and a nondecreasing function
$\omega_\mu\colon[0,r_0)\to[0,+\infty)$, with
$\omega_\mu(r)\to0$ as $r\downarrow0$, such that
\begin{equation}\label{eq:flat-wfr-strong-estimate}
\begin{split}
    \left|
        F(\nu)-F(\mu)
        -
        \left\langle
            \nabla_\WFR F(\mu),
            \barlog_{\bfgamma}
        \right\rangle_\mu
    \right|
    \le
    K\,\mathcal C(\bfgamma)
    +
    \omega_\mu\bigl(\WFR(\mu,\nu)\bigr)
    \sqrt{\mathcal C(\bfgamma)}
\end{split}
\end{equation}
for every $\nu\in\M$ with $\WFR(\mu,\nu)<r_0$ and every
$\bfgamma\in\operatorname{SC}(\mu,\nu)$.
\end{theorem}

\begin{proof}
For $\nu\in U$, write
\begin{equation*}
    \mathcal V_\nu
    \coloneqq
    \bigl(\nabla f(\nu,\cdot),4f(\nu,\cdot)\bigr).
\end{equation*}
We have
\begin{equation*}
    \|\mathcal V_\nu-\mathcal V_\mu\|_\mu^2
    =
    \int_{\R^d}
    \left(
        \norm{\nabla f(\nu,x)-\nabla f(\mu,x)}^2
        +
        4|f(\nu,x)-f(\mu,x)|^2
    \right)\,d\mu(x).
\end{equation*}
We notice that $\mathcal V_\mu\in T_\mu\M$ due to Lemma \ref{lemma:C1,1-tangency}. Next, the continuity assumptions and
\eqref{eq:flat-uniform-spatial-bound} imply
\begin{equation}\label{eq:flat-field-continuity}
    \|\mathcal V_\nu-\mathcal V_\mu\|_\mu
    \longrightarrow0
    \qquad\text{as }\WFR(\nu,\mu)\to0.
\end{equation}
Indeed, for any sequence $\nu_n\to\mu$ in WFR, both
$f(\nu_n,x)\to f(\mu,x)$ and
$\nabla f(\nu_n,x)\to\nabla f(\mu,x)$ hold for every $x$.
The integrand defining the squared norm is uniformly bounded
by \eqref{eq:flat-uniform-spatial-bound}, so
\eqref{eq:flat-field-continuity} follows from dominated convergence.

Choose $r_0>0$ such that $B_{\WFR}(\mu,r_0)\subset U$,
and, for $0\le r<r_0$, define
\begin{equation}\label{eq:flat-continuity-modulus}
    \omega_\mu(r)
    \coloneqq
    \sup_{\WFR(\nu,\mu)\le r}
    \|\mathcal V_\nu-\mathcal V_\mu\|_\mu.
\end{equation}
This function is finite and nondecreasing, and
\eqref{eq:flat-field-continuity} implies
$\omega_\mu(r)\to0$ as $r\downarrow0$.

Fix now $\nu$ with $\WFR(\mu,\nu)<r_0$ and set
\begin{equation*}
    \mu_t\coloneqq(1-t)\mu+t\nu,
    \qquad t\in[0,1].
\end{equation*}
Convexity of the squared WFR distance \cite[Theorem 6.2 and Theorem 7.20]{liero-mielke-savare2018} under linear interpolation
gives
\begin{equation}\label{eq:flat-segment-distance}
    \WFR^2(\mu,\mu_t)
    \le (1-t)\WFR^2(\mu,\mu)+t\WFR^2(\mu,\nu)=t\WFR^2(\mu,\nu).
\end{equation}
In particular, $\mu_t\in U$ and for every $t\in[0,1]$
\begin{equation}\label{eq:flat-segment-field-bound}
    \|\mathcal V_{\mu_t}-\mathcal V_\mu\|_\mu
    \le \omega_\mu(\WFR(\mu_t,\mu))\leq
    \omega_\mu\bigl(\WFR(\mu,\nu)\bigr)
\end{equation}

Let $\bfgamma\in\operatorname{SC}(\mu,\nu)$.
By Proposition~\ref{propo:strong-linear-observables}, applied
to $f(\mu_t)$, and the uniform bound
\eqref{eq:flat-uniform-spatial-bound}, we have
\begin{equation}\label{eq:flat-test-expansion}
    \int_{\R^d}f(\mu_t)\,d(\nu-\mu)
    =
    \left\langle
        \mathcal V_{\mu_t},
        \barlog_{\bfgamma}
    \right\rangle_\mu
    +
    R_t,
    \qquad
    |R_t|\le K\,\mathcal C(\bfgamma).
\end{equation}
Here we use the quantitative dependence of the constant in \eqref{eq:quantitative-dependence}. Notice that the same semicoupling $\bfgamma$ is used for every $t$. By \eqref{eq:flat-field-continuity}, $t\mapsto \mathcal V_{\mu_t}$ is continuous in $H(\mu_)$, hence measurable. Therefore, \eqref{eq:flat-test-expansion} yields measurability of $t\mapsto R_t$ as well.  Integrating \eqref{eq:flat-test-expansion} and using
\eqref{eq:flat-ftc}, we obtain
\begin{align*}
    F(\nu)-F(\mu)
    -
    \left\langle
        \mathcal V_\mu,
        \barlog_{\bfgamma}
    \right\rangle_\mu=
    \int_0^1
    \left\langle
        \mathcal V_{\mu_t}-\mathcal V_\mu,
        \barlog_{\bfgamma}
    \right\rangle_\mu\,dt
    +
    \int_0^1R_t\,dt.
\end{align*}
These integrals are well defined by the continuity assumptions
and the uniform bounds. By Cauchy--Schwarz and the logarithm
energy bound \eqref{eq:energy-bound} we deduce
\begin{align*}
    \left|
        F(\nu)-F(\mu)
        -
        \left\langle
            \mathcal V_\mu,
            \barlog_{\bfgamma}
        \right\rangle_\mu
    \right|\le
    \sqrt{\mathcal C(\bfgamma)}
    \int_0^1
    \|\mathcal V_{\mu_t}-\mathcal V_\mu\|_\mu\,dt
    +
    K\,\mathcal C(\bfgamma).
\end{align*}
Together with \eqref{eq:flat-segment-field-bound}, this proves
\eqref{eq:flat-wfr-strong-estimate}.

Finally, if $\bfgamma\in\operatorname{OptSC}(\mu,\nu)$, then
$\mathcal C(\bfgamma)=\WFR^2(\mu,\nu)$, and hence
\begin{align*}
    \left|
        F(\nu)-F(\mu)
        -
        \left\langle
            \mathcal V_\mu,
            \barlog_{\bfgamma}
        \right\rangle_\mu
    \right|\le
    K\,\WFR^2(\mu,\nu)
    +
    \omega_\mu\bigl(\WFR(\mu,\nu)\bigr)\WFR(\mu,\nu)
    =
    o\bigl(\WFR(\mu,\nu)\bigr).
\end{align*}
The estimate is uniform over all optimal semicouplings, which
proves WFR differentiability at $\mu$.
\end{proof}

Nonetheless, as we show in Section \ref{subsec:flat-derivatives}, there are functionals that are WFR differentiable but do not admit flat derivatives. 
This reflects the fact that certain flat perturbations
of order $\varepsilon$ have WFR size of order
$\sqrt{\varepsilon}$.

We end this section with two chain rules.

\begin{proposition}[Chain rule along WFR geodesics]
\label{propo:wfr-geodesic-chain-rule}
Let $(\bfmu_t)_{t\in[0,1]}$ be a constant-speed WFR geodesic
and let $F\colon\M\to\R$.
For $s\in[0,1)$, denote by $\mathcal L_s$ the regular barycentric
logarithm associated with the rescaled geodesic
\begin{equation*}
    \tau\mapsto\bfmu_{s+(1-s)\tau},
    \qquad \tau\in[0,1],
\end{equation*}
and set $\mathcal V_s\coloneqq{\mathcal L_s}/{(1-s)}$.

If $F$ is WFR differentiable at $\mu_s$, then
\begin{equation*}
    \lim_{\delta\to0}
    \frac{F(\bfmu_{s+\delta})-F(\bfmu_s)}{\delta}
    =
    \langle\nabla_{\WFR}F(\bfmu_s),\mathcal V_s\rangle_{\bfmu_s}.
\end{equation*}
At $s=0$, the limit is understood as a right derivative.
In particular,
\begin{equation*}
    F(\mu_t)
    =
    F(\mu_0)
    +t\langle\nabla_{\WFR}F(\bfmu_0),\mathcal L_0\rangle_{\bfmu_0}
    +o(t)
\end{equation*}
whenever $F$ is WFR differentiable at $\mu_0$.

If, moreover, $F$ is WFR differentiable at every $\mu_t$ and
\begin{equation*}
    C_F\coloneqq
    \sup_{t\in[0,1]}
    \norm{\nabla_{\WFR}F(\bfmu_t)}_{\bfmu_t}
    <\infty,
\end{equation*}
then $t\mapsto F(\bfmu_t)$ is Lipschitz, with Lipschitz
constant at most $C_F\WFR(\bfmu_0,\bfmu_1)$, and
\begin{equation*}
    F(\bfmu_{t_2})-F(\bfmu_{t_1})
    =
    \int_{t_1}^{t_2}
    \langle\nabla_{\WFR}F(\bfmu_t),\mathcal V_t\rangle_{\mu_t}\,dt,
    \qquad 0\le t_1\le t_2\le1.
\end{equation*}
\end{proposition}

\begin{proof}
Set $D\coloneqq\WFR(\mu_0,\mu_1)^2$. If $D=0$, the geodesic is constant and the conclusions
are immediate. Otherwise, let $\boldsymbol{\Pi}$ be an optimal cone-geodesic representation of
$(\bfmu_t)_{t\in[0,1]}$. Then, for $\varphi\in C_c^\infty(\R^d)$, we have
\begin{equation*}
    f_\varphi(t)\coloneqq\int_{\R^d}\varphi\,d\bfmu_t=\int_{\Geo^*(\cone)}
    \bfr_t^2\varphi(\bfx_t)\,d\boldsymbol{\Pi}(\bfy),
\end{equation*}
where we recall $\bfy_t=[\bfx_t,\bfr_t]$. Along each reduced cone geodesic, the integrand has a continuous derivative, with one-sided derivatives at the endpoints. By Lemma~\ref{lemma:cone-observable-estimates}, its derivative is bounded in absolute value by
$B_\varphi(\bfr_0^2+\bfr_1^2)$, an integrable function.
Differentiation under the integral and dominated convergence
therefore gives $f_\varphi\in C^1([0,1])$, with one-sided
derivatives at the endpoints.

By the definition of $\mathcal L_s$,
\begin{equation*}
\begin{split}
    \langle(\nabla\varphi,4\varphi),\mathcal L_s\rangle_{\bfmu_s}
    =
    \lim_{\tau\downarrow0}
    \frac{f_\varphi(s+(1-s)\tau)-f_\varphi(s)}{\tau}=(1-s)f_\varphi'(s).
\end{split}
\end{equation*}
Consequently,
\begin{equation*}
    f_\varphi'(s)=\langle(\nabla\varphi,4\varphi),\mathcal V_s\rangle_{\bfmu_s}.
\end{equation*}
Moreover, the logarithm norm bound gives
\begin{equation*}
    \norm{\mathcal V_s}_{\bfmu_s}
    \le
    \frac{\WFR(\bfmu_s,\bfmu_1)}{1-s}
    =\sqrt D.
\end{equation*}

Fix $s\in[0,1)$ at which $F$ is WFR differentiable.
For $\delta\ne0$ sufficiently small with $s+\delta\in[0,1]$,
choose
$\gamma^\delta\in\operatorname{OptSC}(\bfmu_s,\bfmu_{s+\delta})$
and set
\begin{equation*}
    Z_\delta\coloneqq 
    \frac1\delta\barlog_{\gamma^\delta}.
\end{equation*}
Since $\WFR(\bfmu_s,\bfmu_{s+\delta})=|\delta|\sqrt D$, we have $\norm{Z_\delta}_{\bfmu_s}\le\sqrt  D$.

For every $\varphi\in C_c^\infty(\R^d)$,
Proposition~\ref{propo:strong-linear-observables} gives
\begin{equation*}
    \left|
        \langle(\nabla\varphi,4\varphi),Z_\delta\rangle_{\bfmu_s}
        -
        \frac{f_\varphi(s+\delta)-f_\varphi(s)}{\delta}
    \right|
    \le C_\varphi |\delta|D.
\end{equation*}
Hence
\begin{equation*}
    \langle(\nabla\varphi,4\varphi),Z_\delta-\mathcal V_s\rangle_{\bfmu_s}
    \longrightarrow0.
\end{equation*}
The uniform bounds on $Z_\delta$ and $\mathcal V_s$, together with
density of the fields $(\nabla\varphi,4\varphi)$ in
$\mathrm{T}_{\bfmu_s}\M$, imply
\begin{equation*}
    \langle\nabla_{\WFR}F(\bfmu_s),Z_\delta-\mathcal V_s\rangle_{\bfmu_s}
    \longrightarrow0.
\end{equation*}
WFR differentiability at $\bfmu_s$ now yields
\begin{equation*}
\begin{split}
    F(\bfmu_{s+\delta})-F(\bfmu_s)
    =
    \left\langle
        \nabla_{\WFR}F(\bfmu_s),
        \barlog_{\gamma^\delta}
    \right\rangle_{\bfmu_s}
    +o(|\delta|)=
    \delta\langle\nabla_{\WFR}F(\bfmu_s),\mathcal V_s\rangle_{\bfmu_s}
    +o(|\delta|).
\end{split}
\end{equation*}
This proves the pointwise chain rule, including the
right derivative at $s=0$.

Under the additional assumptions, $g(t)\coloneqq F(\bfmu_t)$ is
continuous on $[0,1]$ and differentiable on $(0,1)$, with
\begin{equation*}
    |g'(t)|
    \le
    \norm{\nabla_{\WFR}F(\bfmu_t)}_{\bfmu_t}
    \norm{\mathcal V_t}_{\bfmu_t}
    \le C_F\sqrt D.
\end{equation*}
The mean value theorem shows that $g$ is Lipschitz.
The integral identity follows from the fundamental
theorem of calculus.
\end{proof}

\begin{proposition}[Chain rule]
Let $F\colon\M\to\R$ be WFR differentiable at $\mu$ and let
$G\colon\R\to\R$ be differentiable at $F(\mu)$. Then $G\circ F$ is WFR
differentiable at $\mu$ and
\begin{equation*}
    \nabla_{\WFR}(G\circ F)(\mu)
    =
    G'(F(\mu))\,\nabla_{\WFR}F(\mu).
\end{equation*}
\end{proposition}

\begin{proof}
For any sequence
$\mu_n\to\mu$ in $\WFR$ and any
$\bfgamma^n\in\operatorname{OptSC}(\mu,\mu_n)$, we have
\begin{equation*}
    F(\mu_n)-F(\mu)
    =
    \left\langle
        \nabla_{\WFR}F(\mu),\barlog_{\bfgamma^n}
    \right\rangle_\mu
    +o\bigl(\WFR(\mu,\mu_n)\bigr).
\end{equation*}
Since by Cauchy-Schwarz and \eqref{eq:energy-bound}
\begin{equation*}
    \left|
    \left\langle
        \nabla_{\WFR}F(\mu),\barlog_{\bfgamma^n}
    \right\rangle_\mu
    \right|
    \leq
    \|\nabla_{\WFR}F(\mu)\|_{\mu}
    \WFR(\mu,\mu_n),
\end{equation*}
it follows that
\begin{equation*}
    |F(\mu_n)-F(\mu)|
    =
    O\bigl(\WFR(\mu,\mu_n)\bigr).
\end{equation*}
By differentiability of $G$ at $F(\mu)$,
\begin{equation*}
    G(F(\mu_n))-G(F(\mu))
    =
    G'(F(\mu))
    \bigl(F(\mu_n)-F(\mu)\bigr)
    +o\bigl(|F(\mu_n)-F(\mu)|\bigr).
\end{equation*}
Combining the previous estimates gives
\begin{equation*}
    G(F(\mu_n))-G(F(\mu))
    =
    \left\langle
        G'(F(\mu))\nabla_{\WFR}F(\mu),
        \barlog_{\bfgamma^n}
    \right\rangle_\mu
    +o\bigl(\WFR(\mu,\mu_n)\bigr),
\end{equation*}
and the conclusion follows.
\end{proof}

\subsubsection{Some differentiable functions}

As a consequence of the criterion developed in Theorem \ref{thm:differentiability-criterion} we obtain WFR differentiability for smooth enough expectations, cylinder functionals and interaction energies.

\begin{corollary}[Linear observables]
Let $f\in C_b^{1,1}(\R^d)$ and define
$F\colon\M\to\R$ by
\begin{equation*}
    F(\mu)
    \coloneqq
    \int_{\R^d} f(x)\,d\mu(x).
\end{equation*}
Then $F$ is everywhere WFR differentiable and
\begin{equation*}
    \nabla_{\WFR}F(\mu)
    =
    \bigl(\nabla f,4f\bigr).
\end{equation*}
\end{corollary}

\begin{corollary}[Cylinder functionals]
Let $N\in\N$, let $f_1,\ldots,f_N\in C_b^{1,1}(\R^d)$, and let
$G\colon\R^N\to\R$ be $C^1$. Define
$F\colon\M\to\R$ by
\begin{equation*}
    F(\mu)
    \coloneqq
    G\left(
        \int_{\R^d}f_1\,d\mu,\ldots,
        \int_{\R^d}f_N\,d\mu
    \right).
\end{equation*}
Then $F$ is everywhere WFR differentiable and
\begin{equation*}
    \nabla_{\WFR}F(\mu)
    =
    \left(
        \sum_{i=1}^N
        \partial_i G\left(
            \int f_1\,d\mu,\ldots,\int f_N\,d\mu
        \right)\nabla f_i,
        \;
        4\sum_{i=1}^N
        \partial_i G\left(
            \int f_1\,d\mu,\ldots,\int f_N\,d\mu
        \right)f_i
    \right).
\end{equation*}
\end{corollary}

\begin{corollary}[Interaction energies]
Let $W\in C_b^{1,1}(\R^d\times\R^d)$ be symmetric, namely
\begin{equation*}
    W(x,y)=W(y,x)
    \qquad
    \text{for every }x,y\in\R^d,
\end{equation*}
and define
\begin{equation*}
    F(\mu)
    \coloneqq
    \frac12
    \int_{\R^d}\int_{\R^d}
    W(x,y)\,d\mu(x)\,d\mu(y).
\end{equation*}
Then $F$ is everywhere WFR differentiable and
\begin{equation*}
    \nabla_{\WFR}F(\mu)(x)
    =
    \left(
        \int_{\R^d}\nabla_x W(x,y)\,d\mu(y),
        \;
        4\int_{\R^d}W(x,y)\,d\mu(y)
    \right).
\end{equation*}
\end{corollary}

\subsubsection{Compatibility with finite-dimensional derivatives}

A classical property of Wasserstein differentiability is its
compatibility with empirical measures. If $F$ is Wasserstein
differentiable at
\begin{equation*}
    \mu_{(x_1,\dots,x_N)}
    \coloneqq
    \frac1N\sum_{i=1}^N\delta_{x_i},
\end{equation*}
then the map $F_N(\boldsymbol x)\coloneqq F(\mu_{\boldsymbol x})$
is Fréchet differentiable at $(x_1,\dots,x_N)$, with
\begin{equation*}
    \nabla_{x_i}F_N(x_1,\dots,x_N)
    =
    \frac1N\nabla_{W_2}F(\mu_{(x_1,\dots,x_N)})(x_i),
    \qquad i=1,\ldots,N.
\end{equation*}
See, for example, \cite[Proposition~3.1]{carmona-delarue2015}.

The WFR differential extends this compatibility to weighted
empirical measures, permitting simultaneous perturbations of the
particle positions and weights. This provides a finite-dimensional
interpretation of both components of the WFR gradient.

To formulate this property, we consider the space
\begin{equation*}
    M_N
    \coloneqq
    (\R^d)^N\times\R^N.
\end{equation*}
For $(\hat x,\hat \ell)\coloneqq(x_1,\dots,x_N,\ell_1,\dots,\ell_N)\in M_N$,
define the empirical-measure map
\begin{equation}\label{eq:log-weight-empirical-embedding}
    \iota_N(\hat x,\hat \ell)
    \coloneqq
    \sum_{i=1}^N e^{\ell_i}\delta_{x_i}.
\end{equation}
We write increments as
\begin{equation*}
    \hat\xi\coloneqq(\xi_1,\dots,\xi_N),
    \qquad
    \hat\zeta\coloneqq(\zeta_1,\dots,\zeta_N).
\end{equation*}

\begin{theorem}[Compatibility with weighted empirical measures]
\label{thm:finite-dimensional-differentiability}
Let $F\colon\M\to\R$ be WFR differentiable at
$\mu=\iota_N(\hat x,\hat \ell)$.
Then the map
$\widetilde F\coloneqq F\circ\iota_N\colon M_N\to\R$
is Fréchet differentiable at $(\hat x,\hat\ell)$, with
\begin{equation}\label{eq:log-weight-finite-dimensional-differential}
    D\widetilde F(\hat x, \hat \ell)[\hat\xi,\hat\zeta]
    =
    \sum_{i=1}^N e^{\ell_i}
    \left(
        \nabla^TF(\mu)(x_i)
        \cdot\xi_i
        +
        \frac14\nabla^RF(\mu)(x_i)\zeta_i
    \right).
\end{equation}
In particular, its Euclidean partial derivatives are
\begin{equation*}
    \nabla_{x_i}\widetilde F(\hat x,\hat\ell)
    =
    e^{\ell_i}\nabla^TF(\mu)(x_i),
    \qquad
    \partial_{\ell_i}\widetilde F(\hat x,\hat\ell)
    =
    \frac14e^{\ell_i}\nabla^RF(\mu)(x_i).
\end{equation*}
\end{theorem}

\begin{proof}
Since $\mu$ has finite support, we can choose
$\varphi\in C_c^\infty(\R^d)$ such that
\begin{equation*}
    \nabla\varphi(x)
    =
    {\nabla^TF(\mu)(x)},
    \qquad
    4\varphi(x)
    =
    {\nabla^RF(\mu)(x)},
    \qquad x\in\operatorname{supp}\mu.
\end{equation*}
Indeed, the value and gradient of a smooth function can be
prescribed independently at finitely many distinct points.
The smooth expectation
$F_\varphi(\nu)\coloneqq\int_{\R^d}\varphi\,d\nu$
therefore has the same WFR gradient as $F$ at $\mu$.
As a consequence, for every $\nu\in\M$ and optimal
semicoupling $\bfgamma\in\oSC(\mu,\nu)$,
\begin{equation*}
\begin{split}
    \left\langle
        \nabla_\WFR F(\mu),
        \barlog_{\bfgamma}
    \right\rangle_\mu
    &=
    \left\langle
        (\nabla\varphi,4\varphi),
        \barlog_{\bfgamma}
    \right\rangle_\mu\\
    &=
    \int_{\R^d}\varphi\,d(\nu-\mu)
    +o\bigl(\WFR(\mu,\nu)\bigr)
\end{split}
\end{equation*}
as $\nu\to\mu$ in WFR, due to
Proposition~\ref{propo:strong-linear-observables}.
Hence, WFR differentiability of $F$ yields
\begin{equation*}
    F(\nu)-F(\mu)
    =
    \int_{\R^d}\varphi\,d(\nu-\mu)
    +o\bigl(\WFR(\mu,\nu)\bigr).
\end{equation*}

For an increment
$(\hat\xi,\hat\zeta)$,
set
\begin{equation*}
    \nu
    \coloneqq
    \iota_N(\hat x+\hat\xi,\hat\ell+\hat\zeta)
    \quad\text{and}\quad
    \delta\coloneqq
    \left(
        \sum_{i=1}^N\norm{\xi_i}^2
        +
        \sum_{i=1}^N\zeta_i^2
    \right)^{1/2}.
\end{equation*}
For sufficiently small increments, the cone coupling
matching
{$[x_i,e^{\ell_i/2}]$}
with
{$[x_i+\xi_i,e^{(\ell_i+\zeta_i)/2}]$}
yields
\begin{equation*}
    \WFR^2(\mu,\nu)
    \leq
    \sum_{i=1}^N e^{\ell_i}
    \left(
        1+e^{\zeta_i}
        -2e^{\zeta_i/2}\cos(\norm{\xi_i})
    \right)
    =O(\delta^2).
\end{equation*}
Consequently, Taylor expansion gives
\begin{align*}
    \widetilde F(\hat x+\hat\xi,\hat\ell+\hat\zeta)
    -
    \widetilde F(\hat x, \hat \ell)
    &=
    \sum_{i=1}^N e^{\ell_i}
    \left(
        e^{\zeta_i}\varphi(x_i+\xi_i)-\varphi(x_i)
    \right)
    +o(\delta)\\
    &=
    \sum_{i=1}^N e^{\ell_i}
    \left(
        {\nabla^TF(\mu)(x_i)}
        \cdot\xi_i
        +
        \frac14
        {\nabla^RF(\mu)(x_i)}
        \zeta_i
    \right)
    +o(\delta).
\end{align*}
This proves the claimed Fréchet differentiability.
\end{proof}

The preceding formula also admits a Riemannian interpretation.
Writing $w_i=e^{\ell_i}$, the identity
$\dot \ell_i=\dot w_i/w_i$ identifies logarithmic weight
velocities with relative mass variations.
We therefore endow $M_N$ with the Riemannian metric
defined, for
$\hat\xi^a=(\xi_i^a)_{i=1}^N$ and
$\hat\zeta^a=(\zeta_i^a)_{i=1}^N$, $a=1,2$, by
\begin{equation}\label{eq:log-weight-particle-metric}
\begin{split}
    g_{(\hat x,\hat\ell)}
    \bigl(
        {
        (\hat\xi^1,\hat\zeta^1),
        (\hat\xi^2,\hat\zeta^2)
        }
    \bigr)
    \coloneqq
    \sum_{i=1}^N e^{\ell_i}
    \left(
        \xi_i^1\cdot\xi_i^2
        +
        \frac14\zeta_i^1\zeta_i^2
    \right).
\end{split}
\end{equation}
By \eqref{eq:log-weight-finite-dimensional-differential},
the corresponding Riemannian gradient is
\begin{equation}\label{eq:log-weight-finite-dimensional-gradient}
    \operatorname{grad}_g\widetilde F(\hat x, \hat \ell)
    =
    \left(
        \bigl(
            {\nabla^TF(\mu)(x_i)}
        \bigr)_{i=1}^N,
        \bigl(
            {\nabla^RF(\mu)(x_i)}
        \bigr)_{i=1}^N
    \right),
\end{equation}
where $\mu\coloneqq\iota_N(\hat x, \hat \ell)$.

When the positions are pairwise distinct,
the map assigning to each tangent vector
$(\hat\xi,\hat\zeta)$ the field taking the value
$(\xi_i,\zeta_i)$ at $x_i$ identifies this metric with
the scalar product of $H(\mu)$.

\section{Differentiability of the squared Wasserstein--Fisher--Rao distance}
\label{sec:squared-wfr-distance}

We apply the preceding calculus to the functional
\begin{equation*}
    F_{\hat\mu}(\mu)\coloneqq\frac{1}{2}\WFR^2(\mu,\hat\mu),
    \qquad \hat\mu\in\M.
\end{equation*}
The main obstructions are source splitting, which entails a loss
of energy under barycentric averaging; target creation, whose cost
is invisible to source tangent fields; and non-uniqueness at the
transport cutoff.
The proof adapts the two-competitor argument of
Alfonsi and Jourdain for the squared Wasserstein distance
\cite{alfonsi-jourdain2020}, gluing cone couplings after
suitable radial dilations. Its main additional ingredient is
the stability of source logarithms and creation mass for
asymptotically optimal cone couplings.

\subsection{The differentiability theorem}
\label{subsec:squared-wfr-theorem}

We formulate the Monge condition after source normalization
by radial dilation. Let
\begin{equation*}
    \vartheta(\y_0,\y_2)
    \coloneqq
    \begin{cases}
        r_0, & r_0>0,\\
        1, & r_0=0.
    \end{cases}
\end{equation*}
We say that $(\mu,\hat\mu)$ admits a unique
source-normalized optimal Monge representation without
target creation if there exist Borel maps
$T\colon\R^d\to\R^d$ and
$q\colon\R^d\to[0,\infty)$ such that, writing
\begin{equation}\label{eq:squared-wfr-monge-plan}
    A^\ast
    \coloneqq
    \bigl(x\mapsto([x,1],[T(x),q(x)])\bigr)\pf\mu,
\end{equation}
every $A\in\operatorname{Opt\Gamma}_{\cone}(\mu,\hat\mu)$
satisfies
\begin{equation}\label{eq:squared-wfr-monge-assumption}
    \operatorname{dil}_{\vartheta,2}
    \bigl(A\llcorner\{r_0>0\}\bigr)
    =A^\ast,
    \qquad
    \int_{\{r_0=0\}}r_2^2\,dA=0.
\end{equation}
The value of $T$ on $\{q=0\}$ is irrelevant.
The coupling $A^\ast$ is optimal, and
$T\pf(q^2\mu)=\hat\mu$.
In particular, $q\in L^2(\mu)$.

Moreover, the set
\begin{equation*}
    R_{\mathrm{Cut}}\coloneqq\Bigl\{x:q(x)>0,\ \norm{T(x)-x}=\frac{\pi}{2}\Bigr\}
\end{equation*}
satisfies $\mu(R_{\mathrm{Cut}})=0$. Indeed, optimal cone couplings are concentrated on reduced
pairs by Theorem~\ref{thm:wfr-cone-formulation}. If $R_{\mathrm{Cut}}$
had positive $\mu$-measure, replacing each pair
$([x,1],[T(x),q(x)])$ over $R_{\mathrm{Cut}}$ by
$([x,1],\vertex)$ and $(\vertex,[T(x),q(x)])$
would preserve the homogeneous marginals and the cost, since
\begin{equation*}
    \dc^2([x,1],[T(x),q(x)])
    =1+q(x)^2
    =\dc^2([x,1],\vertex)+\dc^2(\vertex,[T(x),q(x)]).
\end{equation*}
The resulting optimal coupling $\widetilde A$ would satisfy
\begin{equation*}
    \int_{\{r_0=0\}}r_2^2\,d\widetilde A
    =\int_{R_{\mathrm{Cut}}} q(x)^2\,d\mu(x)>0,
\end{equation*}
contradicting \eqref{eq:squared-wfr-monge-assumption}.

Define
\begin{equation}\label{eq:squared-wfr-reference-log}
    \mathcal L(x)
    \coloneqq
    \mathcal L_{\cone}([x,1],[T(x),q(x)])
    =
    \left(
        q(x)\nabla K_{T(x)}(x),
        2\bigl(q(x)K_{T(x)}(x)-1\bigr)
    \right).
\end{equation}
The semicoupling
\begin{equation*}
    \gamma^\ast
    \coloneqq
    \left(
        (\operatorname{Id},T)\pf\mu,
        (\operatorname{Id},T)\pf(q^2\mu)
    \right)
\end{equation*}
is optimal and satisfies
$\mathcal L=\barlog_{\gamma^\ast}$.
Hence Theorem~\ref{propo:tangency} gives
\begin{equation}\label{eq:squared-wfr-reference-energy}
    \mathcal L\in\mathrm{T}_\mu\M,
    \qquad
    \norm{\mathcal L}_\mu^2
    =
    \int_{\cone^2}\dc^2(\y_0,\y_2)\,dA^\ast(y_0,\y_2)
    =
    \WFR^2(\mu,\hat\mu).
\end{equation}
The energy identity follows from
\eqref{eq:squared-wfr-cone-log-identities} and
\eqref{eq:squared-wfr-monge-plan}.

\begin{theorem}[Differentiability of the squared distance]
\label{thm:squared-wfr-differentiability}
Let $\mu,\hat\mu\in\M$. The following conditions are equivalent.
\begin{enumerate}
    \item[(i)] $F_{\hat\mu}$ is WFR differentiable at $\mu$.
    \item[(ii)] There is exactly one regular barycentric logarithm $\mathcal L$
    between $\mu$ and $\hat\mu$, and
    \begin{equation*}
        \norm{\mathcal L}_\mu=\WFR(\mu,\hat\mu).
    \end{equation*}
    \item[(iii)] The pair $(\mu,\hat\mu)$ has the unique
    source-normalized Monge representation without target creation
    in \eqref{eq:squared-wfr-monge-assumption}.
\end{enumerate}
Whenever these conditions hold, $\mathcal L$ is the field defined in
\eqref{eq:squared-wfr-reference-log}, and

\begin{equation*}
    \nabla_{\WFR}F_{\hat\mu}(\mu)=-\mathcal L.
\end{equation*}
More precisely, for every $\mu_n\to\mu$ in WFR and every
$\gamma^n\in\operatorname{OptSC}(\mu,\mu_n)$,
\begin{equation}
\label{eq:squared-wfr-differentiability}
    \frac{1}{2}\WFR^2(\mu_n,\hat\mu)
    =\frac{1}{2}\WFR^2(\mu,\hat\mu)
    -\left\langle \mathcal L,\barlog_{\gamma^n}
    \right\rangle_\mu
    +o\bigl(\WFR(\mu,\mu_n)\bigr).
\end{equation}
\end{theorem}

\begin{example}[Dirac configurations]\label{ex:squared-wfr-dirac}
Let $\mu=\delta_0\in\mathcal M(\R)$ and fix
$0<\rho<\pi/2<R$. We consider one differentiable case
and three obstructions. All uniqueness statements concern
source-normalized cone couplings (vertex--vertex mass being immaterial).

\begin{enumerate}
    \item[(a)] \emph{Transport of a single atom.}
    Let $\hat\mu=\delta_\rho$. The unique source-normalized
    optimal cone coupling is
    $\delta_{([0,1],[\rho,1])}$, and its barycentric logarithm is
    \begin{equation*}
        \mathcal L
        =(\sin(\rho),2(\cos(\rho)-1)).
    \end{equation*}
    In particular,
    \begin{equation*}
        \norm{\mathcal L}_\mu^2
        =2(1-\cos(\rho)) =\WFR^2(\mu,\hat\mu).
    \end{equation*}
    Thus $F_{\hat\mu}$ is differentiable at $\mu$.

    \item[(b)] \emph{Splitting of the source mass.}
    Let $\hat\mu=\tfrac12(\delta_{-\rho}+\delta_\rho)$.
    The unique source-normalized optimal cone coupling is
    \begin{equation*}
        \frac{1}{2}\delta_{([0,1],[-\rho,1])}
        +\frac{1}{2}\delta_{([0,1],[\rho,1])}.
    \end{equation*}
    Indeed, strict concavity of the square root in the
    correlation term forces the source mass to be divided
    equally between the two target atoms. The transport
    components cancel under barycentric averaging, giving
    \begin{equation*}
        \mathcal L=(0,2(\cos(\rho)-1)).
    \end{equation*}
    Consequently,
    \begin{equation*}
        \WFR^2(\mu,\hat\mu)-\norm{\mathcal L}_\mu^2
        =2(1-\cos(\rho)) -(1-\cos(\rho))^2
        =\sin^2(\rho)>0.
    \end{equation*}
    Hence $F_{\hat\mu}$ is not differentiable at $\mu$:
    uniqueness after source normalization does not exclude
    energy loss through barycentric averaging.

    \item[(c)] \emph{Non-uniqueness at the transport threshold.}
    Let $\hat\mu=\delta_{\pi/2}$. Both direct transport and
    separate annihilation and creation are optimal, as
    represented by the cone couplings
    \begin{equation*}
        A^T=\delta_{([0,1],[\pi/2,1])},
        \qquad
        A^R=\delta_{([0,1],\vertex)}
            +\delta_{(\vertex,[\pi/2,1])}.
    \end{equation*}
    Both have cost $2$, while their barycentric logarithms are
    \begin{equation*}
        \mathcal L^T=(1,-2),
        \qquad
        \mathcal L^R=(0,-2).
    \end{equation*}
    Thus $F_{\hat\mu}$ is not differentiable at $\mu$.
    In particular, the equality
    $\norm{\mathcal L^T}_\mu^2=\WFR^2(\mu,\hat\mu)$
    for one optimal logarithm does not suffice.

    \item[(d)] \emph{Creation at the target.}
    Let $\hat\mu=\delta_0+\delta_R$. An  optimal cone
    coupling is
    \begin{equation*}
        \delta_{([0,1],[0,1])}
        +\delta_{(\vertex,[R,1])}.
    \end{equation*}
    Its source-normalized part is unique, and its barycentric
    logarithm is $\mathcal L=0$, whereas
    \begin{equation*}
        \WFR^2(\mu,\hat\mu)-\norm{\mathcal L}_\mu^2=1.
    \end{equation*}
    The energy defect is entirely due to creation at $R$,
    and $F_{\hat\mu}$ is not differentiable at $\mu$.
\end{enumerate}
\end{example}

\begin{remark}[Absolutely continuous sources: a sufficient condition] \label{rem:squared-wfr-monge}
The normalized Monge condition can hold for absolutely
continuous sources. In this setting, the optimality and
uniqueness results of
\cite[Theorems~6.3, 6.6 and 7.20]{liero-mielke-savare2018}
yield this representation when the target has no
pure-creation part; see also
\cite[Corollary~3.5]{liero-mielke-savare2023}.
In particular, a  sufficient condition is that $\mu$ is absolutely
continuous with respect to the Lebesgue measure on $\R^d$ and
\begin{equation*}
    \hat\mu\left(
        \left\{
            y\in\R^d:
            \operatorname{dist}(y,\operatorname{supp}\mu)
            <\frac{\pi}{2}
        \right\}
    \right)
    =
    \hat\mu(\R^d).
\end{equation*}
The unique active transport map and marginal densities
then determine $q$, while the remaining source part is
sent to the vertex.
No positive lower bound on $q$ or uniform gap from the
transport cutoff is required.
\end{remark}

\begin{remark}[Absolutely continuous sources: target creation] \label{rem:squared-wfr-creation}
Absolute continuity of the source alone does not ensure
differentiability.
Let $\mu$ be a nonzero absolutely continuous measure
with compact support, choose $z\in\R^d$ such that
\begin{equation*}
    \operatorname{dist}(z,\operatorname{supp}\mu)
    >
    \frac{\pi}{2},
\end{equation*}
and set $\hat\mu=\mu+m\delta_z$ , with $m>0$.
For $\mu^\varepsilon=\mu+\varepsilon^2\delta_z$,
$\varepsilon>0$, separation of the supports gives
\begin{equation*}
    \WFR(\mu,\mu^\varepsilon)=\varepsilon,
    \qquad
    \WFR^2(\mu^\varepsilon,\hat\mu )-\WFR^2(\mu,\hat\mu )
    =
    -2\sqrt m\,\varepsilon+\varepsilon^2.
\end{equation*}
The semicoupling
\begin{equation*}
    \gamma^\varepsilon
    \coloneqq
    \left(
        (\operatorname{Id},\operatorname{Id})\pf\mu,
        (\operatorname{Id},\operatorname{Id})\pf\mu
        +\varepsilon^2\delta_{(z,z)}
    \right)
\end{equation*}
is optimal and satisfies
$\barlog_{\gamma^\varepsilon}=0$.
Hence the displayed linear variation cannot be represented
by a source tangent field.

The no-creation condition concerns only the fixed pair
$(\mu,\hat\mu)$: the perturbations $\mu_n$ in
Theorem~\ref{thm:squared-wfr-differentiability} may contain
creation mass, whose contribution is controlled by
Lemma~\ref{lemma:squared-wfr-first-variation}.
\end{remark}

\subsection{Stability and first variation}
\label{subsec:squared-wfr-preliminaries}

We first show that, under the source-normalized Monge condition,
asymptotically optimal cone couplings have source logarithms
converging to $\mathcal L$, with vanishing target creation.
We then use this stability to obtain an integrated first-variation
formula for small perturbations. These two lemmas are the main
ingredients in the sufficiency part of Theorem~\ref{thm:squared-wfr-differentiability}.

The stability argument also involves non-reduced cone pairs.
We therefore use the kernel $K_y$ from Section~\ref{subsec:wfr-geodesics} to extend the normalized cone logarithm of
Remark~\ref{rem:cone-log-normalization} to arbitrary cone pairs.
By Lemma~\ref{lemma:cone-geodesic-velocities}, for
$\y_i=[x_i,r_i]$ with $r_i>0$ and $\y_j=[x_j,r_j]$,
\begin{equation}\label{eq:squared-wfr-normalized-cone-log}
    \mathcal L_{\cone}(\y_i,\y_j)
    =
    \left(
        \frac{r_j}{r_i}\nabla K_{x_j}(x_i),
        2\left(\frac{r_j}{r_i}K_{x_j}(x_i)-1\right)
    \right).
\end{equation}
When $\y_j=\vertex$, this means
$\mathcal L_{\cone}(\y_i,\vertex)=(0,-2)$.
At a source vertex we retain the convention
$\mathcal L_{\cone}(\vertex,\y_j)=0$.
The identity
$|K_y|^2+\norm{\nabla K_y}^2=1$ gives, for all cone pairs,
\begin{equation}\label{eq:squared-wfr-cone-log-identities}
    \dc^2(\y_i,\y_j)
    =
    r_i^2\norm{\mathcal L_{\cone}(\y_i,\y_j)}_{\mathrm{WFR}}^2
    +\mathbf 1_{\{r_i=0\}}r_j^2,
    \qquad
    \norm{\mathcal L_{\cone}(\y_i,\y_j)+(0,2)}_{\mathrm{WFR}}
    =\frac{r_j}{r_i}\quad(r_i>0).
\end{equation}
Thus the source logarithmic energy is bounded by the cone cost;
its deficit at the source vertex is exactly the target-creation cost.
Barycentric averaging can produce a further loss.

The identification with semicoupling logarithms remains
restricted to reduced pairs, with the choice $u=1$ at the cutoff;
the identities above also hold for non-reduced pairs.

\begin{lemma}[Stability of asymptotically optimal cone couplings]
\label{lemma:squared-wfr-stability}
Assume \eqref{eq:squared-wfr-monge-assumption}, and let
$E_n\in\Gamma_{\cone}(\mu,\hat\mu)$ satisfy
\begin{equation*}
    c_n\coloneqq
    \int_{\cone^2}\dc^2(\y_0,\y_2)\,dE_n
    \longrightarrow
    \WFR^2(\mu,\hat\mu).
\end{equation*}
Then
\begin{equation}\label{eq:squared-wfr-stability-1}
    \int_{\{r_0>0\}}
    r_0^2
    \norm{
        \mathcal L_{\cone}(\y_0,\y_2)-\mathcal L(x_0)
    }_{\mathrm{WFR}}^2\,dE_n
    \longrightarrow0.
\end{equation}
and
\begin{equation}\label{eq:squared-wfr-stability-2}
    \int_{\{r_0=0\}}r_2^2\,dE_n
    \longrightarrow0.
\end{equation}
\end{lemma}

\begin{proof}
Discard vertex--vertex mass and apply the dilation
\begin{equation*}
    \vartheta([x_0,r_0],[x_2,r_2])=(r_0^2+r_2^2)^{1/2}.
\end{equation*}
All quantities in the statement are unchanged due to Lemma \ref{lemma:cone-dilation-normalization}, so we
may assume without loss of generality that
\begin{equation*}
    r_0^2+r_2^2=1
    \quad E_n\text{-a.e.},
    \qquad
    E_n(\cone^2)=\mu(\R^d)+\hat\mu(\R^d).
\end{equation*}
If both measures vanish, the conclusion is immediate.

The homogeneous-marginal tightness criterion
\cite[Lemma~7.3]{liero-mielke-savare2018}
and the fixed total mass give narrow relative compactness.
Every narrow limit is concentrated on
$\{r_0^2+r_2^2=1\}$ and has homogeneous marginals
$\mu$ and $\hat\mu$.
Indeed, for $\varphi\in C_b(\R^d)$, the functions
$(r_i^2\wedge1)\varphi(x_i)$ extend continuously by zero
at the corresponding vertex and agree with
$r_i^2\varphi(x_i)$ on this set.
The cone cost is also continuous and bounded there.
Consequently, every narrow limit belongs to
$\operatorname{Opt\Gamma}_{\cone}(\mu,\hat\mu)$.

For $\mathcal V=(\xi,\zeta)\in H(\mu)$, set
\begin{equation*}
    J_n(\mathcal V)
    \coloneqq
    \int_{\{r_0>0\}}
    r_0^2
    \left\langle
        \mathcal L_{\cone}(\y_0,\y_2),\mathcal V(x_0)
    \right\rangle_{\mathrm{WFR}}\,dE_n.
\end{equation*}
If $\mathcal V$ is bounded and continuous, the integrand equals
\begin{equation*}
    r_0r_2\nabla K_{x_2}(x_0)\cdot \xi(x_0)
    +
    \frac12
    \bigl(r_0r_2K_{x_2}(x_0)-r_0^2\bigr)\zeta(x_0).
\end{equation*}
It is bounded and continuous on the normalized set,
including at the vertices.
Every narrow limit is optimal, so
\eqref{eq:squared-wfr-monge-assumption} identifies its
normalized source logarithm with $\mathcal L(x_0)$.
It follows that
\begin{equation*}
    J_n(\mathcal V)\longrightarrow\langle \mathcal L,\mathcal V\rangle_\mu.
\end{equation*}
Here convergence holds for the whole sequence because
every subsequence admits a further narrowly convergent
subsequence giving the same limit.

By Cauchy--Schwarz and
\eqref{eq:squared-wfr-cone-log-identities},
\begin{equation*}
    |J_n(\mathcal V)|\leq\sqrt{c_n}\,\norm{\mathcal V}_\mu,
\end{equation*}
in particular, the linear mappings $J_n$ are uniformly bounded. Since bounded continuous fields are dense in $H(\mu)$,
the preceding convergence extends to every
$\mathcal V\in H(\mu)$.
In particular,
\begin{equation*}
    J_n(\mathcal L)\longrightarrow\norm{\mathcal L}_\mu^2.
\end{equation*}

Finally, using the fixed source homogeneous marginal
and \eqref{eq:squared-wfr-cone-log-identities}, we obtain
the exact identity
\begin{equation*}
    \begin{aligned}
        \int_{\{r_0>0\}}
        r_0^2
        \norm{
            \mathcal L_{\cone}(\y_0,\y_2)-\mathcal L(x_0)
        }_{\mathrm{WFR}}^2\,dE_n
        +
        \int_{\{r_0=0\}}r_2^2\,dE_n
        =
        c_n+\norm{\mathcal L}_\mu^2-2J_n(\mathcal L).
    \end{aligned}
\end{equation*}
The right-hand side converges to zero by
\eqref{eq:squared-wfr-reference-energy}. Since both terms on the left hand side are nonnegative the claim follows.
\end{proof}

We next establish the first-variation estimate used
for both competitors.
Only the perturbation pairs are assumed to be reduced;
the pairs connecting the source to the fixed target
may be non-reduced.

\begin{lemma}[Integrated first variation]
\label{lemma:squared-wfr-first-variation}
Let $\mathcal L$ be defined by
\eqref{eq:squared-wfr-reference-log}, and let
$A_n\in\mathcal M_2(\cone^3)$ have homogeneous marginals $\mu,\mu_n,\hat\mu$ at endpoints $0,1,2$, respectively.
Assume that the $(0,1)$-pairs are reduced and
\begin{equation*}
    \varepsilon_n^2
    \coloneqq
    \int_{\cone^3}\dc^2(\y_0,\y_1)\,dA_n
    \longrightarrow0.
\end{equation*}
Suppose also that
\begin{equation}\label{eq:squared-wfr-first-variation-assumptions}
    \begin{aligned}
        \delta_n^2
        &\coloneqq
        \int_{\{r_0>0\}}
        r_0^2\norm{
            \mathcal L_{\cone}(\y_0,\y_2)-\mathcal L(x_0)
        }_{\mathrm{WFR}}^2\,dA_n
        \longrightarrow0,\\
        m_n
        &\coloneqq
        \int_{\{r_0=0\}}r_2^2\,dA_n
        \longrightarrow0.
    \end{aligned}
\end{equation}
Then
\begin{equation}\label{eq:squared-wfr-integrated-first-variation}
    \begin{aligned}
        \int_{\cone^3}
        \bigl(
            \dc^2(\y_1,\y_2)-\dc^2(\y_0,\y_2)
        \bigr)\,dA_n
        =
        -2\int_{\{r_0>0\}}
        r_0^2
        \left\langle
            \mathcal L(x_0),\mathcal L_{\cone}(\y_0,\y_1)
        \right\rangle_{\mathrm{WFR}}\,dA_n
        +o(\varepsilon_n).
    \end{aligned}
\end{equation}
\end{lemma}

\begin{proof}
If $\varepsilon_n=0$, then $\y_0=\y_1$ for $A_n$-a.e.
triple, so both the cost difference and the logarithmic pairing
vanish. We may therefore restrict to $\varepsilon_n>0$.

On $\{r_0>0\}$, put
\begin{equation*}
    q_1\coloneqq\frac{r_1}{r_0},
    \qquad
    q_2\coloneqq\frac{r_2}{r_0},
    \qquad
    \alpha\coloneqq\frac{\dc(\y_0,\y_1)}{r_0}.
\end{equation*}
If $r_i=0$, we choose the spatial representative $x_i=x_0$.
Reducedness then gives $\norm{x_1-x_0}\leq\pi/2$ and
\begin{gather*}
    \mathcal L_{\cone}(\y_0,\y_1)
    =
    \bigl(
        q_1{\mathbf{sin}}(x_1-x_0),
        2(q_1\cos(\norm{x_1-x_0})-1)
    \bigr),
    \\
    \alpha^2
    =
    (q_1-1)^2+2q_1(1-\cos(\norm{x_1-x_0}))
    =
    \norm{\mathcal L_{\cone}(\y_0,\y_1)}_{\mathrm{WFR}}^2.
\end{gather*}
Expanding the cone costs and the scalar product yields
\begin{equation}\label{eq:squared-wfr-pointwise-identity}
\begin{split}
    \frac{
        \dc^2(\y_1,\y_2)-\dc^2(\y_0,\y_2)
    }{r_0^2}
    +2\left\langle
        \mathcal L_{\cone}(\y_0,\y_2),
        \mathcal L_{\cone}(\y_0,\y_1)
    \right\rangle_{\mathrm{WFR}}
    =\alpha^2-2q_1q_2B,
\end{split}
\end{equation}
where
\begin{equation*}
    B\coloneqq
    \cos(\norm{x_2-x_1}\wedge\pi)
    -K_{x_1}(x_0)K_{x_2}(x_0)
    -\nabla K_{x_1}(x_0)\cdot\nabla K_{x_2}(x_0).
\end{equation*}

To estimate $B$, write
\begin{equation*}
\begin{split}
    B={}&
    K_{x_2}(x_1)-K_{x_2}(x_0)
    -\nabla K_{x_2}(x_0)\cdot(x_1-x_0)\\
    &+\bigl(1-K_{x_1}(x_0)\bigr)K_{x_2}(x_0)\\
    &+\nabla K_{x_2}(x_0)\cdot
    \bigl(x_1-x_0-\nabla K_{x_1}(x_0)\bigr).
\end{split}
\end{equation*}
Since $\operatorname{Lip}(\nabla K_{x_2})\leq1$,
the Taylor expansion of $K_{x_2}$ bounds the first line in absolute value by
$\frac12\norm{x_1-x_0}^2$.
Reducedness gives $\norm{x_1-x_0}\leq\pi/2$, so
\begin{gather*}
    0\leq1-K_{x_1}(x_0)
    =
    1-\cos(\norm{x_1-x_0})
    \leq\frac12\norm{x_1-x_0}^2,\\
    \norm{x_1-x_0-\nabla K_{x_1}(x_0)}
    =
    \norm{x_1-x_0}-\sin(\norm{x_1-x_0})
    \leq\frac{\pi}{12}\norm{x_1-x_0}^2.
\end{gather*}
Together with $|K_{x_2}|\leq1$ and
$\norm{\nabla K_{x_2}}\leq1$, these estimates yield
$|B|\leq C\norm{x_1-x_0}^2$.

Moreover, Cauchy--Schwarz and \eqref{eq:kernel-bounds} give
\begin{equation*}
\begin{split}
    &\left|
        K_{x_1}(x_0)K_{x_2}(x_0)
        +\nabla K_{x_1}(x_0)\cdot\nabla K_{x_2}(x_0)
    \right|\leq 1.
\end{split}
\end{equation*}
Therefore, we obtain $|B|\leq2$.

The identity defining $\alpha$ implies
\begin{equation*}
    \alpha^2
    \geq2q_1(1-\cos(\norm{x_1-x_0}))
    \geq\frac{4q_1}{\pi^2}\norm{x_1-x_0}^2.
\end{equation*}
Hence the quadratic estimate for $B$ gives
$q_1|B|\leq C\alpha^2$.
On the other hand, $\alpha^2\geq(q_1-1)^2$ implies
$q_1\leq1+\alpha$. Thus, on $\{\alpha\geq1\}$,
the bound $|B|\leq2$ gives
\begin{equation*}
    q_1|B|\leq2(1+\alpha)\leq4\alpha.
\end{equation*}
Using the quadratic bound for $\alpha\leq1$ and the linear
bound for $\alpha\geq1$, we conclude that
\begin{equation*}
    q_1|B|\leq C\alpha\min\{\alpha,1\}.
\end{equation*}

Therefore
\begin{equation}\label{eq:squared-wfr-pointwise-remainder}
\begin{split}
    \left|
        \frac{
            \dc^2(\y_1,\y_2)-\dc^2(\y_0,\y_2)
        }{r_0^2}
        +2\left\langle
            \mathcal L_{\cone}(\y_0,\y_2),
            \mathcal L_{\cone}(\y_0,\y_1)
        \right\rangle_{\mathrm{WFR}}
    \right|
    \leq
    \alpha^2+Cq_2\alpha\min\{\alpha,1\}.
\end{split}
\end{equation}

Set
\begin{equation*}
    d\lambda_n
    \coloneqq
    r_0^2\mathbf1_{\{r_0>0\}}\,dA_n.
\end{equation*}
Then $(x_0)\pf\lambda_n=\mu$ and
$\int\alpha^2\,d\lambda_n\leq\varepsilon_n^2$.
By \eqref{eq:squared-wfr-cone-log-identities} and the
definition of $\mathcal L$, the reverse triangle inequality
gives
\begin{equation*}
\begin{split}
    |q_2-q(x_0)|
    &=
    \left|
        \norm{\mathcal L_{\cone}(\y_0,\y_2)+(0,2)}_{\mathrm{WFR}}
        -
        \norm{\mathcal L(x_0)+(0,2)}_{\mathrm{WFR}}
    \right|\\
    &\leq
    \norm{
        \mathcal L_{\cone}(\y_0,\y_2)-\mathcal L(x_0)
    }_{\mathrm{WFR}}.
\end{split}
\end{equation*}
Define
\begin{equation*}
    \eta_n
    \coloneqq
    (x_0)\pf\bigl(\min\{\alpha^2,1\}\lambda_n\bigr).
\end{equation*}
Since $(x_0)\pf\lambda_n=\mu$, we have
\begin{equation*}
    0\leq\eta_n\leq\mu,
    \qquad
    \eta_n(\R^d)
    \leq\int\alpha^2\,d\lambda_n
    \leq\varepsilon_n^2
    \longrightarrow0.
\end{equation*}
As $q^2\in L^1(\mu)$, Lemma~\ref{lemma:vanishing-dominated-measures} yields
$\int q^2\,d\eta_n\to0$. Consequently,
\begin{equation*}
    \int q_2^2\min\{\alpha^2,1\}\,d\lambda_n
    \leq
    2\delta_n^2+2\int q^2\,d\eta_n
    \longrightarrow0.
\end{equation*}
Cauchy--Schwarz therefore gives
\begin{equation*}
    \int q_2\alpha\min\{\alpha,1\}\,d\lambda_n
    \leq
    \varepsilon_n
    \left(
        \int q_2^2\min\{\alpha^2,1\}\,d\lambda_n
    \right)^{1/2}
    =
    o(\varepsilon_n).
\end{equation*}

Multiplying \eqref{eq:squared-wfr-pointwise-remainder}
by $r_0^2$ and integrating now yields
\begin{equation*}
\begin{split}
    &\left|
        \int_{\{r_0>0\}}
        \bigl(\dc^2(\y_1,\y_2)-\dc^2(\y_0,\y_2)\bigr)\,dA_n
        +2\int
        \left\langle
            \mathcal L_{\cone}(\y_0,\y_2),
            \mathcal L_{\cone}(\y_0,\y_1)
        \right\rangle_{\mathrm{WFR}}\,d\lambda_n
    \right|\\
    &\qquad
    \leq
    \varepsilon_n^2
    +C\int q_2\alpha\min\{\alpha,1\}\,d\lambda_n
    =
    o(\varepsilon_n).
\end{split}
\end{equation*}
Replacing $\mathcal L_{\cone}(\y_0,\y_2)$ by
$\mathcal L(x_0)$ in this pairing changes the expression
by at most
\begin{equation*}
\begin{split}
    &2\left|
        \int
        \left\langle
            \mathcal L_{\cone}(\y_0,\y_2)-\mathcal L(x_0),
            \mathcal L_{\cone}(\y_0,\y_1)
        \right\rangle_{\mathrm{WFR}}\,d\lambda_n
    \right|\\
    &\qquad
    \leq
    2\delta_n
    \left(
        \int
        \norm{\mathcal L_{\cone}(\y_0,\y_1)}_{\mathrm{WFR}}^2
        \,d\lambda_n
    \right)^{1/2}
    \leq
    2\delta_n\varepsilon_n
    =
    o(\varepsilon_n).
\end{split}
\end{equation*}

Finally, on $\{r_0=0\}$, the cost difference is
\begin{equation*}
    \dc^2(\y_1,\y_2)-\dc^2(\y_0,\y_2)
    =
    r_1^2-2r_1r_2K_{x_2}(x_1),
\end{equation*}
where the product term is interpreted as zero if
$r_1r_2=0$.
Since $\dc^2(\y_0,\y_1)=r_1^2$ on this set,
\begin{equation*}
    \int_{\{r_0=0\}}r_1^2\,dA_n\leq\varepsilon_n^2.
\end{equation*}
Using $|K_{x_2}|\leq1$, Cauchy--Schwarz, and the definition
of $m_n$, we conclude that
\begin{equation*}
\begin{split}
    \int_{\{r_0=0\}}
    \left|
        \dc^2(\y_1,\y_2)-\dc^2(\y_0,\y_2)
    \right|\,dA_n&
    \leq
    \int_{\{r_0=0\}}r_1^2\,dA_n
    +
    2\left(
        \int_{\{r_0=0\}}r_1^2\,dA_n
    \right)^{1/2}\sqrt{m_n}\\
    &
    \leq
    \varepsilon_n^2+2\varepsilon_n\sqrt{m_n}
    =
    o(\varepsilon_n).
\end{split}
\end{equation*}
Combining the estimates on $\{r_0>0\}$ and $\{r_0=0\}$
proves \eqref{eq:squared-wfr-integrated-first-variation}.
\end{proof}

\subsection{Proof of differentiability and characterization}

\begin{proof}[Proof of Theorem~\ref{thm:squared-wfr-differentiability}]
We first prove $(i)\Rightarrow(ii)$. Set
$D\coloneqq \WFR^2(\mu,\hat\mu)$ and suppose that $F_{\hat\mu}$ is WFR
differentiable at $\mu$, with gradient $\mathcal V$.

The descending metric slope of $F_{\hat\mu}$ at $\mu$ equals $\sqrt{D}$.
Indeed, the triangle inequality gives
\begin{equation*}
    \frac{\max\{F_{\hat{\mu}}(\mu)-F_{\hat\mu}(\nu),0\}}
    {\WFR(\mu,\nu)}
    \le \frac{\sqrt{D}+\WFR(\nu,\hat\mu)}{2},
    \qquad \nu\ne\mu.
\end{equation*}
For $D>0$, any constant-speed geodesic
$(\bfmu_t)_{t\in[0,1]}$ from $\mu$ to $\hat\mu$ satisfies
\begin{equation*}
    \frac{F_{\hat{\mu}}(\mu)-F_{\hat\mu}(\bfmu_t)}{\WFR(\mu,\bfmu_t)}
    =\frac{(2-t)}{2}\sqrt D\longrightarrow\sqrt D.
\end{equation*}
For $D=0$, the slope is zero since $F_{\hat\mu}$ attains its
minimum at $\mu$.
Proposition \ref{propo:metric-slope} therefore gives
\begin{equation*}
    \norm{\mathcal V}_\mu=\sqrt{D}.
\end{equation*}

Let $\mathcal L$ be any regular barycentric logarithm from $\mu$ to $\hat\mu$,
and let $(\bfmu_t)_{t\in[0,1]}$ be an associated geodesic.
Since $F_{\hat\mu}(\bfmu_t)=(1-t)^2D/2$,
Proposition~\ref{propo:wfr-geodesic-chain-rule} gives
\begin{equation*}
    \frac{(1-t)^2D-D}{2}
    =
    t\langle \mathcal V,\mathcal L\rangle_\mu+o(t).
\end{equation*}
Consequently, $\langle \mathcal V,\mathcal L\rangle_\mu=-D$.
Using $\norm{\mathcal L}_\mu\le\sqrt{D}$, we obtain
\begin{equation*}
    D
    =-\langle \mathcal V,\mathcal L\rangle_\mu
    \le\norm{\mathcal V}_\mu\norm{\mathcal L}_\mu
    \le D.
\end{equation*}
Equality in Cauchy--Schwarz implies
\begin{equation*}
    \mathcal L=-\mathcal V,
    \qquad \norm{\mathcal L}_\mu=\sqrt D.
\end{equation*}
Since $\mathcal L$ was arbitrary, this proves (ii).
If $D=0$, the bound $\norm{\mathcal L}_\mu\le \sqrt D$ directly gives
the same conclusion.

We next prove $(ii)\Rightarrow(iii)$. Fix
$A\in\operatorname{Opt\Gamma}_{\cone}(\mu,\hat\mu)$ and write
$\y_0=[x_0,r_0]$, $\y_2=[x_2,r_2]$.
If $A(\cone^2)=0$ the
claim is trivial; otherwise we may assume that $A$ is a
probability measure, by constant dilation, which preserves the homogeneous marginals,
the cost, and all source-weighted logarithmic integrals
below.
Let $\boldsymbol\Pi$ be a probability measure on
$\Geo(\cone)$ that generates $A$, which exists by
Theorem~\ref{thm:wfr-geodesics-cone}.
Let
$Q(d\y_0,d\y_2)=r_0^2 A(d\y_0,d\y_2)$.
Its spatial source marginal is $(x_0)\pf Q=\mu$, so we
may disintegrate it as
$Q=\int_{\R^d}Q_x\,d\mu(x)$.
Let
\begin{equation*}
    \mathcal L_A(x)\coloneqq
    \int_{\cone^2}
    \mathcal L_{\cone}(\y_0,\y_2)\,dQ_x(\y_0,\y_2)
    =\mathcal V^{\boldsymbol\Pi}(x)
\end{equation*}
be the corresponding barycenter, where the last equality
follows from $(e_0,e_1)\pf\boldsymbol\Pi=A$.

Proposition~\ref{propo:cone-barycentric-initial-derivative}
and Theorem~\ref{propo:tangency} show that $\mathcal L_A$ is a regular barycentric logarithm
from $\mu$ to $\hat\mu$.
Thus, by (ii),
\begin{equation*}
    \mathcal L_A=\mathcal L,
    \qquad
    \norm{\mathcal L_A}_\mu=\WFR(\mu,\hat\mu).
\end{equation*}
Therefore, there is no energy loss due to barycentric
averaging. The cone energy identity gives
\begin{equation*}
    \int_{\{r_0>0\}}r_0^2
    \norm{\mathcal L_{\cone}(\y_0,\y_2)}_{\mathrm{WFR}}^2\,dA
    =
    \WFR(\mu,\hat\mu)^2
    -\int_{\{r_0=0\}}r_2^2\,dA.
\end{equation*}
Moreover, the definition of the barycenter and the
homogeneous marginal constraint imply
\begin{equation*}
\begin{split}
    \int_{\{r_0>0\}}r_0^2
    \left\langle
        \mathcal L_{\cone}(\y_0,\y_2),\mathcal L_A(x_0)
    \right\rangle_{\mathrm{WFR}}\,dA
    &=\norm{\mathcal L_A}_\mu^2,\\
    \int_{\{r_0>0\}}r_0^2
    \norm{\mathcal L_A(x_0)}_{\mathrm{WFR}}^2\,dA
    &=\norm{\mathcal L_A}_\mu^2.
\end{split}
\end{equation*}

Expanding the square therefore yields the exact norm deficit
\begin{equation}
\label{eq:squared-wfr-logarithm-deficit}
\begin{split}
    \WFR(\mu,\hat\mu)^2-\norm{\mathcal L_A}_\mu^2
    ={}&
    \int_{\{r_0>0\}}r_0^2
    \norm{\mathcal L_{\cone}(\y_0,\y_2)
          -\mathcal L_A(x_0)}_{\mathrm{WFR}}^2\,dA+\int_{\{r_0=0\}}r_2^2\,dA.
\end{split}
\end{equation}
Its left-hand side vanishes.
Both terms on the right are nonnegative,
and $r_0^2>0$ on $\{r_0>0\}$.
Hence every optimal $A$ satisfies
\begin{equation*}
    \mathcal L_{\cone}(\y_0,\y_2)=\mathcal L(x_0)
    \quad A\text{-a.e. on }\{r_0>0\},
    \qquad
    \int_{\{r_0=0\}}r_2^2\,dA=0.
\end{equation*}

We recover common Borel maps $T,q$ from $\mathcal L=(\xi,\zeta)$. Set
\begin{equation*}
    a(x)\coloneqq 1+\frac{\zeta(x)}2,
    \qquad
    q(x)\coloneqq\sqrt{a(x)^2+\norm{\xi(x)}^2}.
\end{equation*}
Optimal cone pairs are reduced. Therefore, on
$\{r_0>0,r_2>0\}$, writing
$\norm{x_2-x_0}\le\pi/2$, we have
\begin{equation*}
    a(x_0)=\frac{r_2}{r_0}\cos(\norm{x_2-x_0}),
    \qquad
    \xi(x_0)=\frac{r_2}{r_0}
    \mathbf{\sin}(x_2-x_0)
\end{equation*}
On $\{r_0>0,r_2=0\}$, both $a(x_0)$ and $\xi(x_0)$ vanish.
Consequently,
\begin{equation*}
    a\ge0\quad\mu\text{-a.e.},
    \qquad
    q(x_0)=\frac{r_2}{r_0}
    \quad A\text{-a.e. on }\{r_0>0\}.
\end{equation*}
After modifying the representative on a $\mu$-null set,
assume $a\ge0$ everywhere and define
\begin{equation*}
    T(x)\coloneqq
    \begin{cases}
        x+\displaystyle
        \frac{\arccos(a(x)/q(x))}{\norm{\xi(x)}}\xi(x),
        &\xi(x)\ne0,\\
        x,&\xi(x)=0.
    \end{cases}
\end{equation*}
This defines a Borel map: when $\xi(x)\ne0$, the identity
$q(x)^2=a(x)^2+\norm{\xi(x)}^2$ ensures that $q(x)>0$
and $a(x)/q(x)\in[0,1]$.

To verify that $T$ recovers the target position, consider
a pair with $r_0>0$ and $\xi(x_0)\ne0$.
Then $r_2=q(x_0)r_0>0$, and the preceding identities give
\begin{equation*}
    \frac{a(x_0)}{q(x_0)}=\cos(\norm{x_2-x_0}),
    \qquad
    \frac{\xi(x_0)}{\norm{\xi(x_0)}}
    =\frac{x_2-x_0}{\norm{x_2-x_0}}.
\end{equation*}
Since reducedness gives $\norm{x_2-x_0}\in(0,\pi/2]$,
we have $\arccos(a(x_0)/q(x_0))=\norm{x_2-x_0}$, hence
$T(x_0)=x_2$.

If $q(x_0)>0$ and $\xi(x_0)=0$, then
$q(x_0)\sin(\norm{x_2-x_0})=\norm{\xi(x_0)}=0$.
Since $\norm{x_2-x_0}\in[0,\pi/2]$, this implies $\norm{x_2-x_0}=0$,
so again $T(x_0)=x_0=x_2$.

If $q(x_0)=0$, the target is the vertex and the value of
$T(x_0)$ is irrelevant.
Thus every optimal $A$ satisfies
\begin{equation*}
    \y_2=[T(x_0),q(x_0)r_0]
    \quad A\text{-a.e. on }\{r_0>0\}.
\end{equation*}
The maps are common to all optimizers because they depend only
on $\mathcal L$. Together with the vanishing creation term,
this proves \eqref{eq:squared-wfr-monge-assumption}, hence (iii).

Finally, we prove $(iii)\Rightarrow(i)$ and the stated expansion.
Assume (iii), and use the corresponding maps $T,q$, optimal
cone plan $A^\star$, and field $\mathcal L$ defined above.
In particular, $\mathcal L\in\mathrm{T}_\mu\M$ and
$\norm{\mathcal L}_\mu^2=\WFR^2(\mu,\hat\mu)$.

Fix $\mu_n\to\mu$ in WFR and arbitrary
$\gamma^n\in\operatorname{OptSC}(\mu,\mu_n)$. Write
\begin{equation*}
    d_n\coloneqq\WFR(\mu,\mu_n),
    \qquad
    D\coloneqq \WFR^2(\mu,\hat\mu),
    \qquad D_n\coloneqq \WFR^2(\mu_n,\hat\mu).
\end{equation*}
Terms with $d_n=0$ are immediate and can be omitted.
By Proposition~\ref{propo:semicoupling-log-realization}, with the default choice $u=1$, choose a
reduced optimal cone coupling $B_n$ representing the prescribed
barycentric logarithm $\barlog_{\gamma^n}$. Thus
\begin{equation}
\label{eq:squared-wfr-prescribed-pairing}
\begin{split}
    \int_{\cone^2}\dc^2(\y_0,\y_1)\,dB_n
    &=d_n^2,\\
    \int_{\{r_0>0\}}r_0^2
    \left\langle
        \mathcal L(x_0),\mathcal L_{\cone}(\y_0,\y_1)
    \right\rangle_{\mathrm{WFR}}\,dB_n
    &=
    \left\langle
        \mathcal L,\barlog_{\gamma^n}
    \right\rangle_\mu.
\end{split}
\end{equation}

These identities are preserved by the dilations and
vertex--vertex modifications in
Lemma~\ref{lemma:squared-wfr-gluing}.

\emph{Upper bound.}
Glue $B_n$ and $A^\star$ at their common source $\mu$,
using Lemma~\ref{lemma:squared-wfr-gluing} after relabelling
the coordinates, and denote the resulting triple plan by
$A_n^+$.
Its $(0,2)$ marginal is a dilation of $A^\star$, up to
vertex--vertex mass. Consequently,
\begin{equation*}
    \mathcal L_{\cone}(\y_0,\y_2)=\mathcal L(x_0)
    \quad A_n^+\text{-a.e. on }\{r_0>0\},
    \qquad
    \int_{\{r_0=0\}}r_2^2\,dA_n^+=0.
\end{equation*}
Moreover,
\begin{equation*}
    \int_{\cone^3}\dc^2(\y_0,\y_1)\,dA_n^+=d_n^2,
    \qquad
    \int_{\cone^3}\dc^2(\y_0,\y_2)\,dA_n^+=D.
\end{equation*}
The $(1,2)$ marginal is admissible for $(\mu_n,\hat\mu)$.
Applying Lemma~\ref{lemma:squared-wfr-first-variation}
and \eqref{eq:squared-wfr-prescribed-pairing} gives
\begin{equation}
\label{eq:squared-wfr-upper-bound}
\begin{split}
    \frac{D_n-D}{2}
    \le\frac{1}{2}
    \int_{\cone^3}
    \bigl(\dc^2(\y_1,\y_2)-\dc^2(\y_0,\y_2)\bigr)\,dA_n^+=
    -\left\langle
        \mathcal L,\barlog_{\gamma^n}
    \right\rangle_\mu
    +o(d_n).
\end{split}
\end{equation}

\emph{Lower bound.}
Choose
$C_n\in\operatorname{Opt\Gamma}_{\cone}(\mu_n,\hat\mu)$
and glue it to $B_n$ at their common endpoint $\mu_n$,
again using Lemma~\ref{lemma:squared-wfr-gluing}.
Denote the resulting triple plan by $A_n^-$ and set
\begin{equation*}
    E_n:=(\pi_{02})\pf A_n^-.
\end{equation*}
Then $E_n\in\Gamma_{\cone}(\mu,\hat\mu)$, and the prescribed
pair costs satisfy
\begin{equation*}
    \int_{\cone^3}\dc^2(\y_0,\y_1)\,dA_n^-=d_n^2,
    \qquad
    \int_{\cone^3}\dc^2(\y_1,\y_2)\,dA_n^-=D_n.
\end{equation*}
The cone triangle inequality and Minkowski's inequality yield
\begin{equation*}
    \sqrt D
    \le
    \left(
        \int_{\cone^2}\dc^2(\y_0,\y_2)\,dE_n
    \right)^{1/2}
    \le d_n+\sqrt{D_n}.
\end{equation*}
Since
\begin{equation*}
    \bigl|\sqrt{D_n}-\sqrt D\bigr|\le d_n\longrightarrow0,
\end{equation*}
the couplings $E_n$ are asymptotically optimal:
\begin{equation*}
    \int_{\cone^2}\dc^2(\y_0,\y_2)\,dE_n\longrightarrow D.
\end{equation*}
Lemma~\ref{lemma:squared-wfr-stability} therefore gives
\begin{equation*}
    \int_{\{r_0>0\}}r_0^2
    \norm{\mathcal L_{\cone}(\y_0,\y_2)-\mathcal L(x_0)}_{\mathrm{WFR}}^2
    \,dA_n^-\longrightarrow0,
    \qquad
    \int_{\{r_0=0\}}r_2^2\,dA_n^-\longrightarrow0.
\end{equation*}
Thus Lemma~\ref{lemma:squared-wfr-first-variation}
applies to $A_n^-$ as well.
Using admissibility of $E_n$ and
\eqref{eq:squared-wfr-prescribed-pairing}, we obtain
\begin{equation}
\label{eq:squared-wfr-lower-bound}
\begin{split}
    \frac{D_n-D}{2}
    \ge\frac{1}{2}
    \int_{\cone^3}
    \bigl(\dc^2(\y_1,\y_2)-\dc^2(\y_0,\y_2)\bigr)\,dA_n^-=
    -\left\langle
        \mathcal L,\barlog_{\gamma^n}
    \right\rangle_\mu
    +o(d_n).
\end{split}
\end{equation}

Combining \eqref{eq:squared-wfr-upper-bound} and
\eqref{eq:squared-wfr-lower-bound} proves
\eqref{eq:squared-wfr-differentiability}.
Since $\mathcal L\in\mathrm{T}_\mu\M$, this is WFR differentiability
at $\mu$, with
\begin{equation*}
    \nabla_{\WFR}F_{\hat\mu}(\mu)=-\mathcal L.
\end{equation*}
\end{proof}

\section{Restriction to probability measures}

In this last section, we consider the restriction of our calculus to the space of probability measures.

A time-dependent tangent vector
$\mathcal V_t=(\xi_t,\zeta_t)$ acts on a narrowly continuous curve of measures through the continuity equation with reaction
\begin{equation}\label{eq:continuity+reaction}
\partial_t\mu_t
+
\operatorname{div}(\xi_t\mu_t) =
\zeta_t\mu_t,
\qquad t\in[0,1],
\end{equation}
with the standard integrability conditions
\begin{equation*}
    \int_0^1\int_{\R^d}\norm{\xi_t(x)}^2\,d\mu_t(x)dt<+\infty,\quad
    \int_0^1\int_{\R^d}|\zeta_t(x)|^2\,d\mu_t(x)dt<+\infty.
\end{equation*}
In particular, testing against compactly supported cutoffs converging to $\varphi\equiv1$ yields
\begin{equation}
\frac{d}{dt}\mu_t(\R^d) = \int_{\R^d}\zeta_t(x)\,d\mu_t(x)
\end{equation}
for almost every $t$, and hence a curve taking values in $\prob(\R^d)$ must satisfy
\begin{equation}\label{eq:mass-invariance}
\int_{\R^d}\zeta_t(x)\,d\mu_t(x)=0
\end{equation}
for almost every $t$. This motivates defining, for $\mu\in\prob(\R^d)$,
\begin{equation}\label{eq:tangent-space-P}
T_\mu\prob(\R^d)
\coloneqq
\left\{
(\xi,\zeta)\in T_\mu\M
\colon
\int_{\R^d}\zeta(x)\,d\mu(x)=0
\right\}.
\end{equation}
Since the map $(\xi,\zeta)
\longmapsto
\int_{\R^d}\zeta\,d\mu$ is continuous on $T_\mu\M$, the space $T_\mu\prob(\R^d)$ is a closed subspace of $T_\mu\M$. Moreover, one has $(0,c)\in T_\mu\M$ for every constant $c\in\R$, as follows by approximating constant potentials by compactly supported smooth cutoffs. As a consequence, when $\mu\in\prob(\R^d)$, the WFR-orthogonal projection of $T_\mu\M$ onto $T_\mu\prob(\R^d)$ is given by
\begin{equation}\label{eq:projection}
\Pi_\mu(\xi,\zeta)
\coloneqq
\left(
\xi,
\zeta-\int_{\R^d}\zeta(x)\,d\mu(x)
\right).
\end{equation}
Thus, the projection leaves the transport component unchanged and only acts on the reaction component by removing its mean.

The proposed construction is analogous to the standard Riemannian-submanifold setting, where the gradient of the restriction of a smooth function to a submanifold is obtained by orthogonally projecting the ambient gradient onto the tangent space of the submanifold.

\begin{definition}\label{def:probability-gradients}
We say that  $F\colon\prob(\R^d)\to\R$ is WFR differentiable at $\mu\in\prob(\R^d)$ if there exists
\begin{equation*}
    \mathcal V=(\xi,\zeta)\in T_\mu\prob(\R^d)
\end{equation*}
such that, for every sequence $(\mu_n)_{n\in\N}\subset\prob(\R^d)$ satisfying $\WFR(\mu_n,\mu)\to0$, 
and every sequence of optimal semicouplings
\begin{equation*}
    \bfgamma^n=(\gamma_0^n,\gamma_1^n)
    \in\operatorname{OptSC}(\mu,\mu_n),
\end{equation*}
one has
\begin{equation}\label{eq:wfr-differentiability-probability}
\begin{split}
    F(\mu_n)& = F(\mu) +\langle \mathcal{V}, \Pi_\mu \barlog_{\bfgamma^n}\rangle_{\mu}
    +o\bigl(\WFR(\mu,\mu_n)\bigr)
\end{split}
\end{equation}
as $n\to+\infty$.
\end{definition}

The previous definition offers an advantage over the classical Wasserstein differentiability notion in that it works in the full space of probabilities without requiring any finiteness conditions on the second moments. Moreover, differently from what happens for flat derivatives, the reaction component of the gradient is canonical, since adding a constant to it would make the gradient not tangent anymore. 

We notice that, since the gradient is assumed to lie in $T_\mu\prob(\R^d)$, by standard properties of the projection
\begin{equation*}
    \langle \mathcal{V}, \Pi_\mu \barlog_{\bfgamma^n}\rangle_{\mu}=\langle \Pi_\mu\mathcal{V}, \barlog_{\bfgamma^n}\rangle_{\mu}=\langle \mathcal{V}, \barlog_{\bfgamma^n}\rangle_{\mu}
\end{equation*}
so that taking the projection also of the logarithm maps to make them tangent to $\prob(\R^d)$ is irrelevant and one can rewrite \eqref{eq:wfr-differentiability-probability} as
\begin{equation*}
        F(\mu_n) = F(\mu) +\langle \mathcal{V}, \barlog_{\bfgamma^n}\rangle_{\mu}
    +o\bigl(\WFR(\mu,\mu_n)\bigr)
\end{equation*}

Similarly to the case of the ambient gradients in $\M$ we have the following result. 

\begin{proposition}[Strength and uniqueness on probability measures]
\label{propo:probability-strongness-uniqueness}
Let $F\colon\prob(\R^d)\to\R$, let $\mu\in\prob(\R^d)$,
and let $\widetilde{\mathcal V}\in T_\mu\prob(\R^d)$
satisfy \eqref{eq:wfr-differentiability-probability}.
Then, for every sequence $(\mu_n)_{n\in\N}\subset\prob(\R^d)$
and every sequence of semicouplings
$\bfgamma^n\in\operatorname{SC}(\mu,\mu_n)$ such that
$\mathcal C(\bfgamma^n)\to0$, one has
\begin{equation}\label{eq:probability-strong-expansion}
    F(\mu_n)
    =F(\mu)+
    \langle
        \widetilde{\mathcal V},
        \barlog_{\bfgamma^n}
    \rangle_\mu
    +o\left(\sqrt{\mathcal C(\bfgamma^n)}\right).
\end{equation}
In particular, if $F$ is differentiable at $\mu$ according
to Definition~\ref{def:probability-gradients}, its gradient
in $T_\mu\prob(\R^d)$ is unique, and we denote it by
\begin{equation}\label{eq:probability-gradient}
    \widetilde\nabla_\WFR F(\mu)
    =
    \bigl(
        \widetilde\nabla^T F(\mu),
        \widetilde\nabla^R F(\mu)
    \bigr).
\end{equation}
\end{proposition}

\begin{proof}
The proof of strength in the ambient space, see Theorem \ref{thm:strongness}, applies unchanged
when the endpoints are restricted in $\prob(\R^d)$.

We now prove uniqueness.
Let $\widetilde{\mathcal V}_1$ and
$\widetilde{\mathcal V}_2$ be two elements of
$T_\mu\prob(\R^d)$ satisfying
\eqref{eq:wfr-differentiability-probability}.
For $\varphi\in C_c^\infty(\R^d)$, set
\begin{equation*}
    T\coloneqq\nabla\varphi,
    \qquad
    \rho\coloneqq
    4\left(\varphi-\int_{\R^d}\varphi\,d\mu\right).
\end{equation*}
For $\varepsilon>0$, consider
\begin{equation*}
    \mu_\varepsilon
    \coloneqq
    (\operatorname{Id}+\varepsilon T)\pf(q_\varepsilon^2\mu),
\end{equation*}
where
\begin{equation*}
        q_\varepsilon
    \coloneqq
    \frac{1+\varepsilon\rho/2}{\sqrt{Z_\varepsilon}},
    \qquad
    Z_\varepsilon
    \coloneqq
    1+\frac{\varepsilon^2}{4}
        \int_{\R^d}\rho^2\,d\mu.
\end{equation*}
Since $\int\rho\,d\mu=0$, one has
$\int q_\varepsilon^2\,d\mu=1$, and hence
$\mu_\varepsilon\in\prob(\R^d)$.
Define
\begin{equation*}
    G_\varepsilon\coloneqq(\operatorname{Id},\operatorname{Id}+\varepsilon T),
    \qquad
    \bfgamma^\varepsilon
    \coloneqq
    \bigl(
        (G_\varepsilon)\pf\mu,
        (G_\varepsilon)\pf(q_\varepsilon^2\mu)
    \bigr)
    \in\operatorname{SC}(\mu,\mu_\varepsilon).
\end{equation*}

The normalized square-root weight can be written as
\begin{equation*}
    q_\varepsilon
    =1+\frac{\varepsilon}{2}\rho_\varepsilon,
    \qquad
    \rho_\varepsilon
    \coloneqq
    \frac{\rho}{\sqrt{Z_\varepsilon}}
    +\frac{2}{\varepsilon}
       \left(\frac{1}{\sqrt{Z_\varepsilon}}-1\right).
\end{equation*}
Since $Z_\varepsilon=1+O(\varepsilon^2)$ and $\mu$ is a probability measure, we have $Z_\varepsilon^{-1/2}-1=O(\varepsilon^2)$ and
\begin{equation*}
    \norm{\rho_\varepsilon-\rho}_{L^2(\mu)}
    \leq
    |Z_\varepsilon^{-1/2}-1|
    \left(\norm{\rho}_{L^2(\mu)}+\frac{2}{\varepsilon}\right)
    =O(\varepsilon).
\end{equation*}
Thus $\rho_\varepsilon\to\rho$ in $L^2(\mu)$.
Lemma~\ref{lemma:graph-linearization}, with
$T_\varepsilon=T$, therefore yields
\begin{equation*}
    \frac{\barlog_{\bfgamma^\varepsilon}}
         {\varepsilon}
    \longrightarrow(T,\rho)
    \quad\text{in }H(\mu),
    \qquad
    \sqrt{\mathcal C(\bfgamma^\varepsilon)}=O(\varepsilon).
\end{equation*}
Applying \eqref{eq:probability-strong-expansion} to each
$\widetilde{\mathcal V}_i$, we obtain
\begin{equation*}
        F(\mu_\varepsilon)-F(\mu)
    =
    \varepsilon
    \langle\widetilde{\mathcal V}_i,(T,\rho)\rangle_\mu
    +o(\varepsilon),
    \qquad i=1,2.
\end{equation*}
Subtracting these identities, dividing by $\varepsilon$,
and passing to the limit gives
\begin{equation*}
        \left\langle
        \widetilde{\mathcal V}_1-\widetilde{\mathcal V}_2,
        \left(
            \nabla\varphi,
            4\left(\varphi-\int_{\R^d}\varphi\,d\mu\right)
        \right)
    \right\rangle_\mu
    =0.
\end{equation*}
The centered potential fields appearing here are the
images of $(\nabla\varphi,4\varphi)$ under $\Pi_\mu$.
They are therefore dense in $T_\mu\prob(\R^d)$.
Since
$\widetilde{\mathcal V}_1-\widetilde{\mathcal V}_2
\in T_\mu\prob(\R^d)$,
we conclude that
$\widetilde{\mathcal V}_1=\widetilde{\mathcal V}_2$.
\end{proof}

\begin{remark}[Metric slope on probability measures]
The identification of the metric slope extends to the
restriction of WFR to $\mathcal P(\R^d)$. Indeed,
if $F\colon\mathcal P(\R^d)\to\R$ is WFR differentiable
at $\mu$, then
\begin{equation*}
    |\partial^-F|_{\mathcal P}(\mu)
    =
    \|\widetilde\nabla_{\WFR}F(\mu)\|_\mu,
\end{equation*}
where
\begin{equation*}
    |\partial^-F|_{\mathcal P}(\mu)
    \coloneqq
    \limsup_{\substack{\nu\to\mu,\ \nu\neq\mu\\
                       \nu\in\mathcal P(\R^d)}}
    \frac{(F(\mu)-F(\nu))_+}{\WFR(\mu,\nu)}.
\end{equation*}
The proof follows as in the ambient case, using
Proposition~\ref{propo:probability-strongness-uniqueness} and the mass-normalized graph
perturbations constructed in its proof for the
lower bound.
\end{remark}

WFR differentiable functions are differentiable according to Definition \ref{def:probability-gradients} when restricted to $\prob(\R^d)$.

\begin{proposition}
Assume $F\colon\M\to\R$ is Wasserstein--Fisher--Rao
differentiable at $\mu\in\prob(\R^d)$. Then its restriction
$F_{\vert\prob}\colon\prob(\R^d)\to\R$ is
Wasserstein--Fisher--Rao differentiable at $\mu$ and
\begin{equation}\label{eq:gradient-identification}
    \widetilde\nabla_\WFR F(\mu)
    =\Pi_\mu\nabla_{\WFR}F(\mu).
\end{equation}
\end{proposition}
\begin{proof}
Write
\begin{equation*}
    \widetilde{\mathcal V}
    \coloneqq\Pi_\mu\nabla_{\WFR}F(\mu)
    \in T_\mu\prob(\R^d).
\end{equation*}
Let $\mu_n\in\prob(\R^d)$ satisfy
$\WFR(\mu_n,\mu)\to0$, and let
$\bfgamma^n=(\gamma_{0}^n,\gamma_{1}^n)
\in\operatorname{OptSC}(\mu,\mu_n)$.
By the WFR differentiability of $F$ in $\M$,
\begin{equation*}
    F(\mu_n)-F(\mu)
    =
    \left\langle
        \nabla_{\WFR}F(\mu),
        \barlog_{\bfgamma^n}
    \right\rangle_\mu
    +
    o\bigl(\WFR(\mu,\mu_n)\bigr).
\end{equation*}
The definition of $\Pi_\mu$ gives
\begin{equation*}
\begin{aligned}
    \left\langle
        \nabla_{\WFR}F(\mu),
        \barlog_{\bfgamma^n}
    \right\rangle_\mu
    -
    \left\langle
        \Pi_\mu\nabla_{\WFR}F(\mu),
        \Pi_\mu\barlog_{\bfgamma^n}
    \right\rangle_\mu=
    \frac14
    \left(\int_{\R^d}\nabla^RF(\mu)\,d\mu\right)
    \left(
        \int_{\R^d}
        \barlog_{\bfgamma^n}^R\,d\mu
    \right).
\end{aligned}
\end{equation*}
Choose a common dominating measure $\lambda^n$ and write
$a_n=d\gamma^n_{0}/d\lambda^n$ and
$b_n=d\gamma^n_{1}/d\lambda^n$.
Since $\mu$ and $\mu_n$ are probability measures,
\eqref{eq:total-mass-formualtion} yields
\begin{align*}
    \WFR(\mu,\mu_n)^2
    &=
    2-2\int_{\R^d\times\R^d}
        \sqrt{a_n(x,y)b_n(x,y)}
        \cos\left(\norm{y-x}\wedge\frac{\pi}{2}\right)
    \,d\lambda^n(x,y)\\
    &=
    -\int_{\R^d\times\R^d}
        \Log_{\bfgamma^n}^R(x,y)
    \,d\gamma^n_{0}(x,y)\\
    &=
    -\int_{\R^d}
        \barlog_{\bfgamma^n}^R
    \,d\mu.
\end{align*}
Consequently, the difference of the two pairings is
$O(\WFR(\mu,\mu_n)^2)=o(\WFR(\mu,\mu_n))$, and hence
\begin{equation*}
    F(\mu_n)-F(\mu)
    =
    \left\langle
        \widetilde{\mathcal V},
        \Pi_\mu\barlog_{\bfgamma^n}
    \right\rangle_\mu
    +
    o\bigl(\WFR(\mu,\mu_n)\bigr),
\end{equation*}
which is precisely the WFR differentiability of
$F_{\vert\prob}$ at $\mu$.
\end{proof}

\begin{remark}[Relation with spherical Hellinger--Kantorovich geometry]
Although WFR geodesics between probability measures need not
preserve mass, the intrinsic geometry on $\prob(\R^d)$ is
described by the spherical Hellinger--Kantorovich distance.
With our normalization, this distance satisfies
\begin{equation*}
    \operatorname{SHK}(\mu,\nu)
    =
    \arccos\left(1-\frac12\WFR^2(\mu,\nu)\right),
    \qquad \mu,\nu\in\prob(\R^d);
\end{equation*}
see \cite{laschos-mielke2019}.
In particular,
\begin{equation*}
    \operatorname{SHK}(\mu,\nu)
    =
    \WFR(\mu,\nu)+O\bigl(\WFR^3(\mu,\nu)\bigr)
    \qquad\text{as }\WFR(\mu,\nu)\to0.
\end{equation*}
Thus replacing $\WFR$ by $\operatorname{SHK}$ in the
little-$o$ remainder of our differentiability definition
on $\prob(\R^d)$ leaves the definition unchanged.
The zero-mean condition on the reaction component and
the projection formula \eqref{eq:gradient-identification}
reflect the unit-mass constraint at the tangent level.
\end{remark}

The main result of this section is the following one.

\begin{theorem}[Wasserstein--flat decomposition of the WFR gradient]
\label{thm:wfr-wasserstein-flat-decomposition}
Let $\mu\in\prob_2(\R^d)$ and let
$F\colon\prob(\R^d)\to\R$ be WFR differentiable at $\mu$, with
\begin{equation*}
    \widetilde\nabla_{\mathrm{WFR}}F(\mu)
=
\bigl(
\widetilde\nabla^T F(\mu),
\widetilde\nabla^R F(\mu)
\bigr).
\end{equation*}
Then $F|_{\prob_2(\R^d)}$ is Wasserstein differentiable at $\mu$ and
\begin{equation*}
    \widetilde\nabla^T F(\mu)=\nabla_{W_2}F(\mu).
\end{equation*}
If, in addition, $F$ is flat differentiable at $\mu$ (cf. Definition \ref{def:flat-differentiability}), then
\begin{equation*}
    \widetilde\nabla_{\mathrm{WFR}}F(\mu)
=
\left(
\nabla_{W_2}F(\mu),
4\left[
\frac{\delta F}{\delta\mu}(\mu,\cdot)
-
\int_{\R^d}
\frac{\delta F}{\delta\mu}(\mu,y)\,d\mu(y)
\right]
\right).
\end{equation*}
\end{theorem}

The proof is split into the next two subsections. We first identify the
transport component with the Wasserstein gradient. We then study
reaction perturbations and show that, whenever a flat derivative exists,
the reaction component is its canonical zero-mean representative.

\subsection{The Wasserstein component}

If we restrict to $\probtwo$, we recover Wasserstein
differentiability. In this sense, we can think of
Wasserstein--Fisher--Rao calculus as an extension of
Wasserstein calculus, and of the Wasserstein gradient
as the partial derivative corresponding to spatial
transport.

For completeness we recall the definition of a Wasserstein differentiable function.

\begin{definition}[Wasserstein differentiability]
Let $F\colon\probtwo\to\R$. We say that $F$ is Wasserstein differentiable at $\mu$ if there exists
\begin{equation*}
    \xi \in T_\mu\probtwo \coloneqq \overline{\left\{\nabla\psi\colon\psi\in C_c^\infty(\R^d)\right\}}^{L^2(\mu)}
\end{equation*}
such that, for every sequence $(\mu_n)_{n\in\N}\subset\probtwo$ satisfying $W_2(\mu_n,\mu)\to0$, and every sequence of optimal couplings $\gamma^n$
one has
\begin{equation}\label{eq:w-differentiability}
    F(\mu_n)=F(\mu)+\int_{\R^d\times\R^d} \langle \xi(x),y-x\rangle \,d\gamma_n(x,y)+o\bigl(W_2(\mu,\mu_n)\bigr)
\end{equation}
as $n\to+\infty$. If such a $\xi$ exists, it is unique, we denote it by $\nabla_{W_2}F(\mu)$ and call it \emph{Wasserstein gradient} of $F$ at $\mu$.
\end{definition}

The main result is the following, which follows from Proposition \ref{propo:probability-strongness-uniqueness} by using specific semicouplings built from optimal couplings.

\begin{theorem}
Let $\mu\in\mathcal P_2(\R^d)$ and let
$F\colon\prob(\R^d)\to\R$ be WFR differentiable at $\mu$. Then the restriction
of $F$ to $\mathcal P_2(\R^d)$ is Wasserstein differentiable at $\mu$ and 
\begin{equation*}
    \nabla_{W_2}F(\mu)=\widetilde\nabla^TF(\mu).
\end{equation*}
\end{theorem}

\begin{proof}
Let $(\mu_n)_{n\in\N}\subset\prob_2(\R^d)$ satisfy
\begin{equation*}
    d_n\coloneqq W_2(\mu,\mu_n)\longrightarrow0,
\end{equation*}
and let $\gamma_n$ be a sequence of optimal couplings between $\mu$ and $\mu_n$. That is
\begin{equation*}
    d_n^2
=
\int_{\R^d\times\R^d}\norm{x-y}^2\,d\gamma_n(x,y).
\end{equation*}

Consider the pure-transport semicoupling $\bfgamma^n=(\gamma_n,\gamma_n)$.
Taking $\lambda_n=\gamma_n$, the corresponding densities satisfy
$a_n=b_n=1$ and its semicoupling cost satisfies
\begin{align*}
\mathcal C(\bfgamma^n)
=
\int_{\R^d\times\R^d}
2\left(
1-\cos(\norm{y-x}\wedge\frac{\pi}{2})
\right)
\,d\gamma_n(x,y)
\leq
\int_{\R^d\times\R^d}|x-y|^2\,d\gamma_n(x,y)
=
d_n^2.
\end{align*}
In particular, the cost of $\bfgamma^n$ converges to zero. Therefore,
by Proposition~\ref{propo:probability-strongness-uniqueness},
\begin{align}\label{eq:wasserstein-wfr-expansion}
F(\mu_n)-F(\mu)
&=
\int_{\R^d\times\R^d}
\widetilde\nabla^TF(\mu)(x)\cdot
\Log_{\bfgamma^n}^T(x,y)
\,d\gamma_n(x,y)
\nonumber\
\\ &+
\frac14
\int_{\R^d\times\R^d}
\widetilde\nabla^RF(\mu)(x)
\Log_{\bfgamma^n}^R(x,y)
\,d\gamma_n(x,y)
+o(d_n).
\end{align}

We first show that the reaction contribution is negligible. Set
$\delta=\norm{x-y}$ and define
\begin{equation*}
    h_R(\delta)
\coloneqq
\begin{cases}
0,&\delta=0,\\[1mm]
\displaystyle
\frac{2(\cos(\delta)-1)}\delta,
&0<\delta<\frac{\pi}{2},\\[1mm]
\displaystyle-\frac2\delta,
&\delta\geq\frac{\pi}{2}.
\end{cases}
\end{equation*}
Then $h_R$ is bounded, $h_R(\delta)\to0$ as $\delta\to0$, and
\begin{equation*}
    \Log_{\bfgamma^n}^R(x,y)
=
|x-y|h_R(|x-y|).
\end{equation*}
Consequently, by the Cauchy--Schwarz inequality, by setting
\begin{equation*}
    H^F_R(x,y)\coloneqq \norm{\widetilde\nabla^RF(\mu)(x)} ^2
\lvert h_R(\norm{x-y})\rvert ^2
\end{equation*}
we have
\begin{align*}
\left|
\int_{\R^d\times\R^d}
\widetilde\nabla^RF(\mu)(x)
\Log_{\bfgamma^n}^R(x,y)
\,d\gamma^nn(x,y)
\right|
\leq
d_n
\left(
\int_{\R^d\times\R^d}
H_R^F(x,y)
\,d\gamma^n(x,y)
\right)^{1/2}.
\end{align*}
The second factor converges to zero by Corollary \ref{cor:vanishing-weights}.
Therefore
\begin{equation}\label{eq:wasserstein-reaction-negligible}
\int_{\R^d\times\R^d}
\widetilde\nabla^RF(\mu)(x)
\Log_{\bfgamma^n}^R(x,y)
\,d\gamma_n(x,y)
=
o(d_n).
\end{equation}

We next consider the transport contribution. We have
\begin{equation*}
    \Log_{\bfgamma^n}^T(x,y)
=
\kappa(\norm{x-y})(y-x),
\end{equation*}
where
\begin{equation*}
    \kappa(\delta)
\coloneqq
\begin{cases}
1&\delta=0\\[1mm]
\displaystyle\frac{\sin \delta}{\delta},
&0< \delta \leq\frac{\pi}{2}\\[1mm]
0,
&\delta>\frac{\pi}{2}.
\end{cases}
\end{equation*}
Hence
\begin{equation*}
    \Log_{\bfgamma^n}^T(x,y)-(y-x)
=
h_T(\norm{x-y})(y-x),
\end{equation*}
where $h_T(\delta)\coloneqq\kappa(\delta)-1.$
The function $h_T$ is bounded and satisfies
\begin{equation*}
    h_T(\delta)\longrightarrow0
\qquad\text{as }\delta\to0.
\end{equation*}
Therefore, for
\begin{equation*}
    H^F_T(x,y)\coloneqq \norm{\widetilde\nabla^TF(\mu)(x)} ^2
\lvert h_T(\norm{x-y})\rvert ^2
\end{equation*}
it holds
\begin{align*}
\left|
\int_{\R^d\times\R^d}
\widetilde\nabla^TF(\mu)(x)\cdot
\left(
\Log_{\bfgamma^n}^T(x,y)-(y-x)
\right)
\,d\gamma_n(x,y)
\right|
\leq
d_n
\Biggl(
\int_{\R^d\times\R^d}
H^F_T(x,y)
\,d\gamma_n(x,y)
\Biggr)^{1/2}.
\end{align*}
Again, the second factor converges to zero by Corollary \ref{cor:vanishing-weights}.
It follows that
\begin{equation}\label{eq:wasserstein-transport-expansion}
\int_{\R^d\times\R^d}
\widetilde\nabla^TF(\mu)(x)\cdot
\Log_{\bfgamma^n}^T(x,y)
\,d\gamma_n(x,y)
=
\int_{\R^d\times\R^d}
\widetilde\nabla^TF(\mu)(x)\cdot(y-x)
\,d\gamma_n(x,y)
+
o(d_n).
\end{equation}

Combining
\eqref{eq:wasserstein-wfr-expansion},
\eqref{eq:wasserstein-reaction-negligible}, and
\eqref{eq:wasserstein-transport-expansion}, we obtain
\begin{equation*}
    F(\mu_n)-F(\mu)
=
\int_{\R^d\times\R^d}
\widetilde\nabla^TF(\mu)(x)\cdot(y-x)
\,d\gamma_n(x,y)
+
o\bigl(W_2(\mu,\mu_n)\bigr).
\end{equation*}
Finally, since
\begin{equation*}
    \bigl(
\widetilde\nabla^TF(\mu),
\widetilde\nabla^RF(\mu)
\bigr)
\in T_\mu\M,
\end{equation*}
there exists a sequence $(\varphi_k)_{k\in\N}\subset C_c^\infty(\R^d)$
such that
\begin{equation*}
    (\nabla\varphi_k,4\varphi_k)
\longrightarrow
\bigl(
\widetilde\nabla^TF(\mu),
\widetilde\nabla^RF(\mu)
\bigr)
\end{equation*}
in $H(\mu)$. In particular,
\begin{equation*}
    \nabla\varphi_k
\longrightarrow
\widetilde\nabla^TF(\mu)
\qquad\text{in }L^2(\mu;\R^d).
\end{equation*}
Thus $\widetilde\nabla^TF(\mu)\in T_\mu\mathcal P_2(\R^d)$ and the proof is completed.
\end{proof}

We have the following immediate corollary.
\begin{corollary}
Let $\mu\in\prob_2(\R^d)$ and $\hat\mu\in\M$
satisfy the equivalent conditions of
Theorem~\ref{thm:squared-wfr-differentiability}.
Then the restriction of $F_{\hat\mu}$ to $\prob_2(\R^d)$
is Wasserstein differentiable at $\mu$ and
\begin{equation*}
    \nabla_{W_2}
    \bigl(F_{\hat\mu}|_{\prob_2(\R^d)}\bigr)(\mu)
    =-\mathcal L^T,
\end{equation*}
where $\mathcal L$ is the unique regular barycentric logarithm
from $\mu$ to $\hat\mu$.
Equivalently, with the maps $T,q$ of Theorem \ref{thm:squared-wfr-differentiability},
\begin{equation*}
    \nabla_{W_2}
    \bigl(F_{\hat\mu}|_{\prob_2(\R^d)}\bigr)(\mu)(x)
    =
    q(x)\frac{\sin\bigl(\norm{T(x)-x}\bigr)}{\norm{T(x)-x}}(x-T(x))
    \qquad \mu\text{-a.e.},
\end{equation*}
where
$\sin\theta/\theta=1$ at $\theta=0$.
\end{corollary}

\subsection{The reaction component and flat directions}\label{subsec:flat-derivatives}

We now compare WFR differentiability with flat differentiability. We first recall the local notion of flat differentiability that we shall use.

\begin{definition}[Flat differentiability]\label{def:flat-differentiability}
Let $F\colon\prob(\R^d)\to\R$ and $\mu\in\prob(\R^d)$. We say that $F$ is flat differentiable at $\mu$ if there exists a function
\begin{equation*}
    \frac{\delta F}{\delta\mu}(\mu,\cdot)\in C_b(\R^d)
\end{equation*}
such that, for every $\nu\in\mathcal P(\R^d)$,
\begin{equation}\label{eq:flat-differentiability}
    \lim_{\varepsilon\downarrow0}
    \frac{
        F\bigl((1-\varepsilon)\mu+\varepsilon\nu\bigr)-F(\mu)
    }{\varepsilon}
    =
    \int_{\R^d}
        \frac{\delta F}{\delta\mu}(\mu,x)\,
        d(\nu-\mu)(x).
\end{equation}
\end{definition}

\begin{remark}
The flat derivative is determined only up to an additive constant.
\end{remark}

Every first-order Wasserstein perturbation is also a first-order WFR
perturbation, which allowed us to derive Wasserstein differentiability from WFR differentiability. Flat perturbations, on the other hand, form a larger
class: first-order flat perturbations may create mass on sets that are
negligible for the reference measure, and such creation need not occur
at order $\varepsilon$ in the WFR distance.

The following result identifies the flat directions that are directly
controlled by WFR differentiability.

\begin{proposition}[Absolutely continuous flat directions]
Let $F\colon\prob(\R^d)\to\R$ be WFR differentiable at
$\mu\in\prob(\R^d)$. Let $\nu\in\prob(\R^d)$, and consider its Lebesgue decomposition
with respect to $\mu$ to be $\nu=h\mu+\nu^\perp$.
For $\varepsilon\in[0,1]$, set
\begin{equation*}
    \mu^\varepsilon\coloneqq(1-\varepsilon)\mu+\varepsilon\nu
\end{equation*}
and define the diagonal semicoupling $\bfgamma^\varepsilon$ as
\begin{equation*}
    \gamma_0^\varepsilon
    \coloneqq (\operatorname{Id},\operatorname{Id})\pf\mu,
    \qquad
    \gamma_1^\varepsilon
    \coloneqq (\operatorname{Id},\operatorname{Id})\pf\mu^\varepsilon,
\end{equation*}
then
\begin{equation}\label{eq:diagonal-mixture-cost}
    \mathcal C(\bfgamma^\varepsilon)
    =
    \int_{\R^d}
        \left(
            \sqrt{1+\varepsilon(h-1)}-1
        \right)^2
        d\mu
    +
    \varepsilon\nu^\perp(\R^d).
\end{equation}
In particular, if
\begin{equation*}
    \nu=h\mu,
    \qquad
    h-1\in L^2(\mu),
\end{equation*}
then
\begin{equation}\label{eq:flat-directional-derivative}
    \lim_{\varepsilon\downarrow0}
    \frac{F(\mu^\varepsilon)-F(\mu)}{\varepsilon}=
    \frac14
    \int_{\R^d}\widetilde\nabla^RF(\mu)(x)\,d(\nu-\mu)(x).
\end{equation}
Thus $F$ admits a flat directional derivative along every
absolutely continuous direction whose density belongs to $1+L^2(\mu)$.
\end{proposition}

\begin{proof}
The measure
\begin{equation*}
    \lambda\coloneqq(\operatorname{Id},\operatorname{Id})_\#(\mu+\nu^\perp)
\end{equation*}
dominates both components of $\bfgamma^\varepsilon$. On the part corresponding to $\mu$, the densities are
$a=1$ and $b=1+\varepsilon(h-1)$, whereas on a
$\mu$-negligible set carrying $\nu^\perp$ they are
$a=0$ and $b=\varepsilon$. Since all pairs lie on the diagonal, the first part contributes
\begin{equation*}
        \int_{\R^d}
        \left(\sqrt{1+\varepsilon(h(x)-1)}-1
        \right)^2
        d\mu(x)
\end{equation*}
to the semicoupling cost, while the singular part is pure target
creation and contributes $\varepsilon\,\nu^\perp(\R^d)$. This proves
\eqref{eq:diagonal-mixture-cost}.

Suppose that $\nu=h\mu$ and $h-1\in L^2(\mu)$.
Set $T_\varepsilon=0$ and
\begin{equation*}
    \rho_\varepsilon
    \coloneqq
    \frac{2(h-1)}{\sqrt{1+\varepsilon(h-1)}+1}.
\end{equation*}
Then
\begin{equation*}
    \sqrt{1+\varepsilon(h-1)}
    =1+\frac{\varepsilon}{2}\rho_\varepsilon.
\end{equation*}
Moreover, $\rho_\varepsilon\to h-1$ pointwise and
$|\rho_\varepsilon|\leq2|h-1|$, so the convergence also
holds in $L^2(\mu)$.
Lemma~\ref{lemma:graph-linearization} gives
\begin{equation*}
    \frac{\barlog_{\bfgamma^\varepsilon}}
         {\varepsilon}
    \longrightarrow(0,h-1)
    \quad\text{in }H(\mu),
    \qquad
    \frac{\mathcal C(\bfgamma^\varepsilon)}{\varepsilon^2}
    \longrightarrow\frac14\int_{\R^d}(h-1)^2\,d\mu.
\end{equation*}
Since $\sqrt{\mathcal C(\bfgamma^\varepsilon)}=O(\varepsilon)$,
Proposition~\ref{propo:probability-strongness-uniqueness} gives
\begin{equation*}
    F(\mu^\varepsilon)-F(\mu)
    =
    \left\langle
        \widetilde\nabla_\WFR F(\mu),
        \barlog_{\bfgamma^\varepsilon}
    \right\rangle_\mu
    +o(\varepsilon).
\end{equation*}
Dividing by $\varepsilon$ and passing to the limit yields
\begin{equation*}
\begin{split}
    \lim_{\varepsilon\downarrow0}
    \frac{F(\mu^\varepsilon)-F(\mu)}{\varepsilon}
    &=
    \left\langle
        \widetilde\nabla_\WFR F(\mu),(0,h-1)
    \right\rangle_\mu\\
    &=
    \frac14\int_{\R^d}
        \widetilde\nabla^R F(\mu)(x)(h(x)-1)\,d\mu(x)\\
    &=
    \frac14\int_{\R^d}
        \widetilde\nabla^R F(\mu)(x)\,d(\nu-\mu)(x).\qedhere
\end{split}
\end{equation*}
\end{proof}

The preceding result also clarifies the relation between the reaction component of the WFR gradient and the flat derivative. Suppose that $F$ is additionally flat differentiable at $\mu$. Comparing \eqref{eq:flat-directional-derivative} with \eqref{eq:flat-differentiability}, we obtain
\begin{equation*}
    \int_{\R^d}
\left(
\frac{\delta F}{\delta\mu}(\mu,x)
-\frac14\widetilde\nabla^R F(\mu)(x)
\right)(h(x)-1)\,d\mu(x)=0
\end{equation*}
for every $h\ge 0$ with $\int h\,d\mu=1$ and $h-1\in L^2(\mu)$.
In particular, taking $h=1+\varepsilon\psi$ for any bounded zero-mean function $\psi$ and sufficiently small $\varepsilon>0$, we deduce that the expression in parentheses is orthogonal to every bounded zero-mean function. By density in the zero-mean subspace of $L^2(\mu)$, it must be constant $\mu$-a.e. Thus
\begin{equation}\label{eq:flat-reaction-identification}
\frac{\delta F}{\delta\mu}(\mu,x)
=\frac14\widetilde\nabla^R F(\mu)(x)+c_\mu,
\qquad \mu\text{-a.e.},
\end{equation}
for some constant $c_\mu$.
It is important to note that the ambiguity in
\eqref{eq:flat-reaction-identification} belongs exclusively to the flat
derivative. Indeed, flat variations inside $\prob(\R^d)$ have zero total mass, and therefore
\begin{equation*}
        \int_{\R^d}c\,d(\nu-\mu)=0
\end{equation*}
for every constant $c$. Consequently, the flat derivative cannot
distinguish between two representatives differing by an additive
constant.

By contrast, the reaction component $\widetilde\nabla^RF(\mu)$ of the WFR gradient is
uniquely determined by the fact that $\widetilde\nabla_{\WFR}F(\mu)\in T_\mu \prob(\R^d)$. In this sense, $\frac14\widetilde\nabla^RF(\mu)$
provides a canonical representative, in $L^2(\mu)$, of the flat
derivative along absolutely continuous directions.

The preceding proposition concerns flat perturbations that are also
first-order WFR perturbations. Not every flat perturbation has this
property. Although singular mass can be generated by
transport at order $\varepsilon$, this becomes impossible when new
mass is placed beyond the transport threshold. The following example
shows that, in this case, the intrinsic WFR scale is
$\sqrt{\varepsilon}$, and consequently WFR differentiability does
not imply flat differentiability.

\begin{example}
Fix $1<p<2$ and $\hat\mu\in\prob(\R^d)$, and consider $F\colon\prob(\R^d)\to\R$ defined by
\begin{equation*}
        F(\mu)\coloneqq\WFR(\mu,\hat\mu)^p.
\end{equation*}
Then $F$ is WFR differentiable at $\hat\mu$, with
\begin{equation*}
    \widetilde\nabla_{\WFR}F(\hat\mu)=(0,0)
\end{equation*}
given that
\begin{equation*}
    F(\mu)-F(\hat\mu)
    =
    \WFR(\mu,\hat\mu)^p
    =
    o\bigl(\WFR(\mu,\hat\mu)\bigr).
\end{equation*}

We now fix $d=1$, set
\begin{equation*}
    \widehat\mu=\delta_0,
    \qquad
    \nu=\delta_L,
    \qquad
    |L|>\frac{\pi}{2},
\end{equation*}
and consider the segment
\begin{equation*}
    \mu^\varepsilon=(1-\varepsilon)\delta_0+\varepsilon\delta_L.
\end{equation*}
Since $0$ and $L$ lie beyond the transport threshold, mass at
$0$ cannot be effectively transported to $L$. It follows that
\begin{equation}\label{eq:dirac-mixture-wfr}
    \WFR^2(\delta_0,\mu^\varepsilon)
    =
    \left(1-\sqrt{1-\varepsilon}\right)^2+\varepsilon
    =
    2\left(1-\sqrt{1-\varepsilon}\right).
\end{equation}
Consequently,
\begin{equation*}
    \WFR^2(\delta_0,\mu^\varepsilon)
    =
    \varepsilon+O(\varepsilon^2),
    \qquad
    \WFR(\delta_0,\mu^\varepsilon)
    \sim\sqrt{\varepsilon}.
\end{equation*}
Therefore
\begin{equation*}
    \frac{F(\mu^\varepsilon)-F(\hat\mu)}{\varepsilon}
    =
    \frac{
        \left[
            2\left(1-\sqrt{1-\varepsilon}\right)
        \right]^{p/2}
    }{\varepsilon}
    \sim
    \varepsilon^{p/2-1}
    \longrightarrow+\infty.
\end{equation*}
Thus $F$ does not possess a finite flat directional derivative at
$\widehat\mu$ in the direction $\delta_L-\delta_0$, and in
particular it is not flat differentiable at $\hat\mu$.
\end{example}

\bibliographystyle{plain}
\bibliography{bibliography}

\appendix

\section{Expansions of costs and logarithms}

\subsection{Calculations along cone geodesics}
\label{app:cone-geodesic-calculations}

We prove Lemmas~\ref{lemma:cone-geodesic-velocities}
and~\ref{lemma:cone-observable-estimates}, using the notation of
Section~\ref{subsec:wfr-geodesics}.

\begin{proof}[Proof of Lemma~\ref{lemma:cone-geodesic-velocities}]
In the first case $r,s>0$, direct differentiation of
$\bfm_t$ yields the formula for $\dot{\bfm}_t$, and the
chain rule proves the formula for $\dot{\bfr}_t$ whenever
$\bfr_t>0$. If $\norm{y-x}=0$, then $\bfx_t\equiv x$, and the formula for
$\dot{\bfx}_t$ is immediate.
If $0<\norm{y-x}<\pi$, then $v\neq0$ and
\begin{equation*}
    \bfx_t
    =
    x+
    \operatorname{Arg}\left(1+tu+\mathrm{i}t\norm{v}\right)
    \frac{v}{\norm{v}},
\end{equation*}
where the complex number inside the argument never vanishes.

Consequently, direct differentiation gives
\begin{equation*}
    \dot{\bfx}_t=\frac{v}{(1+tu)^2+t^2\norm{v}^2}.
\end{equation*}
If $\norm{y-x}\geq\pi$, then $v=0$ and $\bfx_t$ is constant
on each side of the vertex time $\tau=r/(r+s)$.
Thus the same formula holds whenever $\bfr_t>0$.
Since $s>0$, the velocity formulas also extend continuously
to $t=1$ and give the corresponding left derivatives there.

When $s>0$, substituting the velocity formulas and the expression for
$\bfr_t$ gives
\begin{equation*}
\begin{aligned}
    \bfr_t^2\norm{\dot{\bfx}_t}^2+|\dot{\bfr}_t|^2
    =r^2\bigl(u^2+\norm{v}^2\bigr)
    =r^2+s^2-2rs
    K_y(x)
    =\dc^2(\bfy_0,\bfy_1).
\end{aligned}
\end{equation*}
If $\bfy_1=\vertex$, the radial formula gives
$\dot{\bfx}_t=0$ and $\dot{\bfr}_t=-r$, so the same identity follows
from $\dc^2(\bfy_0,\vertex)=r^2$.
\end{proof}

\begin{proof}[Proof of Lemma~\ref{lemma:cone-observable-estimates}]
Suppose first that $r,s>0$ and $\norm{y-x}<\pi$.
From the formula for $\bfm_t$ we have $\bfm_t>0$ in $[0,1]$.
Hence the geodesic never reaches the vertex. Since $\Phi_\varphi(\bfy_t)=\bfm_t\varphi(\bfx_t)$
we obtain from the chain rule
\begin{equation*}
    \frac{d}{dt}\Phi_\varphi(\bfy_t)
    =
    \dot{\bfm}_t\varphi(\bfx_t)
    +
    rs\,\nabla\varphi(\bfx_t)
    \cdot\mathbf{sin}(y-x).
\end{equation*}
At $t=0$, using the expression for
$\Log_\cone(\bfy_0,\bfy_1)$ in Section~\ref{subsec:wfr-geodesics}, this is
precisely the claimed initial derivative.

Moreover, $\dot{\bfm}_t$ is affine in $t$ and
$|rs\cos(\norm{y-x})|\leq rs$, so
\begin{equation*}
    |\dot{\bfm}_t|
    \leq\max\{|\dot{\bfm}_0|,|\dot{\bfm}_1|\}
    \leq2\bigl(\max\{r^2,s^2\}+rs\bigr)
    \leq3(r^2+s^2),
    \qquad
    rs\norm{\mathbf{sin}(y-x)}
    \leq rs\leq\frac12(r^2+s^2).
\end{equation*}
Thus the formula for the derivative yields the claimed derivative bound. Integrating over $[0,t]$ gives the difference-quotient bound.

For the quadratic estimate, assume first that $\varphi\in C_b^2(\R^d)$.
Differentiating the velocity formulas of
Lemma~\ref{lemma:cone-geodesic-velocities} gives, for $t\in(0,1)$,
\begin{equation*}
    \ddot{\bfr}_t
    =
    \bfr_t\norm{\dot{\bfx}_t}^2,
    \qquad
    \ddot{\bfx}_t
    =
    -2\frac{\dot{\bfr}_t}{\bfr_t}\dot{\bfx}_t.
\end{equation*}
Differentiating $\Phi_\varphi(\bfy_t)=\bfr_t^2\varphi(\bfx_t)$ twice
and using \eqref{eq:constant-speed-norm} yields
\begin{equation*}
\begin{aligned}
    \frac{d^2}{dt^2}\Phi_\varphi(\bfy_t)
    =
    \bfr_t^2
    \dot{\bfx}_t^\top
    \bigl(D^2\varphi(\bfx_t)+2\varphi(\bfx_t)I\bigr)
    \dot{\bfx}_t
    +2\bfr_t\dot{\bfr}_t
        \nabla\varphi(\bfx_t)\cdot\dot{\bfx}_t
    +2\varphi(\bfx_t)|\dot{\bfr}_t|^2.
\end{aligned}
\end{equation*}
Since
\begin{equation*}
    2\bfr_t\norm{\dot{\bfx}_t}|\dot{\bfr}_t|
    \leq
    \bfr_t^2\norm{\dot{\bfx}_t}^2
    +
    |\dot{\bfr}_t|^2,
\end{equation*}
the constant-speed identity gives
\begin{equation*}
    \left|
        \frac{d^2}{dt^2}\Phi_\varphi(\bfy_t)
    \right|
    \leq
    \left(
        \norm{D^2\varphi}_{\infty,\mathrm{op}}
        +2\norm{\varphi}_\infty
        +\norm{\nabla\varphi}_\infty
    \right)\dc(\bfy_0,\bfy_1)^2.
\end{equation*}
Taylor's formula with integral remainder therefore yields the claimed quadratic estimate.

If $r,s>0$ and $\norm{y-x}\geq\pi$, set
$\tau\coloneqq r/(r+s)$. Then
\begin{equation*}
    \Phi_\varphi(\bfy_t)
    =
    \begin{cases}
        (r+s)^2(t-\tau)^2\varphi(x), & 0\leq t\leq\tau,\\
        (r+s)^2(t-\tau)^2\varphi(y), & \tau<t\leq1.
    \end{cases}
\end{equation*}
Its initial right derivative is $-2r(r+s)\varphi(x)$,
which agrees with the claimed formula since
$\Log_\cone(\bfy_0,\bfy_1)=(0,-r-s)$.
The map is continuously differentiable across $t=\tau$,
with derivative zero there, and
\begin{equation*}
    \left|\frac{d}{dt}\Phi_\varphi(\bfy_t)\right|
    \leq
    2(r+s)\max\{r,s\}\norm{\varphi}_\infty
    \leq
    3(r^2+s^2)\norm{\varphi}_\infty.
\end{equation*}
This proves the derivative and difference-quotient bounds.
Moreover, its first derivative is absolutely continuous and
\begin{equation*}
    \left|\frac{d^2}{dt^2}\Phi_\varphi(\bfy_t)\right|
    \leq
    2\norm{\varphi}_\infty(r+s)^2
    =
    2\norm{\varphi}_\infty\dc(\bfy_0,\bfy_1)^2
\end{equation*}
for almost every $t\in(0,1)$. Taylor's formula with
integral remainder again yields the claimed quadratic estimate.

If one endpoint is the vertex, the geodesic is radial and
\begin{equation*}
    \Phi_\varphi(\bfy_t([x,r],\vertex))
    =(1-t)^2r^2\varphi(x),
    \qquad
    \Phi_\varphi(\bfy_t(\vertex,[y,s]))
    =t^2s^2\varphi(y).
\end{equation*}
In particular, the initial right derivatives are respectively
$-2r^2\varphi(x)$ and $0$. The former agrees with the claimed initial derivative since $\Log_\cone(\bfy_0,\vertex)=(0,-r)$, while the derivative, difference-quotient, and quadratic bounds follow directly. If both endpoints are the vertex, all claims are immediate.

Finally, standard mollification extends the quadratic estimate to
$\varphi\in C_b^{1,1}(\R^d)$ with the stated constant;
see \cite[Appendix~C.5 and Section~5.3.1]{evans2010partial}.
\end{proof}

\subsection{Graph semicouplings}

\begin{lemma}[Linearization of graph semicouplings]
\label{lemma:graph-linearization}
Let $\mu\in\M$ and suppose that
\begin{equation*}
    (T_\varepsilon,\rho_\varepsilon)
    \longrightarrow (T,\rho)
    \qquad\text{in }H(\mu)
\end{equation*}
as $\varepsilon\to0$. Set
\begin{equation*}
    X_\varepsilon\coloneqq\operatorname{Id}+\varepsilon T_\varepsilon,
    \qquad
    q_\varepsilon\coloneqq1+\frac{\varepsilon}{2}\rho_\varepsilon,
    \qquad
    G_\varepsilon\coloneqq(\operatorname{Id},X_\varepsilon),
\end{equation*}
and define
\begin{equation*}
    \mu^\varepsilon
    \coloneqq(X_\varepsilon)_\#(q_\varepsilon^2\mu),
    \qquad
    \bfgamma^\varepsilon
    \coloneqq
    \bigl((G_\varepsilon)_\#\mu,\,
    (G_\varepsilon)_\#(q_\varepsilon^2\mu)\bigr).
\end{equation*}
Then $\bfgamma^\varepsilon\in\operatorname{SC}(\mu,\mu^\varepsilon)$
and
\begin{equation}\label{eq:graph-log-linearization}
    \frac{\barlog_{\bfgamma^\varepsilon}}
         {\varepsilon}
    \longrightarrow(T,\rho)
    \qquad\text{in }H(\mu).
\end{equation}
Moreover,
\begin{equation}\label{eq:graph-cost-linearization}
    \frac{\mathcal C(\bfgamma^\varepsilon)}{\varepsilon^2}
    \longrightarrow
    \norm{(T,\rho)}_\mu^2.
\end{equation}
The same conclusions hold for one-sided limits.
\end{lemma}

\begin{proof}
It is straightforward to check that
$\bfgamma^\varepsilon\in\operatorname{SC}(\mu,\mu^\varepsilon)$.
Let $\lambda^\varepsilon\coloneqq(G_\varepsilon)_\#\mu$.
Then $\lambda^\varepsilon$ dominates both
$\gamma_0^\varepsilon$ and $\gamma_1^\varepsilon$.
Writing
$a^\varepsilon=d\gamma_0^\varepsilon/d\lambda^\varepsilon$
and
$b^\varepsilon=d\gamma_1^\varepsilon/d\lambda^\varepsilon$,
the graph structure gives the pulled-back densities
\begin{equation*}
    \tilde a^\varepsilon(x)
    \coloneqq a^\varepsilon(G_\varepsilon(x))=1,
    \qquad
    \tilde b^\varepsilon(x)
    \coloneqq b^\varepsilon(G_\varepsilon(x))=q_\varepsilon^2(x),
    \qquad\mu\text{-a.e.}
\end{equation*}
Thus
\begin{equation*}
\begin{split}
    \mathcal C(\bfgamma^\varepsilon)
    =
    \int_{\R^d\times\R^d}
    c(x,y,a^\varepsilon(x,y),b^\varepsilon(x,y))
    \,d\lambda^\varepsilon(x,y)=
    \int_{\R^d}c_\varepsilon(x)\,d\mu(x),
\end{split}
\end{equation*}
where
\begin{equation*}
    c_\varepsilon(x)
    \coloneqq
    \tilde a^\varepsilon(x)+\tilde b^\varepsilon(x)
    -2\sqrt{\tilde a^\varepsilon(x)\tilde b^\varepsilon(x)}
    \cos\left(\norm{\varepsilon T_\varepsilon(x)}
              \wedge\frac{\pi}{2}\right).
\end{equation*}
Similarly, the graph structure gives
\begin{equation*}
    \mathcal L_\varepsilon(x)
    \coloneqq
    \Log_{\bfgamma^\varepsilon}
        (x,X_\varepsilon(x))
    =
    \barlog_{\bfgamma^\varepsilon}(x).
\end{equation*}
We can rewrite the pointwise cost as
\begin{equation*}
    c_\varepsilon(x)
    =
    (|q_\varepsilon(x)|-1)^2
    +2|q_\varepsilon(x)|
    \left(
        1-\cos\left(\norm{\varepsilon T_\varepsilon(x)}
                    \wedge\frac{\pi}{2}\right)
    \right).
\end{equation*}
Using
\begin{equation*}
\begin{split}
    \bigl||q_\varepsilon(x)|-1\bigr|
    \leq\frac{|\varepsilon|}{2}|\rho_\varepsilon(x)|,\qquad
    0\leq
    2\left(
        1-\cos\left(\norm{\varepsilon T_\varepsilon(x)}
                    \wedge\frac{\pi}{2}\right)
    \right)
    \leq\min\{\varepsilon^2\norm{T_\varepsilon(x)}^2,2\},
\end{split}
\end{equation*}
together with Young's inequality and the pointwise estimate
in the proof of Lemma~\ref{lemma:energy-bound}, we obtain
\begin{equation}\label{eq:graph-linearization-bound}
    \frac{\norm{\mathcal L_\varepsilon}_{\WFR}^2}{\varepsilon^2}
    \leq
    \frac{c_\varepsilon}{\varepsilon^2}
    \leq
    \frac32\norm{T_\varepsilon}^2
    +\frac12|\rho_\varepsilon|^2.
\end{equation}

Fix any sequence $\varepsilon_n\to0$ with $\varepsilon_n\neq0$.
By strong convergence in $H(\mu)$ and
\cite[Theorem~4.9]{brezis2011}, applied componentwise,
we may extract a subsequence, not relabeled, such that
\begin{equation*}
    (T_{\varepsilon_n},\rho_{\varepsilon_n})
    \longrightarrow(T,\rho)
    \qquad\mu\text{-a.e.},
\end{equation*}
and there exists $\Psi\in L^2(\mu)$ satisfying
\begin{equation*}
    \norm{T_{\varepsilon_n}}+|\rho_{\varepsilon_n}|
    \leq\Psi
    \qquad\mu\text{-a.e. for every }n.
\end{equation*}

At $\mu$-almost every point, eventually
$q_{\varepsilon_n}>0$ and
$\norm{\varepsilon_nT_{\varepsilon_n}}<\pi/2$.
Therefore,
\begin{equation*}
    \mathcal L_{\varepsilon_n}(x)
    =
    \left(
        q_{\varepsilon_n}(x)
        \frac{\sin\left(\norm{\varepsilon_n T_{\varepsilon_n}(x)}\right)}
             {\norm{\varepsilon_n T_{\varepsilon_n}(x)}}
        \varepsilon_n T_{\varepsilon_n}(x),
        \,
        2\left(
            q_{\varepsilon_n}(x)
            \cos\left(\norm{\varepsilon_n T_{\varepsilon_n}(x)}\right)-1
        \right)
    \right),
\end{equation*}
with the continuous convention $\sin\theta/\theta=1$
at $\theta=0$, and
\begin{equation*}
    \frac{c_{\varepsilon_n}}{\varepsilon_n^2}
    =
    \frac14|\rho_{\varepsilon_n}|^2
    +
    2q_{\varepsilon_n}
    \frac{1-\cos\left(\norm{\varepsilon_n T_{\varepsilon_n}}\right)}
         {\varepsilon_n^2}.
\end{equation*}
Consequently,
\begin{equation*}
    \frac{\mathcal L_{\varepsilon_n}}{\varepsilon_n}
    \longrightarrow(T,\rho),
    \qquad
    \frac{c_{\varepsilon_n}}{\varepsilon_n^2}
    \longrightarrow\norm{T}^2+\frac14|\rho|^2
    \qquad\mu\text{-a.e.}
\end{equation*}
Together with \eqref{eq:graph-linearization-bound}
and the common majorant $\Psi$,
dominated convergence proves both conclusions along
this subsequence. Since every sequence tending to zero
admits such a subsequence, the full limits follow.
\end{proof}

\section{Radial dilations and gluing}\label{app:radial-dilations-moments}

\subsection{Radial dilations}

We collect the dilation identities used in the cone constructions.
For a finite measure on $\Geo(\cone)$ with finite endpoint radial
second moments, homogeneous time marginals are defined by the same
formula as for probability measures. They are finite, since
$r_t\leq (1-t)r_0+tr_1$ along every cone geodesic.

\begin{lemma}[Dilation and probability normalization]
\label{lemma:cone-dilation-normalization}
Let $A\in\mathcal M_2(\cone^N)$ and let
$\vartheta\colon\cone^N\to(0,\infty)$ be Borel measurable with
$\int\vartheta^2\,dA<\infty$. Define
\begin{equation}\label{eq:appendix-cone-dilation}
    D_\vartheta\bigl((\y_i)_{i=1}^N\bigr)
    \coloneqq\bigl([x_i,r_i/\vartheta]\bigr)_{i=1}^N,
    \qquad
    \operatorname{dil}_{\vartheta,2}A
    \coloneqq(D_\vartheta)\pf(\vartheta^2A),
\end{equation}
where $\y_i=[x_i,r_i]$ and $\vartheta$ is evaluated at the whole
tuple. Then $\operatorname{dil}_{\vartheta,2}A\in\mathcal M_2(\cone^N)$ and
\begin{equation}\label{eq:squared-wfr-dilation-invariance}
    \int\Phi\,d(\operatorname{dil}_{\vartheta,2}A)
    =\int\vartheta^2\Phi(D_\vartheta(\y))\,dA(\y)
    =\int\Phi\,dA
\end{equation}
for every integrable Borel function $\Phi$ that is jointly
$2$-homogeneous in the radial variables. In particular, dilation
preserves homogeneous marginals and integrated squared cone distances.

The normalized cone logarithm is jointly $0$-homogeneous:
\begin{equation}\label{eq:appendix-log-homogeneity}
    \mathcal L_\cone([x,\theta r],[y,\theta s])
    =\mathcal L_\cone([x,r],[y,s]),\qquad\theta>0.
\end{equation}
Consequently, dilation also preserves integrable source-weighted
logarithmic pairings, squared logarithmic errors, and creation masses,
that is, integrals of the forms
\begin{equation*}
\begin{gathered}
    r_i^2\left\langle W(x_i),
           \mathcal L_\cone(\y_i,\y_j)\right\rangle_{\WFR},
    \qquad
    r_i^2\left\|\mathcal L_\cone(\y_i,\y_j)-W(x_i)\right\|_{\WFR}^2,
    \\
    r_j^2\mathbf 1_{\{r_i=0\}},
\end{gathered}
\end{equation*}
where $W\colon\R^d\to\R^d\times\R$ is Borel measurable.
The first two expressions are set to zero when $r_i=0$.
These statements apply on reduced pairs as in
Remark~\ref{rem:cone-log-normalization}, and also on arbitrary cone
pairs when the same normalization of the cone logarithm is used.

Finally, every finite measure $\widetilde{\boldsymbol{\Pi}}$ on
$\Geo(\cone)$ with finite endpoint radial second moments admits
a probability normalization $\boldsymbol{\Pi}$ preserving all
homogeneous time marginals, endpoint cost, and integrable
source-weighted logarithmic pairings. The normalization preserves
concentration on $\Geo^*(\cone)$.
\end{lemma}

\begin{proof}
The dilation map is Borel measurable and
$\operatorname{dil}_{\vartheta,2}A(\cone^N)=\int\vartheta^2\,dA<\infty$.
Joint homogeneity gives
\begin{equation*}
    \int\Phi\,d(\operatorname{dil}_{\vartheta,2}A)
    =\int\vartheta^2\Phi\bigl(D_\vartheta(\y)\bigr)\,dA(\y)
    =\int\Phi\,dA.
\end{equation*}
Applying this to $\sum_i r_i^2$ proves the finite second-moment
assertion. Applying it to $r_i^2\varphi(x_i)$ and
$\dc^2(\y_i,\y_j)$ proves the marginal and cost identities.

For $r>0$, the cone-geodesic formulas give
\begin{equation*}
    \Log_{\cone}([x,\theta r],[y,\theta s])
    =\left(\Log_{\cone}^T([x,r],[y,s]),
           \theta\Log_{\cone}^R([x,r],[y,s])\right).
\end{equation*}
Since $N_{\theta r}(\xi,\theta\zeta)=N_r(\xi,\zeta)$,
this proves \eqref{eq:appendix-log-homogeneity}; the convention at
the source vertex gives the same identity there. The displayed
source-weighted expressions are therefore jointly $2$-homogeneous,
so \eqref{eq:squared-wfr-dilation-invariance} applies.

For path measures, use the same convention
\begin{equation*}
    (D_\vartheta\bfy)_t
    \coloneqq[\bfx_t,\bfr_t/\vartheta(\bfy)].
\end{equation*}
The factor is constant along each path, so this map preserves cone
geodesics and reducedness. Set
$M\coloneqq\widetilde{\boldsymbol{\Pi}}(\Geo(\cone))$.
If $M>0$, define
\begin{equation}\label{eq:appendix-probability-normalization}
    \boldsymbol{\Pi}
    \coloneqq\frac1M(D_{M^{-1/2}})\pf
                  \widetilde{\boldsymbol{\Pi}}.
\end{equation}
This is a probability measure: its radii are multiplied by $\sqrt M$
and its weight is divided by $M$. Every jointly $2$-homogeneous
path observable therefore has the same integral before and after
normalization. Apply this to $r_t^2\varphi(x_t)$, the endpoint cost,
and the source-weighted logarithmic pairings. If $M=0$, use the
probability measure concentrated on the constant vertex path;
all these integrals vanish for both measures.
\end{proof}

\subsection{Gluing prescribed cone couplings}

The next result follows by normalizing the common coordinate and
then applying the usual gluing construction; see
\cite[Lemma~7.10 and the proof of Lemma~7.11]{liero-mielke-savare2018}.
We give the proof to record explicitly the preservation of the
prescribed couplings up to dilation and an atom at the vertex pair.

\begin{lemma}[Gluing prescribed cone couplings]
\label{lemma:squared-wfr-gluing}
Let $\mu_0,\mu_1,\mu_2\in\M$ and
\begin{equation*}
    A_{01}\in\Gamma_{\cone}(\mu_0,\mu_1),
    \qquad
    A_{12}\in\Gamma_{\cone}(\mu_1,\mu_2).
\end{equation*}
There exist positive Borel functions
$\vartheta_{ij}\colon\cone^2\to(0,\infty)$, constants
$c_{ij}\ge0$, and $A\in\mathcal M_2(\cone^3)$ such that,
for $(i,j)\in\{(0,1),(1,2)\}$,
\begin{equation*}
    (\pi_{ij})\pf A
    =
    \operatorname{dil}_{\vartheta_{ij},2}
    \bigl(A_{ij}+c_{ij}\delta_{(\vertex,\vertex)}\bigr),
    \qquad
    (\pi_1)\pf A=\iota\pf\mu_1+\delta_{\vertex},
\end{equation*}
where $\iota(x):=[x,1]$ and $\pi_i,\pi_{ij}$ are the
coordinate projections.

Moreover, for every Borel function $\Phi$ jointly
$2$-homogeneous in the radial variables and integrable
with respect to $A_{ij}$,
\begin{equation*}
    \int_{\cone^3}\Phi(\y_i,\y_j)\,dA
    =
    \int_{\cone^2}\Phi\,dA_{ij}.
\end{equation*}
In particular, the prescribed pair costs and
source-weighted logarithmic pairings are preserved.
\end{lemma}

\begin{proof}
For each $(i,j)\in\{(0,1),(1,2)\}$, let $\rho_{ij}$ denote
the radius of the coordinate corresponding to $\mu_1$.
Set
\begin{equation*}
    z_{ij}\coloneqq A_{ij}(\{\rho_{ij}=0\}),
    \qquad
    c_{ij}\coloneqq\mathbf 1_{\{z_{ij}=0\}},
    \qquad
    M_{ij}\coloneqq z_{ij}+c_{ij}>0,
\end{equation*}
and define
\begin{equation*}
    \vartheta_{ij}\coloneqq
    \begin{cases}
        \rho_{ij},&\rho_{ij}>0,\\
        M_{ij}^{-1/2},&\rho_{ij}=0.
    \end{cases}
\end{equation*}
These are positive Borel functions satisfying
\begin{equation*}
    \int\vartheta_{ij}^2\,
       d\bigl(A_{ij}+c_{ij}\delta_{(\vertex,\vertex)}\bigr)
    =\mu_1(\R^d)+1<\infty.
\end{equation*}
Thus the dilated plans
\begin{equation*}
    B_{ij}\coloneqq\operatorname{dil}_{\vartheta_{ij},2}
       \bigl(A_{ij}+c_{ij}\delta_{(\vertex,\vertex)}\bigr)
\end{equation*}
belong to $\mathcal M_2(\cone^2)$.
On $\{\rho_{ij}>0\}$, the common coordinate has radius one and
the weight is multiplied by $\rho_{ij}^2$. On
$\{\rho_{ij}=0\}$, its ordinary mass becomes
$M_{ij}^{-1}(z_{ij}+c_{ij})=1$. Hence both plans have the same
ordinary marginal in their common coordinate:
\begin{equation*}
    \sigma\coloneqq\iota\pf\mu_1+\delta_{\vertex}.
\end{equation*}

Disintegrate
\begin{equation*}
    B_{01}(d\y_0,d\y_1)
      =\kappa^0_{\y_1}(d\y_0)\,\sigma(d\y_1),
    \qquad
    B_{12}(d\y_1,d\y_2)
      =\sigma(d\y_1)\,\kappa^2_{\y_1}(d\y_2),
\end{equation*}
where the conditional measures are probabilities, and set
\begin{equation*}
    A(d\y_0,d\y_1,d\y_2)
    \coloneqq\kappa^0_{\y_1}(d\y_0)\,
                \sigma(d\y_1)\,
                \kappa^2_{\y_1}(d\y_2).
\end{equation*}
Its prescribed pair marginals are $B_{01}$ and $B_{12}$, and
its middle marginal is $\sigma$. Moreover,
\begin{equation*}
    \int_{\cone^3}(r_0^2+r_1^2+r_2^2)\,dA
    =\mu_0(\R^d)+\mu_1(\R^d)+\mu_2(\R^d)<\infty,
\end{equation*}
so $A\in\mathcal M_2(\cone^3)$.
Finally, a jointly $2$-homogeneous real-valued function satisfies
$\Phi(\vertex,\vertex)=0$. Applying
\eqref{eq:squared-wfr-dilation-invariance} to the augmented pair
plans therefore gives
\begin{equation*}
    \int\Phi(\y_i,\y_j)\,dA
    =\int\Phi\,dB_{ij}
    =\int\Phi\,dA_{ij},
\end{equation*}
as required.
\end{proof}

\section{Technical results}

\subsection{Vanishing-measure estimates}

\begin{lemma}[Vanishing dominated measures]
\label{lemma:vanishing-dominated-measures}
Let $(X,\mathcal B)$ be a measurable space, let $\mu\in\mathcal M^+(X)$ be a finite measure, and let $(\eta_n)_{n\in\N}\subset\mathcal M^+(X)$ satisfy
\begin{equation}
0\leq\eta_n\leq\mu
\qquad\text{for every }n\in\N
\end{equation}
and
\begin{equation}
\eta_n(X)\longrightarrow0.
\end{equation}
Then, for every $f\in L^1(\mu)$,
\begin{equation}
\int_X |f|\,d\eta_n\longrightarrow0.
\end{equation}
\end{lemma}

\begin{proof}
Let $f\in L^1(\mu)$. For every $M>0$, we have
\begin{align*}
\int_X |f|\,d\eta_n
=
\int_{{|f|\leq M}}|f|\,d\eta_n
+
\int_{{|f|>M}}|f|\,d\eta_n
\leq
M\eta_n(X)
+
\int_{{|f|>M}}|f|\,d\mu,
\end{align*}
where we used $\eta_n\leq\mu$ in the second term. Therefore,
\begin{equation*}
    \limsup_{n\to+\infty}
\int_X|f|\,d\eta_n
\leq
\int_{{|f|>M}}|f|\,d\mu.
\end{equation*}
Since $f\in L^1(\mu)$, the right-hand side converges to zero as
$M\to+\infty$. This proves the claim.
\end{proof}

\begin{corollary}
\label{cor:vanishing-dominated-measures-L2}
Under the assumptions of Lemma~\ref{lemma:vanishing-dominated-measures},
for every $g\in L^2(\mu;\R^m)$ one has
\begin{equation}
\int_X |g|^2\,d\eta_n\longrightarrow0.
\end{equation}
\end{corollary}

\begin{corollary}
\label{cor:vanishing-weights}
Let $\mu\in\mathcal M(\R^d)$ and let $(\gamma^n)_{n\in\N}\subset\mathcal M(\R^d\times\R^d)$ satisfy
\begin{equation*}
    (\pi^1)\pf\gamma^n=\mu,\qquad\int_{\R^d\times\R^d}\norm{x-y}^2\,d\gamma^n(x,y)\longrightarrow 0.
\end{equation*}
Let $h:[0,+\infty)\to\R$ be bounded and Borel measurable, with
\begin{equation*}
    h(0)=0,
\qquad
\lim_{r\downarrow0}h(r)=0.
\end{equation*}
Then, for every $g\in L^2(\mu;\R^m)$,
\begin{equation*}
    \int_{\R^d\times\R^d}
\norm{g(x)}^2 |h(\norm{x-y})|^2
\,d\gamma^n(x,y)\longrightarrow 0.
\end{equation*}
\end{corollary}

\begin{proof}
We may assume $\norm{h}_\infty>0$. Define
\begin{equation*}
    \eta_n\coloneqq (\pi^1)_\#
\left(
\frac{|h(\norm{x-y})|^2}{\norm{h}_\infty^2}
\gamma^n
\right).
\end{equation*}
Clearly, $0\le\eta_n\le\mu$. Moreover, for every $\delta>0$,
\begin{equation*}
    \eta_n(\R^d)
\le
\frac{\mu(\R^d)}{\norm{h}_\infty^2}
\sup_{0\le r\le\delta}|h(r)|^2
+
\frac{1}{\delta^2}
\int_{\R^d\times\R^d}\norm{x-y}^2\,d\gamma^n(x,y).
\end{equation*}
Taking first $n\to+\infty$ and then $\delta\downarrow0$,
we obtain $\eta_n(\R^d)\to0$. The conclusion now follows
from Corollary~\ref{cor:vanishing-dominated-measures-L2},
since
\begin{equation*}
    \int_{\R^d\times\R^d}
\norm{g(x)}^2|h(\norm{x-y})|^2\,d\gamma^n(x,y)
=
\norm{h}_\infty^2
\int_{\R^d}\norm{g(x)}^2\,d\eta_n(x).
\end{equation*}
\end{proof}

\subsection{Approximation and tangency}

\begin{lemma}\label{lemma:bounded-support-approximation}
Let $(\xi,\zeta)\in H(\mu)$ and consider for $R>0$ smooth cutoff functions $\chi_R\in C_c^\infty(\R^d)$ satisfying
$0\leq\chi_R\leq1$ and
\begin{equation}\label{eq:cut-off-function}
    \chi_R\equiv1\;\text{on }B_{R/2},
    \qquad
    \operatorname{supp}\chi_R\Subset B_R,
    \qquad
    \norm{\nabla\chi_R}_\infty\leq\frac{C}{R},
\end{equation}
for some $C>0$ independent of $R$. Then $(\chi_R\xi,\chi_R\zeta)\to (\xi,\zeta)$ in $H(\mu)$.
\end{lemma}
\begin{proof}
Since $\chi_R\to1$ pointwise and
$0\leq\chi_R\leq1$, dominated convergence yields
\begin{equation*}
    \norm{(\chi_R-1)\xi}_{L^2(\mu)}\to0,
    \qquad
    \norm{(\chi_R-1)\zeta}_{L^2(\mu)}\to0.
\end{equation*}
Moreover,
\begin{equation*}
    \norm{\nabla\chi_R\cdot \xi}_{L^2(\mu)}
    \leq
    \frac{C}{R}\norm{\xi}_{L^2(\mu)}
    \longrightarrow0.
\end{equation*}
Consequently, $\xi\chi_R\to \xi$ and $\zeta\chi_R\to\zeta$
in $L^2(\mu)$ as desired.
\end{proof}

\begin{lemma}\label{lemma:C1,1-tangency}
Let $\varphi\in C^{1,1}_b(\R^d)$. For every $\mu\in\M$
\begin{equation*}
    (\nabla\varphi, 4\varphi)\in T_\mu\M.
\end{equation*}
\end{lemma}
\begin{proof}
Let $\varphi_R\coloneqq\chi_R\varphi$, with $\chi_R$ as in
Lemma~\ref{lemma:bounded-support-approximation}.
Mollification approximates $\varphi_R$ and its gradient uniformly
by smooth compactly supported functions and their gradients.
Since $\mu$ is finite, this gives
$(\nabla\varphi_R,4\varphi_R)\in T_\mu\M$.

Moreover,
\begin{equation*}
    (\nabla\varphi_R,4\varphi_R)
    =
    (\chi_R\nabla\varphi,\chi_R4\varphi)
    +(\varphi\nabla\chi_R,0)
    \longrightarrow(\nabla\varphi,4\varphi)
    \quad\text{in }H(\mu),
\end{equation*}
by Lemma~\ref{lemma:bounded-support-approximation} and
$\norm{\varphi\nabla\chi_R}_{L^2(\mu)}
\leq C R^{-1}\norm{\varphi}_{L^2(\mu)}$.
The conclusion follows from the closedness of $T_\mu\M$.
\end{proof}

\begin{lemma}\label{lemma:compact-support-approximation}
Let $\mathcal W=(w,\eta)\in T_\mu\M^\perp$. Then there exists
a family of normal vectors
$\mathcal W_R=(w_R,\eta_R)\in T_\mu\M^\perp$, $R>0$,
with compact support contained in
\begin{equation*}
    B_R\coloneqq
    \left\{x\in\R^d\colon \norm{x}<R\right\},
\end{equation*}
such that $\mathcal W_R\to\mathcal W$ in $H(\mu)$
as $R\to+\infty$.
\end{lemma}

\begin{proof}
For $\chi_R$ as in Lemma \ref{lemma:bounded-support-approximation}, define
\begin{equation*}
    w_R\coloneqq\chi_Rw,
    \qquad
    \eta_R\coloneqq
    \chi_R\eta+\nabla\chi_R\cdot w.
\end{equation*}
Both components belong to $L^2(\mu)$ and vanish
$\mu$-a.e.\ outside the compact set
$\operatorname{supp}\chi_R$.

We first verify that $\mathcal W_R\in T_\mu\M^\perp$.
Indeed for every $\varphi\in C_c^\infty(\R^d)$, orthogonality
of $\mathcal W$, tested against $\chi_R\varphi$, gives
\begin{equation*}
\begin{aligned}
    0
    &=
    \int_{\R^d}
    \left(
        \nabla(\chi_R\varphi)\cdot w
        +\chi_R\varphi\eta
    \right)\,d\mu\\
    &=
    \int_{\R^d}
    \left(
        \nabla\varphi\cdot(\chi_Rw)
        +\varphi(\chi_R\eta+\nabla\chi_R\cdot w)
    \right)\,d\mu\\
    &=
    \int_{\R^d}
    \left(
        \nabla\varphi\cdot w_R+\varphi\eta_R
    \right)\,d\mu.
\end{aligned}
\end{equation*}
Convergence follows from Lemma \ref{lemma:bounded-support-approximation}.
\end{proof}

\end{document}